\documentclass{article}
\usepackage{latexsym, amssymb, amsmath, theorem}
\usepackage{epsfig} 
\usepackage{eufrak}
\numberwithin{equation}{section}

\theoremstyle{plain}
\newtheorem{theorem}{Theorem}[section]
\newtheorem{assumption}[theorem]{Assumption}
\newtheorem{proposition}[theorem]{Proposition}
\newtheorem{lemma}[theorem]{Lemma}
\newtheorem{corollary}[theorem]{Corollary}

\newtheorem{definition}[theorem]{Definition}
\newtheorem{example}[theorem]{Example}
\newtheorem{remark}[theorem]{Remark}
\newenvironment{proof}{{\noindent \textbf{Proof}\,\,}}{\hspace*{\fill}$\Box$\medskip}
\def\G{\Gamma}
\def\go#1{\EuFrak#1}   
\def\rr{\mathbb R}
\def\zz{\mathbb Z}
\def\gh{\go h}
\def\gg{\go g}
\def\gt{\go t}

\def\su{\go{su}}

\def\oc{\overline{\mathbb C}}
\def\rgg{R(\G,G)}
\def\xgg{X(\G,G)}

\def\oys{\Omega_{Y,\Sigma}}
\def\qq{\mathbb Q}
\title{Instability of nondiscrete free subgroups in Lie groups}

\author{Alexey Glutsyuk
\thanks{Laboratoire
J.-V.Poncelet (UMI 2615 du CNRS et l'Universit\'e Ind\'ependante
de Moscou).} \thanks{Permanent address: CNRS, Unit\'e de Math\'ematiques
Pures et Appliqu\'ees, M.R., \'Ecole Normale Sup\'erieure de Lyon,
46 all\'ee d'Italie, 69364 Lyon 07, France.  Email:
aglutsyu@umpa.ens-lyon.fr}
\thanks{Supported by part by RFBR grants 02-01-00482, 02-01-22002, 
07-01-00017-a  
and NTsNIL\_a (RFBR-CNRS)  05-01-02801, NTsNIL\_a (RFBR-CNRS) 
10-01-93115.}}

\begin{document}
\maketitle

\begin{abstract} We study finitely-generated nondiscrete free  
subgroups  in Lie groups. We address the following question first raised by 
\'Etienne Ghys: is it always possible to make arbitrarily small 
perturbation of the generators of the free 
subgroup in such a way that the new 
group formed by the perturbed generators be not free? In other words, 
 is it possible to approximate generators of a free subgroup  by elements 
 satisfying a nontrivial relation? We prove that the answer to Ghys' question is 
 positive and generalize this result to certain non-free subgroups. We also consider 
the question on the best approximation rate in terms of the minimal length of relation 
in the approximating group. 
We give an upper bound on the optimal approximation rate as $e^{-cl^{\kappa}}$, where 
$c>0$ is a constant, $l$ the minimal length of relation  
and $0.19<\kappa<0.2$. 
\end{abstract}

\tableofcontents

\def\a{\alpha}
\def\ad{\operatorname{ad}}
\def\Ad{\operatorname{Ad}}

\def\qq{\mathbb Q}
\def\rr{\mathbb R}
\def\re{\operatorname{Re}}
\def\im{\operatorname{Im}}
\def\La{\Lambda}
\def\wt#1{\widetilde#1}
\def\cc{\mathbb C}
\def\xgga{X(\G,Ad(G))}
\section{Introduction and the plan of the paper}

\subsection{Main result: instability of freeness. Plan of the paper}

Let $G$ be a connected non-solvable  Lie group. It is well-known (see \cite{E}) that almost each  
(in the sense of the Haar measure) pair of elements 
$(A,B)\in G\times G$ generates a free noncommutative subgroup in $G$. 
At the same time in the case, when $G$ is connected and semisimple, there is a 
neighborhood 
$U\subset G\times G$ of the identity in $G\times G$ such that 
a topologically-generic pair $(A,B)\in U$ 
generates a dense subgroup: the latter pairs form an open dense subset in $U$. 
This was proved in \cite{BG} (theorem 2.1 on page 450). 

Pairs of elements of $G$ satisfying nontrivial relations form a countable union 
of subvarieties (relation subvarieties) in $G\times G$. We show that the relation 
subvarieties are dense in $U$. We prove a similar result for  $M$- tuples of elements 
in $G$ for arbitrary $M\geq2$. 

\begin{definition} A representation of a group to a Lie group is {\it discrete  (dense)}, 
if its image is respectively discrete (dense).
\end{definition}

The main result of the paper is the following 
\begin{theorem} \label{thth} 
Let $\Gamma$ be a finitely-generated free group of rank at least two, and $G$ be a 
Lie group of positive dimension. Then every  injective nondiscrete representation 
$\rho:\Gamma\to G$ is the limit of a sequence of non-injective representations.
\end{theorem}

{\bf Addendum to Theorem \ref{thth}.} {\it In Theorem \ref{thth} the approximating non-injective representations may be chosen to have non-free images.}

The addendum is proved at the end of Section 6.

\begin{remark} \label{freereps} Under the assumptions of Theorem \ref{thth} let 
$\gamma=(\gamma_1,\dots,\gamma_M)$ 
be a collection of free generators of $\Gamma$. 
The space $\rgg$ of representations $\phi:\G\to G$ is identified with 
$G^M$ via the correspondence $\phi\mapsto \phi(\gamma)=(\phi(\gamma_1),
\dots,\phi(\gamma_M))\in G^M$.  The statement of the theorem says that for every 
$\a\in G^M$ generating a nondiscrete free subgroup in $G$  there exist  a sequence 
$\alpha_r\in G^M$, 
$\alpha_r\to\alpha$, as $r\to\infty$,  and a sequence of 
words $w_r$ in $M$ elements\footnote{Everywhere in the 
paper, by a word in given symbols $\gamma_1,\dots,\gamma_M$ 
we mean a {\it reduced} word (i.e., a word without inverse 
neighbors $\gamma_i^{\pm1}$) in the same symbols 
and their inverses.} for which  
\begin{equation} w_r(\gamma_1,\dots,\gamma_M)\neq1 \text{ in } \Gamma  
\text{ and } w_r(\alpha_r)=1 \text{ for all } r.\label{words}
\end{equation}
Theorem \ref{thth} clearly holds (and is almost obvious) for rank 1 groups $\Gamma$, since 
the closure of a nondiscrete cyclic subgroup in $G$ is a finite extension of $U(1)$. 
\end{remark}

\begin{remark} The condition that the group $\rho(\G)$ is  
nondiscrete is natural. There are discrete 
free subgroups of $PSL_2(\cc)$ (e.g., a Schottky 
group) that are structurally-stable, i.e., remain free under small perturbations.  The 
structurally-stable finitely-generated torsion-free  
subgroups in $PSL_2(\cc)$ were described by D.Sullivan 
(see \cite{sul}, p.243, and Theorem \ref{sull} below).
\end{remark}

We fix a left-invariant Riemannian metric on $G$. It induces a distance on $G$: 
the distance between two points of $G$ is the Riemannian distance, if they lie in the 
same connected component, and infinity otherwise. This induces a distance on 
the product $G^M$.   
We measure the distances between representations as between points of $G^M$ 
(see Remark \ref{freereps}). In what follows for each element $w$ of a free group 
we set 
$$|w|=\text{ the length of the reduced word representing } w.$$

\begin{theorem} \label{cdiap0} In  Theorem \ref{thth} there exist a $c=c(\rho)>0$, 
 sequences of representations $\rho_r:\G\to G$ and elements 
 $w_r\in\G\setminus1$ 
 and a geometric progression 
$l_r=l_r(\rho)\in\mathbb N$ such that for every $r\in\mathbb N$  
$$ l_{r+1}=9l_r, \ |w_r|\leq l_r, \ \rho_r(w_r)=1, \ 
\ dist(\rho_r,\rho)<e^{-(cl_r)^{\kappa}}; \ \kappa=\frac{\ln1.5}{\ln9}.$$
\end{theorem}

\begin{remark} Theorem \ref{cdiap0} does not say that the representations $\rho_r$ 
have non-free images.
\end{remark}

\def\var{\varepsilon}

The question of instability of nondiscrete free subgroups  
was first raised by \'E.Ghys. He also suggested to study 
the best rate of approximation of the generators of a nondiscrete free 
subgroup by generators of subgroups 
having relations of lengths no greater than a given $l$.   
(This is analogous to the approximation of irrational numbers 
by rational ones, where the best approximation rate is well-known; 
it is achieved by continued fractions. In our situation the injective nondiscrete 
 representation 
$\rho$ plays the role of an irrational number, the non-injective ones 
play the role of rationals.)  
Theorem \ref{cdiap0} provides an upper bound of the best approximation 
rate. It is expected that the optimal approximation 
rate should 
be exponential in the minimal length of  relation, while the upper bound 
from Theorem \ref{cdiap0} is subexponential. 

The proof of Theorem \ref{cdiap0} will be given in Section 7. It uses 
Theorem \ref{gapprox} (stated in Subsection 7.1), which deals with a semisimple Lie 
group and a finitely-generated dense subgroup. 
It provides an upper bound for the rate of approximation 
of the elements of the unit ball in the Lie group  by elements of the subgroup 
that  satisfy a bound of derivatives in the parameters of the generators. 
These and related results and open problems are discussed in Subsections 1.2, 
3.2, 7.1 and 7.2. 

Theorem \ref{gapprox} follows (see Subsection 7.1) from Lemma \ref{lwsk} and Theorem
 \ref{gapprox1}, both stated in Subsection 7.1 and proved in Section 8. Theorem \ref{gapprox1} 
 proves the statement of Theorem \ref{gapprox} for every 
 Lie group whose Lie algebra 
 satisfies the so-called weak Solovay-Kitaev inequality (see Definition \ref{weaksk}). 
 This inequality says  that each element $z$ of a Lie algebra 
 $\gg$ equals $[x,y]+[x',y']$, with an estimate on the norms of $x$, $y$, $x'$, $y'$ 
 in terms of the norm of $z$. Lemma \ref{lwsk} proves this inequality for all 
 semisimple Lie algebras. 
 
 In Subsection 7.2 we state a generalization of Theorem \ref{cdiap0}: 
 Theorem \ref{diap}. We deduce Theorem \ref{cdiap0} from Theorems \ref{gapprox} 
 and \ref{diap} at the same place. Theorem \ref{diap} is proved in Section 7. 
 
 Theorem \ref{tsk} (stated in Subsection 7.1 and proved by R.Solovay and A.Kitaev in 
 \cite{sk1, sk2, sk3})  concerns those Lie groups whose Lie 
 algebras  satisfy the (strong) Solovay-Kitaev inequality (see Definition \ref{defsk}). 
 This inequality says that  each element of a Lie algebra is a Lie bracket (with an 
 estimate). 
 For those Lie groups Theorem \ref{tsk} provides an upper bound (stronger than in 
 Theorem \ref{gapprox}) for the rate of approximation of its 
 elements in the unit ball by elements of a finitely-generated dense subgroup. 
 Corollary \ref{cdiap} in Subsection 7.2 is a 
stronger version of Theorem \ref{cdiap0} for those Lie 
groups  whose semisimple parts have Lie algebras satisfying the strong 
Solovay-Kitaev inequality. Corollary \ref{cdiap} follows from Theorems \ref{tsk}, 
\ref{diap}  and Remark \ref{remsk} (whose statement is proved at 
the end of Section 8).  

We show that every complex semisimple Lie algebra and every real semisimple split 
Lie algebra satisfy the 
 strong Solovay-Kitaev inequality (Theorem \ref{thsk} stated in Subsection 7.1). 
 Theorem \ref{thsk} implies Lemma \ref{lwsk}, both theorem and lemma are 
 proved in Subsection 8.1.   

\begin{remark} In the case, when the Lie group under 
consideration is $PSL_2(\mathbb R)$, 
Theorem \ref{thth} easily follows from the density of the elliptic elements 
of finite orders in an open domain of $PSL_2(\mathbb R)$: 
the proof is given in Subsection 3.1. The case of $PSL_2(\mathbb C)$ is the  
first nontrivial case in some sense: the elliptic elements in $PSL_2(\cc)$ 
are nowhere dense. Theorem \ref{thth} for $PSL_2(\cc)$ follows from theorem A in 
\cite{sul} (see Theorem \ref{sull} in Subsection 3.2 below). 
The proof in \cite{sul} is based on a very beautiful idea using 
quasiconformal mappings. 
\end{remark}

In Sections 4 and 5 we prove the next generalization of Theorem \ref{thth}  for 
the following  representations of finitely-generated, but not necessarily free groups to linear algebraic groups. 
 
\begin{theorem} \label{tconj} 
Let $\G$ be a finitely-generated group and $G$ be a real semisimple linear algebraic group. 
Let $\rho:\G\to G$ be an injective  representation. Let the Lie subgroup $\overline{\rho(\G)}$ 
 contain the identity component of $G$  and have no nontrivial connected normal Lie 
subgroups. Let the centralizer of  $\overline{\rho(\G)}$ in the complexification 
$G_{\cc}$ of $G$ coincide with the center of $G_{\cc}$. 
Let the character variety $\xgg$ have dimension at least $dim G$ at 
$[\rho]$. Then $\rho$ is a limit of non-injective representations.
\end{theorem}

A background on representation and character varieties is recalled in Subsection 2.2. 

\begin{remark} \label{ad-irr}
A real semisimple Lie group has no nontrivial connected normal Lie 
subgroups, if and only if its adjoint representation is irreducible. 
\end{remark}
\begin{example}  There are disconnected semisimple but not 
simple linear algebraic groups without nontrivial connected normal Lie 
subgroups. For example, the semidirect product of the group 
$S=SL_2(\rr)\times SL_2(\rr)$ and $\zz_2$, the latter  acting by permuting 
the factors of $S$. 
\end{example}

Theorem \ref{tconj} easily implies Theorem \ref{thth}, as is shown in  Section 6. 
To this end, assuming that the representation 
$\rho:\G\to G$ is dense, we consider the canonical semisimple part $G_{ss}$ of 
$G$ and the induced representation $\rho_{ss}:\G\to G_{ss}$. Both 
$G_{ss}$ and $\rho_{ss}$ are introduced in Definition \ref{sss}. 
We show that the representation $\rho_{ss}$ is almost a product of 
representations from Theorem \ref{tconj}. 
Afterwards Theorem \ref{thth} for the initial group $G$ easily follows by elementary 
algebra. 
 
 In the proof of Theorem \ref{tconj} we use Proposition \ref{open} stated and proved in 
Subsection 2.1. It says that 
the  set of $M$- tuples of elements in a semisimple Lie group generating  dense 
subgroups is open (if non-empty).

The main part of the paper is Section 4, where  we start the proof of Theorem 
\ref{tconj} and prove the Main Technical Lemma. In the same section we complete the 
proof of Theorem \ref{tconj} 
 for groups with proximal elements - a class of linear algebraic 
 groups containing all the groups $SL_n(\rr)$ and more generally, 
all the simple split groups (see \cite{VO}, p.288). 
A reader can read the proofs in Section 4 assuming everywhere that  $G=SL_n(\rr)$. 
The proof of Theorem \ref{tconj} for groups without 
proximal elements will be finished in Section 5. 
One of the key statements used in the proof of Theorem \ref{tconj} is 
Proposition \ref{p} in Subsection 2.4, which concerns infinitesimal deformations of 
representations and extension of the corresponding cocycles. The related 
background material is provided in Subsections 2.3 and 2.4. 

To  construct approximating non-injective representations, we use 
appropriate products of powers of 
iterated commutators and study iterations of appropriate 
commutator mappings $\phi_g:G\to G$, $h\mapsto ghg^{-1}h^{-1}$. We use 
 Theorem \ref{dsplit} and its Corollary \ref{cpoc}, both stated and proved in 
 Subsection 2.10. Theorem \ref{dsplit} is a version of Poincar\'e-Dulac theorem 
(\cite{A}, chapter 5, section 25, subsection D).  It provides a normal 
 form for a class of contracting analytic germs $(\rr^n,0)\to(\rr^n,0)$,  
 which include the commutator mappings under question. 
 
 For the proof of Theorem \ref{tconj} we consider a deformation $\rho_u$ of the 
 given representation $\rho=\rho_0$ depending on a parameter $u\in\rr^n$, 
 $n=dim G$, whose projection to the character variety has rank $n$ at 0 as a function 
 of the parameter. We construct appropriate sequences 
 $k_r\in\mathbb N$ and $\omega_r\in\Gamma$ (products of iterated commutators with 
 growing depths $k_r+const$, $k_r\to+\infty$) so that the mappings 
 $\rr^n\to G$: $\wt u\mapsto
 \rho_{k_r^{-1}\wt u}(\omega_r)$ converge uniformly with derivatives on 
 compact sets to a limit mapping $\Psi$. The main technical part of the 
 proof is to show that one can achieve that the limit $\Psi$ be a local diffeomorphism 
 at 0 (Lemmas \ref{lsi} and \ref{lem1} in Subsection 4.1). The local diffeomorphicity 
 easily implies Theorem \ref{tconj}, as is shown at the end of Subsection 4.1. 
 To prove the local diffeomorphicity in the case, when $G$ has no 
 proximal elements, we show that $G$ contains the so-called $\cc$-1-proximal 
 elements introduced in Subsection 2.7. The derivative at 1 of the commutator 
 mapping associated  to a 
$\cc$-1-proximal element $g$ has a ``dominating'' invariant plane $L(g)\subset\gg$ 
equipped with a canonical structure of complex line. The main contribution to the 
derivative 
$\Psi'(0)$ is given by a linear operator $\Omega_{Y,\Sigma}:T_0\rr^n\to\gg$ introduced 
in Section 5. It depends on appropriate $n$- tuples of $\cc$-1-proximal elements 
$A_j=\rho(g_j)\in G$, $g_j\in\G$,  vectors $Y_j\in L(A_j)$ and 
real-linear complex-valued forms $\sigma_j:T_0\rr^n\to\cc$, $j=1,\dots,n$: 
$\Omega_{Y,\Sigma}=\sum Y_j\sigma_j$. 
In the proof of the invertibility of the operator $\Omega_{Y,\Sigma}$ we use  
preparatory purely linear-algebraic Propositions \ref{lem3} and 
\ref{lem4} stated and proved in Subsection 2.8.
 
A brief historical overview is presented in Subsections 1.2 and 3.2. 

In Subsection 2.9 we give a proof of a version 
of Theorem \ref{thth}  for  the simplest solvable noncommutative 
Lie group $Aff_+(\rr)$ (Proposition \ref{paff}). The author is sure that Proposition \ref{paff} is well known to the specialists. The proof gives an illustration 
of the basic ideas used in the proof of Theorem \ref{thth}.  The analytic argument 
from the proof of Proposition \ref{paff} 
is stated  in a generalized form (Proposition \ref{pann}). 
It will be used in the proof of Theorem \ref{tconj}. 

 The background on roots of semisimple Lie algebras   
 is recalled in Subsection 2.5. 
The background on proximal elements is recalled in Subsection 2.6.

\def\nn{\mathbb N}

\subsection{Historical remarks and some open questions} 

In 1951 M.Kuranishi (\cite{kur}, p.71) established existence of 2- generated free dense 
subgroups in all connected semisimple Lie groups.

 The famous Tits' alternative \cite{T} says that every finitely-generated linear group 
 $G$ satisfies one of the two following incompatible statements:

- either it is virtually-solvable, i.e., contains a solvable subgroup of 
a finite index;

- or it contains a non-abelian free subgroup with two generators. 

The same is true if $G$ is a  linear group over a field of characteristic zero (not necessarily finitely-generated, see the same paper \cite{T}). 
 
Every  dense subgroup of a connected 
semisimple real Lie group satisfies the second 
statement: it contains a free subgroup with two generators. 

The question of possibility to choose the latter free subgroup to be 
dense was  stated in \cite{CG} and studied in \cite{BG} and \cite{CG}. 
\'E.Ghys and Y.Carri\`ere \cite{CG} have proved 
the positive answer in a particular case. E.Breuillard and 
T.Gelander \cite{BG} have done it in the general case. 

T.Gelander \cite{gel} has shown that in every compact nonabelian Lie group 
each finite tuple of elements can be approximated arbitrarily well by another tuple 
(of the same number of elements) that generates a nonvirtually free group. 

D.Sullivan's and F.Labourie's results  
concerning stable subgroups in $SL_2(\cc)$ and $SL_n(\rr)$ are discussed in Subsection 3.2.

A question  concerning Diophantine 
properties of an individual pair  $A,B\in SO(3)$ was studied in \cite{KR}.   
We say that a pair $(A,B)\in SO(3)\times SO(3)$ is {\it Diophantine} 
(see \cite{KR}), if there exists a  $D=D(A,B)>1$ 
such that for every word $w_k(a,b)$ of length $k$ 
$$|w_k(A,B)-1|>D^{-k}.$$
 A.Gamburd, D.Jakobson and P.Sarnak  have stated the following

\medskip

{\bf Question 1 \cite{GJS}.} {\it Is it true that almost every 
pair $(A,B)\in SO(3)\times SO(3)$ is Diophantine?}

\medskip 
V.Kaloshin and I.Rodnianski \cite{KR} proved that almost every pair 
$(A,B)$ satisfies a weaker inequality with the latter right-hand 
side replaced by $D^{-k^2}$. 
\medskip

{\bf Question 2.} {\it Is there an analogue of Theorem \ref{thth} for the group of 

- germs of one-dimensional real diffeomorphisms (at their common fixed point)? 

- germs of one-dimensional conformal diffeomorphisms? 

- diffeomorphisms of compact manifold?}

\medskip

The latter question concerning conformal germs is related to 
study of one-dimensional holomorphic foliations.  A related result was obtained in 
the joint paper \cite{IP} by Yu.S.Ilyashenko and A.S.Pyartli, which deals 
with one-dimensional holomorphic foliations on $\mathbb{CP}^2$ with 
isolated singularities and invariant infinity line. They have shown 
that for a typical foliation the holonomy group at infinity is free. 
Here ``typical'' means ``lying outside a set of zero Lebesgue measure''. 
It is not known whether this is true for an open set of foliations. 
 
 \def\R{\mathcal R}

\section{Preliminaries}

\subsection{Dense subgroups in Lie groups}
Let $\G$ be a free group of rank $M\geq2$, $G$ be a Lie group equipped with 
a Riemannian metric. Recall that $\rgg=G^M$ (Remark \ref{freereps}). For every 
$\alpha\in G^M$, the corresponding representation will be denoted $\rho_{\alpha}$. 
An $M$- tuple corresponding to a dense 
representation will be called {\it irrational}. Set 
$$G^M_{irr}=\{\text{the irrational M- tuples}\}\subset G^M, \ G_0=\text{ the 
identity component of } G.$$

\begin{proposition} \label{open} Let $G$ be a semisimple Lie group. The subset 
$G^M_{irr}\subset G^M$ is open. It is nonempty, if in addition $G$ is connected. 
\end{proposition} 
\begin{proof}  
Given an $\alpha\in G^M_{irr}$, we construct its neighborhood $V\subset G^M$ 
contained in $G^M_{irr}$. 
Recall that there exist a neighborhood $U\subset G_0\times G_0$  of the 
identity and an open and dense subset $U'\subset U$ 
of pairs generating dense subgroups in $G_0$ 
(see the beginning of the paper and the citation \cite{BG} therein). 
 There exist  $w_1, w_2\in\G$ such that 
$(\rho_{\alpha}(w_1),\rho_{\alpha}(w_2))\in U'$ (the density of 
$\rho_{\alpha}$). There exists a neighborhood $V\subset G^M$ 
of  $\alpha$ such that for every 
$\alpha'\in V$ one has $(\rho_{\alpha'}(w_1),\rho_{\alpha'}(w_2))
\in U'$  (continuity). The latter pair generates a dense subgroup in $G_0$, by 
definition. 
Hence, for each $\alpha'\in V$  the representation $\rho_{\alpha'}$ is dense,  
since  each connected component of $G$ 
intersects $\rho_{\alpha'}(\G)$ (as $\rho_{\alpha}(\G)$, by continuity).
  Thus, $V\subset G^M_{irr}$ and $G^M_{irr}$ is open. 
  If $G$ is connected, then $G^M_{irr}\supset U'\times G^{M-2}\neq\emptyset$. 
  The proposition is proved. 
\end{proof}

\begin{remark} One has $G^M_{irr}=\emptyset$, e.g., if $G\slash G_0$ is 
a free group with more than $M$ generators. 
\end{remark}

\def\fd{\mathbb F_2}

\def\wkao{w_k(a(0),b(0))}
\def\wkau{w_k(a(u),b(u))}
\def\wkak{w_k(a(u_k),b(u_k))}
\def\orao{\omega_r(a(0),b(0))}

\def\ojk{\omega_{jk}}
\def\okk{\omega_{kk}}
\def\ok{\omega_k}

\subsection{Representation and character varieties} 
The  most of the background material recalled here may be found in \cite{jm}, pages 
53-61. Let $\G$ be a group with $M$ generators 
$\gamma=(\gamma_1,\dots,\gamma_M)$. Let 
$G\subset GL_N(\rr)$ be a semisimple linear algebraic group. 
The space $R(\G,G)=Hom(\G,G)$ 
is called the {\it representation variety}. The mapping $\rgg\to G^M$: 
$\rho\mapsto\rho(\gamma)$ identifies $\rgg$ with the real 
affine algebraic variety in $G^M$ 
defined by the equations $w(g_1,\dots,g_M)=1$, where $w$ runs over the defining 
 relations of $\G$. We have already seen (Remark \ref{freereps}) 
that $\rgg=G^M$, if $\G$ is free.  
The group $G$ acts on  $\rgg$ by conjugations.  There exist a 
canonical real affine algebraic variety $\xgg$ called the {\it character variety} and 
a polynomial map $\pi:\rgg\to\xgg$ 
that is constant on the conjugacy classes, see \cite{jm}, p.53 for details. Namely, there 
 exists a finite collection of real polynomials 
$f_1,\dots,f_l$ on $\rgg$ that generate the algebra of $G$- invariants. Set 
$$\pi:\rgg\to\rr^l: \ (g_1,\dots,g_M)\mapsto(f_1(g_1,\dots,g_M),\dots,
f_l(g_1,\dots,g_M)), \ [\rho]=\pi(\rho),$$
$$\xgg=\text{ the Zariski closure of } \pi(\rgg).$$
The projection $\pi$ induces a continuous mapping $\rgg\slash G\to\xgg$ that is 
not necessarily injective (\cite{jm}, p.53) and not necessarily surjective. Both $\rgg$ and 
$\xgg$ are affine algebraic varieties, thus, 
each of them contains a dense Zariski open subset of smooth points, denoted 
respectively 
$R_{reg}(\G,G)$ and $X_{reg}(\G,G)$. Recall that the {\it dimension} of an analytic 
set $X$ at its (maybe singular) point $y$, denoted $dim_yX$, is the maximum of the 
dimensions of the irreducible components of $X$ passing through $y$. Let $S\subset\rgg$ denote 
the subset of those representations $\rho$ for which $\overline{\rho(\G)}\supset G_0$ and the 
centralizer of the Lie subgroup $\overline{\rho(\G)}$ in the complexification $G_{\cc}$ of $G$ 
coincides with the center of $G_{\cc}$. Set $S_{reg}=S\cap R_{reg}(\G,G)$. 
The subsets $S,S_{reg}\subset\rgg$ are open (Proposition \ref{open} applied to the Lie group 
$\overline{\rho(\G)}$ for every $\rho\in S$),  and $S_{reg}$ is dense in $S$. 
\begin{proposition}  \label{densecenter} Let $G$ be a real semisimple linear algebraic 
group, $S\subset\rgg$ be the above open subset, and let $S$ be non-empty. 
Then the projection $\pi:S\to\xgg$ is open,  
$\pi(S_{reg})\subset X_{reg}(\G,G)$, and $\pi:S_{reg}\to X_{reg}(\G,G)$ is a 
submersion. For every $\rho\in S$ one has 
$dim_{[\rho]}\xgg=dim_{\rho}\rgg-dim G$.
\end{proposition}
\begin{proof}  
The stabilizer  of each $\rho\in S$ in $G_{\cc}$  coincides with the center of $G_{\cc}$, as does 
the centralizer of $\overline{\rho(\G)}$. 
The $G_{\cc}$- orbit of every $\rho\in S$  is closed, since $\overline{\rho(G)}\supset G_0$ and 
by   theorem 1.1 on p.54 in \cite{jm}. Thus,  $\rho$ is good in the sense 
of the definition in \cite{jm}, p.53.  Therefore, the action of $G$ on $\rgg$ admits a local  
analytic slice through every $\rho\in S$ and  the projection $\pi$ induces 
an analytic equivalence between the slice  and 
a neighborhood of $[\rho]$ in $\xgg$ (\cite{jm}, theorem 1.2 and the remark on p.57). 
This implies the proposition. 
\end{proof}
\begin{example} \label{exdense} Let $G$ be the automorphism group of a semisimple Lie algebra. 
The centralizer of $G_0$  both in $G$ and $G_{\cc}=Aut(\gg_{\cc})$ is trivial. 
Therefore, every representation $\rho:\G\to G$ such that 
$\overline{\rho(\G)}\supset G_0$ satisfies the statements of Proposition 
\ref{densecenter}.
\end{example}

\subsection{Infinitesimal deformations}
Let $G$ be a real Lie group. 
Everywhere in the paper for every tangent vector $v\in T_gG$ and any $h\in G$ we 
denote
\begin{equation}
vh\in T_{gh}G \text{ the image of } v \text{ under the right multiplication by } h.
\label{vhg}\end{equation}

Let $\G$ be a group with $M$ generators, $\rho\in\rgg$. 
Consider all the germs $\rho_t$ of 1-parametric  deformations of 
$\rho$ in $\rgg$:
$$(\rr,0)\to(\rgg,\rho), \ t\mapsto\rho_t\in\rgg, \  \rho_0=\rho,$$
$C^1$- smooth as mappings to the ambient manifold $G^M\supset\rgg$. 
Their derivatives $\frac{d\rho_t}{dt}|_{t=0}\in T_{\rho}G^M$ form the 
{\it tangent cone} to $\rgg$ at $\rho$. (This definition of the tangent cone 
is equivalent to that from \cite{jm}, page 60, which deals with analytic deformations.) 
The tangent cone will be denoted $TC_\rho$. This is a tangent space to $\rgg$, if 
$\rgg$ is smooth at $\rho$. Each element of $TC_{\rho}$ is given by a family  of 
tangent vectors $\{ v_a\in T_{\rho(a)}G\}_{a\in\G}$, 
$v_a=\frac{d\rho_t(a)}{dt}|_{t=0}$, satisfying the {\it Leibniz multiplication 
rule:} $v_{ab}$ is the vector of the infinitesimal movement of $\rho(ab)$, while 
$\rho(a)$ moves by $v_a$ and $\rho(b)$ moves by $v_b$. Consider the vector function  
\begin{equation}
c:\G\to\gg;\  c(a)=v_a(\rho(a))^{-1}\in\gg,\label{corcoc}\end{equation}
see (\ref{vhg}). The Leibniz multiplication rule is equivalent to the statement that 
$c(a)$ is a 1-cocycle with coefficients in $Ad\circ\rho$,  see the following definition. 
We call it the {\it cocycle 
associated} to the element  $\frac{d\rho_t}{dt}|_{t=0}\in TC_{\rho}$.  

\begin{definition} (\cite{w}, p.150; \cite{jm}, p.59) Let $\G$ be a group, $V$ be 
a vector space, $\phi:\G\to GL(V)$ be a linear representation. 
A {\it 1- cocycle} with coefficients in $\phi$ is a vector function $c:\G\to V$ satisfying 
the identity
\begin{equation} c(ab)=c(a)+\phi(a)c(b) \text{ for every } a,b\in\G,\label{coc}\end{equation}
in particular, $c(1)=0$. Identity (\ref{coc}) is called the {\it cocycle identity}. 
A 1-cocycle is called a {\it 1- coboundary}, if there exists a vector $v\in V$ such that 
\begin{equation}c(a)=v-\phi(a)v \text{ for every } a\in\G.\label{cobound}\end{equation}
The space of 1-cocycles (1- coboundaries) with coefficients in $\phi$ 
will be denoted $Z^1(\G,\phi)$ (respectively, $B^1(\G,\phi)$). The {\it first cohomology 
group} is $H^1(\G,\phi)=Z^1(\G,\phi)\slash B^1(\G,\phi)$. 
\end{definition}

\def\coc{Z^1(\G,Ad\circ\rho)}

\begin{remark} \label{cohom} Let $\G$ be a finitely-generated group, 
$G$ be a Lie group,  $\rho\in\rgg$. The correspondence 
$\frac{d\rho_t}{dt}|_{t=0}\mapsto c$, see formula (\ref{corcoc}), identifies the tangent 
cone $TC_{\rho}$ with a subset 
in $\coc$. The tangent space to the $G$- orbit of $\rho$ is contained in 
$TC_{\rho}$ and identified with $B^1(\G,Ad\circ\rho)$ (\cite{w}, p.151;  \cite {jm}, p.61): 
each coboundary (\ref{cobound}) with $\phi=Ad\circ\rho$, $v\in\gg$, 
corresponds to the derivative 
$\frac{d\rho_t}{dt}|_{t=0}$, where $\rho_t=\exp(tv)\rho\exp(-tv)$, and vice versa. 
If $G$ is a semisimple linear algebraic group and $\rho\in S_{reg}$, see Proposition 
\ref{densecenter},  
then $T_{[\rho]}\xgg\subset H^1(\G,Ad\circ\rho)$, by the latter statement and Proposition 
\ref{densecenter}. If in the above assumptions $\G$ is free, then each 
cocycle in $\coc$ is associated to the derivative of some $\rho_t$ and 
$T_{\rho}\rgg=\coc$,  $T_{[\rho]}\xgg=H^1(\G,Ad\circ\rho)$. But in general, this is not 
the case (\cite{jm}, p.53).
\end{remark}
We  will use the following immediate corollaries of the cocycle identity 
(\ref{coc}): 
\begin{equation} c(a_1\dots a_l)=c(a_1)+\phi(a_1)c(a_2)+\dots+\phi(a_1\dots a_{l-1})
c(a_l) \text{ for all } l\in\nn \text{ and } a_j\in\G,\label{cocn}\end{equation}
\begin{equation} c(g^{-1})=-\phi(g^{-1})c(g) \text{ for every } g\in\G.\label{cocin}
\end{equation}

%In fact, 
%the space $\coc$ is always the Zariski tangent space to $\rgg$ at $\rho$. If $G$ is a 
%semisimple Lie group, and $\rho:\G \to G$ is an injective representation with dense 
%image, then $H^1(\G,Ad\circ\rho)$ is the Zariski tangent space at $[\rho]$ 
%to the character variety $\xgg$ (utochnit utverzhdenie i ssylku???)

\subsection{Extension of cocycles}
\begin{theorem} \label{h=0} (\cite{gui}, corollary 7.9, ch.III, p.240). 
Let $G$ be a semisimple Lie group, $\phi$ be its finite-dimensional 
representation. Then every 
continuous cocycle in $Z^1(G,\phi)$ is a coboundary.
\end{theorem}

\begin{remark}
We will use only the particular case of Theorem \ref{h=0} for 
{\it locally-Lipschitz} cocycles and $\phi=Ad$. A locally-Lipschitz cocycle with 
coefficients in $Ad$ corresponds to a locally-Lipschitz vector field via the 
construction from the beginning of the previous subsection. In this case the theorem 
can be proved directly by showing that the flow of the latter field 
(which is locally well-defined by Lipschitz property) is a one-parametric 
family of inner automorphisms. This can be deduced from the cocycle identity 
and semisimplicity. 
\end{remark}

One of the key arguments in the paper is the following

\begin{proposition} \label{p} Let $G$ be a Lie group, $\G\subset G$ be a dense 
subgroup. Let $\phi:G\to GL_N(\rr)$ be a continuous 
linear representation of $G$,  
$\zeta\in Z^1(\G,\phi)$ be a cocycle that is Lipschitz at the identity (with respect to 
the metric on $\G$ induced from $G$). Then $\zeta$ extends to a 
locally-Lipschitz cocycle on $G$.
\end{proposition}

\begin{corollary} \label{corp} Under the assumptions of Proposition \ref{p} let 
the group $G$ be semisimple. Then the cocycle $\zeta$ is a coboundary.
\end{corollary}

The corollary follows immediately from Proposition \ref{p} and Theorem \ref{h=0}. 

\begin{proof} {\bf of Proposition \ref{p}.} Fix a left-invariant Riemannian metric on $G$. 
Recall that $\zeta(1)=0$ and $\zeta$ is Lipschitz at 1: there  exist $L,r>0$ 
(let us fix them) such that  
\begin{equation} |\zeta(\gamma)|\leq Ldist(\gamma,1) \text{ for every } 
\gamma\in D_r\cap\G.
\label{zetast}\end{equation}
{\bf Claim.} {\it For each compact set $K\subset G$ there exists a constant $L_K>0$ such 
that for every $a\in K\cap\G$ and $b\in aD_r\cap\G$ one has} 
\begin{equation}|\zeta(b)-\zeta(a)|\leq L_Kdist(a,b).\label{zg}\end{equation}

\begin{proof} Take arbitrary 
$a\in K\cap\G$, $b\in aD_r\cap\G=a(D_r\cap\G)$ and set $\gamma=a^{-1}b$. One has 
$$\gamma\in D_r\cap\G, \ b=a\gamma, \ \zeta(b)-\zeta(a)=\phi(a)\zeta(\gamma),$$
by the cocycle identity. This together with (\ref{zetast}) and the equality 
$dist(a,b)=dist(\gamma,1)$ (the 
left invariance of the metric)  implies (\ref{zg}) with 
$L_K=L\max_{g\in K}||\phi(g)||$. 
\end{proof}

Let $K\subset G$ be an arbitrary compact subset. 
For every $a,b\in K\cap\G$ with 
$dist(a,b)<r$ one has $b\in aD_r\cap\G$, and inequality (\ref{zg}) holds. This implies 
that the cocycle $\zeta$ is Lipschitz on $K\cap\G$. Hence, $\zeta$ 
extends to $K=\overline{K\cap\G}$ (passing to limits)  as a Lipschitz cocycle with coefficients in $\phi$. This proves Proposition \ref{p}.
\end{proof}

\subsection{Roots of semisimple Lie algebras} 
Everywhere the complexification of a real vector space (Lie algebra) 
$\gg$ will be denoted $\gg_{\cc}$. By the {\it multiplicity} of an eigenvalue of a linear 
operator we understand its {\it algebraic} multiplicity as that of complex zero of 
the characteristic polynomial. Recall \cite{VO} that an element of a semisimple Lie 
algebra is called {\it regular}, if its adjoint  has the minimal possible multiplicity of zero 
eigenvalue.

 Let $\gg$ be a complex semisimple Lie algebra, $x\in\gg$ be a regular element, 
 $\gh\subset\gg$ be the centralizer of $x$; then $x\in\gh$. 
 Recall (see, \cite{VO},  pp. 153, 159) that

- a) $\gh$ is a maximal commutative subalgebra called the {\it  Cartan subalgebra} 
associated to $x$; 

- b) the adjoint action of $\gh$ on $\gg$ is diagonalizable in an appropriate basis of 
$\gg$; 

- c) the corresponding eigenvalues  are linear functionals on $\gh$, 
thus, elements of $\gh^*$; the nonidentically zero ones are called {\it roots}, the 
root eigenspaces are complex eigenlines;

- d) for every root $\a$ the only roots complex-proportional to $\a$ are $\pm\a$; 

- e) some roots form a complex basis in $\gh^*$ and moreover, 
an integer  root basis in the following sense: each root is 
an integer linear combination of the basic roots; 

Statement d) follows from the analogous statement in \cite{VO} (theorem 6 on p.159) 
for real-proportional roots and from statement e). 

If $\gg$ is a real semisimple Lie algebra, $x\in\gg$ is a regular element, then 
similarly, the centralizer of $x$ is called the Cartan subalgebra. Its complexification 
is the above-defined complex Cartan subalgebra in $\gg_{\cc}$.

\begin{proposition} \label{ad+-} Let $G$ be a real connected semisimple Lie 
group. For every 
$x\in\gg$ ($g\in G$) and every complex eigenvalue $\lambda$ of $\ad_x$ ($Ad_g$) the number 
$-\lambda$ (respectively, $\lambda^{-1}$) is also an eigenvalue of the corresponding 
adjoint with the same multiplicity, as $\lambda$. 
\end{proposition} 
\begin{proof} It suffices to prove the statement of the proposition for the Lie 
algebra: this would imply its statement for every $g\in G$ 
covered by an exponential chart, and then, for every $g\in G$ 
(connectedness and analyticity). For a regular $x\in\gg$ the nonzero 
eigenvalues of $\ad_x$ are the values at $x$ of the roots of the corresponding 
complex Cartan  subalgebra in $\gg_{\cc}$. The latter  
are split into pairs of opposite values with equal multiplicities 
(statement d)). The same holds for every $x\in\gg$, by density of regular elements. 
This proves the proposition.
\end{proof}

\subsection{Proximal elements}

\begin{definition} \label{ddsplit} 
 A finite-dimensional (real or complex) linear   
operator $E\to E$ has a {\it dominated splitting}, if its complex eigenvalues are split into 
two non-empty collections, $\Lambda_+$ and $\Lambda_-$, such that 
$$|\lambda_+|>|\lambda_-| \text{ for every } \lambda_{\pm}\in\Lambda_{\pm}.$$ 
(In particular,  $0\notin\Lambda_+$.)
The dominated splitting is the decomposition $E=E_+\oplus E_-$ into the 
direct sum of the invariant subspaces corresponding to $\Lambda_{\pm}$. The subspace 
$E_+$ is called the {\it dominating subspace}.  
\end{definition}

\begin{definition} \label{defprox} 
A linear operator $\rr^n\to\rr^n$ is called {\it proximal}, if it has a dominated splitting with 
one-dimensional dominating space. A Lie group element is proximal, if its 
adjoint is.
\end{definition}
\begin{remark} \label{ropen} The complex eigenvalue of maximal modulus of a proximal 
operator is always unique, simple, real and nonzero.  The set of proximal operators (elements) is open. A positive power of a proximal operator (element) is proximal. This implies that if 
 an algebraic group contains a proximal element, then so does its identity component.
\end{remark}

\begin{definition} A {\it maximal connected $\rr$- split torus} in  a real  semisimple linear 
algebraic group 
 $G$ is a maximal Lie subgroup isomorphic to $(\rr^*_+)^k$, $k\in\nn$, with a 
 diagonalizable adjoint action on $\gg$. (Or equivalently, this is the idendity component 
 of the maximal $\rr$- split algebraic torus, see \cite{VO}, pp.123, 288, 292.) The group 
 $G$  is called {\it split} (see \cite{VO}, p. 288),
 if some its maximal 
connected $\rr$- split torus is a maximal connected commutative Lie subgroup. 
\end{definition}

\begin{example} Each group $SL_n(\rr)$ is split and contains proximal elements: 
 the diagonal matrices with positive eigenvalues form a 
maximal connected $\rr$- split torus; a typical diagonal matrix is a proximal element. 
In general, each maximal connected $\rr$- split torus $T$ of a split group $G$ contains proximal elements. Indeed, let $\go t\subset\gg$ denote the Lie algebra of $T$. Consider 
the roots corresponding to the Cartan subalgebra $\go t_{\cc}\subset\gg_{\cc}$.  
The root values at a generic $v\in\go t$ are distinct and nonzero. 
Fix this $v$. The exponent of the maximal root value is the biggest eigenvalue of $Ad_{\exp v}$. It is simple, thus, $\exp v\in G$ is proximal. 
\end{example}

\begin{example} The group $SO(3)$ is compact and hence, 
not split, has trivial maximal connected $\rr$- split torus and no proximal elements. 
 The simplest example of a non-split noncompact simple group is $SO(3,1)$, 
 which has no proximal elements. The group $SU(2,2)$ is simple, non-split 
 and has proximal elements.  
\end{example}

\begin{theorem} \label{lab} Let a semisimple linear algebraic group contain a proximal 
element. Then each its maximal connected $\rr$- split torus contains a proximal element.
\end{theorem}
The proof of Theorem \ref{lab} is implicitly contained in \cite{abels} 
(p.25, proof of theorem 6.3). It was shown in loc. cit. that under the assumptions of 
Theorem \ref{lab} the maximal $\rr$- split algebraic torus contains proximal elements. 
Thus, so does its identity component (the last statement of Remark \ref{ropen}). 

\begin{definition} \label{1pr} An element $g$ of a Lie group $G$ will be called 
{\it 1- proximal,} if the operator $Ad_g-Id$ is proximal. Let $s(g)$ denote its eigenvalue 
with maximal modulus. Set
\begin{equation}\Pi_{\rr,1}=\Pi_{\rr,1}(G)=\{ \text{1- proximal elements} \ g\in G \ | \ |s(g)|<1\}.
\label{prox<1}\end{equation} 
\end{definition}
We use the following characterization of semisimple algebraic groups with 
proximal 
elements.

\begin{theorem} \label{1prpr} A semisimple linear algebraic group $G$ contains a 
proximal element, 
if and only if  $\Pi_{\rr,1}(G_0)\neq\emptyset$. In this case $\Pi_{\rr,1}(G)\subset G$ is an open 
subset accumulating  to the identity.
\end{theorem}
In the proof of this theorem given below we use the following properties of the adjoint 
representation of a semisimple Lie group.

\begin{proposition} \label{pr1pr} In a connected semisimple Lie 
group every 1- proximal element $g$ with $s(g)>-1$ is proximal.
\end{proposition}
\begin{proof} The eigenvalues  $(s(g)+1)^{\pm1}$ of $Ad_g$ are simple   
(Remark \ref{ropen} and Proposition \ref{ad+-}). 
We claim that $(s(g)+1)^{\pm1}$ is the unique eigenvalue of $Ad_g$ with maximal 
modulus, if $s(g)\in\rr_{\pm}$, $s(g)>-1$. Indeed, in both cases 
$(s(g)+1)^{\pm1}\geq|s(g)|+1$. Case ``$+$'' is obvious. Case ``$-$'': 
$-1<s(g)<0$, $(s(g)+1)^{-1}=(1-|s(g)|)^{-1}>1+|s(g)|$. 
For each eigenvalue $\lambda\neq s(g)$ of 
$Ad_g-Id$ one has $|\lambda|<|s(g)|$ (1- proximality). Hence, 
$(s(g)+1)^{\pm1}\geq|s(g)|+1>|\lambda|+1\geq|\lambda+1|.$ 
This proves the proximality of $Ad_g$ and Proposition \ref{pr1pr}.
\end{proof}

\begin{proposition} \label{prox1pr} Let $G$ be a semisimple linear algebraic group, 
$T\subset G$ be a maximal connected 
$\rr$- split torus. Let $g\in T$ be a proximal element of $G$. 
Then $g$ is also 1- proximal.
\end{proposition}
\begin{proof} The eigenvalues of $Ad_g$ are real, since $Ad_T:\gg\to\gg$ 
is diagonalizable. They  are positive, since this is true for $Ad_1=Id$ 
and the torus $T$ is connected. The nonunit eigenvalues are split into pairs of inverses 
(Proposition \ref{ad+-}). Hence, we can order them as follows 
(distinct indices correspond to distinct (may be multiple) eigenvalues):
\begin{equation}0<\lambda_1^{-1}<\lambda_2^{-1}<\dots<\lambda_k^{-1}<1<
\lambda_k<\dots<\lambda_1.\label{lambdy}\end{equation}
The eigenvalue $\lambda_1$ is simple (proximality). One has 
$$\lambda_1-1>\lambda_1^{-1}(\lambda_1-1)=1-\lambda_1^{-1}, \ 
\text{since} \ 0<\lambda_1^{-1}<1,$$
by (\ref{lambdy}). This together with (\ref{lambdy}) implies that $\lambda_1-1$ 
is a simple eigenvalue of $Ad_g-Id$ with maximal modulus. Hence, the operator 
$Ad_g-Id$ is proximal. Proposition \ref{prox1pr} is proved.
\end{proof}

\begin{proof} {\bf of Theorem \ref{1prpr}.} If $\Pi_{\rr,1}(G_0)\neq\emptyset$, then 
each element of $\Pi_{\rr,1}(G_0)$ is proximal  (Proposition \ref{pr1pr}). 
Conversely, let $G$ contain proximal elements. Let $T\subset G$ be a 
maximal connected 
$\rr$- split torus, $g\in T$ be a proximal element of $G_0$ (which exists by Theorem  
\ref{lab}). Consider the 1- parametric subgroup $H\subset T$ passing 
through $g$. Recall that the Lie groups under consideration are real. 
The powers $g^a$, $a\in\rr_+$,  are well-defined elements of $H$. They 
are proximal (Remark \ref{ropen}) and hence, 1- proximal (Proposition \ref{prox1pr}). 
One has $g^a\to1$, as $a\to0$. Therefore, 
$g^a\in\Pi_{\rr,1}(G_0)$, whenever $a$ is small enough. Thus, the set 
$\Pi_{\rr,1}(G)$ is nonempty, open (Remark \ref{ropen}) and accumulates to 1. 
This proves the theorem. 
\end{proof}

\def\la{\lambda}

\subsection{$\cc$-1-proximal elements}

We will show that the identity component of 
every semisimple linear algebraic group without proximal elements 
contains the so-called $\cc$-1-proximal elements introduced below. 

\begin{proposition} \label{p3.1} Let $G$ be the identity component of a 
semisimple linear algebraic  group. There 
exists a nonempty  subset $U\subset G$ open in the induced Zariski topology on $G$ 
such that the adjoint of each 
$g\in U$ satisfies the following statements:

1) the number of its nonunit complex eigenvalues is maximal, nonzero, and 
they are simple;

2) if there are distinct eigenvalues $\Lambda_1,\Lambda_2\neq1$ with 
$|\Lambda_1-1|=|\Lambda_2-1|$, then $\Lambda_1=\overline\Lambda_2$. 
\end{proposition}
The proposition is proved below. 

\begin{definition} \label{defc1} A linear operator $\rr^n\to\rr^n$ 
is called {\it $\cc$- proximal}, 
if  it has a dominated 
splitting (see Definition \ref{ddsplit}) with dominating subspace of dimension two that 
corresponds to a pair of non-real complex-conjugate eigenvalues.
A Lie group element $g$ is {\it $\cc$-1-proximal}, if 
the operator $Ad_g-Id:\gg\to\gg$ is $\cc$- proximal.  
\end{definition}

\begin{remark} \label{r1prox} A linear operator $\rr^n\to\rr^n$ is $\cc$- proximal, if and 
only if it has exactly two complex 
eigenvalues of maximal modulus, and they are nonreal, hence complex-conjugate,   
nonzero and simple.  The set of $\cc$-1-proximal elements is open (it may be empty).
\end{remark}

\begin{proposition} \label{p3.2} Every element of a semisimple Lie group 
 whose adjoint 
satisfies the above statements 1) and 2) in Proposition \ref{p3.1} 
is either 1- proximal (see Definition \ref{1pr}),  
or $\cc$-1-proximal.
\end{proposition}
\begin{proof} Let $Ad_g$ satisfy 1) and 2), $\lambda$ be its eigenvalue for which the 
modulus $|\lambda-1|$ is the maximal possible. Then $\lambda-1\neq0$ and $\lambda$ is 
a simple eigenvalue (statement 1)). For every eigenvalue $\lambda'\neq\lambda,
\bar\lambda$ one has $|\lambda-1|>|\lambda'-1|$ (statement 2)). Therefore, 
$g$ is 1- proximal, if $\lambda\in\rr$, and $\cc$-1-proximal otherwise. Proposition 
\ref{p3.2} is proved.
\end{proof}
 
 We will use the following properties of the adjoint of a $\cc$-1-proximal element.
 \begin{proposition} \label{p3.3} 
  Let $A:\rr^n\to\rr^n$ be a linear operator with a pair of simple 
 complex-conjugate eigenvalues $s,\bar s\notin\rr$. Let $L\subset\rr^n$ be the 
 invariant two-dimensional subspace corresponding to them. The plane $L$ 
 carries exactly two $A$- invariant structures of complex line 
 compatible with its $\rr$- linear structure. These complex structures 
 are  conjugate. The 
 restriction $A:L\to L$ acts by the multiplication by $s$ (respectively, $\bar s$) 
 in these complex structures. 
  \end{proposition}
 \begin{proof} The restriction 
 $A|_L$ is  conjugate to the operator $S:\cc\to\cc$ acting by multiplication by 
 $s$ (basic linear algebra). The only $S$- invariant complex structures on $\cc$ 
  are the standard one and its complex conjugate (where $S$ acts by multiplication by 
 $\bar s$).  Indeed, each complex structure  
defines an ellipse centered at 0 (up to homothety), which is an orbit of a vector under the multiplication by the complex numbers with unit modulus. 
A complex structure different from the two above ones defines an ellipse that is 
not a circle. The latter ellipse is not homothetic to its image under $S$, since 
$s\notin\rr$. Hence, the complex structure under question cannot be $S$- invariant. 
This proves the proposition. 
\end{proof}

 \begin{definition} \label{defstr}
  Let $G$ be a Lie group, $g\in G$ be a $\cc$-1-proximal element, $s(g)$ be an 
  eigenvalue of $Ad_g-Id$ with the maximal modulus, 
 $L(g)\subset\gg$ be the dominating invariant plane (which corresponds to the 
 eigenvalues $s(g)$ and $\bar s(g)$, see Definition \ref{defc1}). The 
 corresponding complex structure on 
 $L(g)$, in which $Ad_g-Id:L(g)\to L(g)$ acts by multiplication by $s(g)$, will be called 
 the {\it $s(g)$- complex structure.}
 \end{definition}  
 
 \def\osg{\overline{s(g)}}
 
\begin{proposition} \label{p3.4} Let $G$ be a Lie group, $V\subset G$ be a connected 
component of the subset of the $\cc$-1-proximal elements. 
The values $s(g)$, $\bar s(g)$ from Definition \ref{defstr} yield two 
real-analytic nonwhere vanishing complex-conjugate functions $s,\bar s:V\to\cc=\rr^2$. 
Moreover, $\arg s$ and $\arg \bar s$ are single-valued on $V$ as well. 
\end{proposition}
\begin{proof} The value $s(g)$ does not vanish for each $g\in V$ 
(Remark \ref{r1prox}). The local real analyticity of the eigenvalues $s(g)$ and 
$\bar s(g)$ follows from their simplicity. The global analyticity (say, of $s(g)$) follows from the fact that its 
analytic continuation along every closed loop in $V$ does not change the analytic 
branch. Indeed, the result of analytic continuation 
of $s(g)$ remains an eigenvalue of 
$Ad_g-Id$ with the maximal modulus, by definition and local 
analyticity. Therefore, given a $g_0\in V$ and a loop $\gamma\subset V$ 
based at $g_0$, 
the result of the analytic continuation of $s(g)$ along $\gamma$ is either $s(g_0)$, 
or $\bar s(g_0)$. In the latter case there exists a $g'\in\gamma$ where 
$s(g')\in\rr$, by continuity. Thus, $s(g')\in\rr$ is a complex 
eigenvalue of $Ad_{g'}-Id$ with maximal modulus, - a contradiction 
to the $\cc$-1- proximality. The above argument shows that the value $s(g)$ 
corresponding to a $\cc$-1-proximal element $g$ cannot cross the real axis. Hence, 
$\arg s(g)$ and $\arg\bar s(g)$ are also single-valued analytic functions. 
Proposition \ref{p3.4} is proved.
\end{proof}
\def\pc1{\Pi_{\cc,1}}

Everywhere below  we fix a real-analytic branch of the 
eigenvalue function $s(g)$ from  Proposition \ref{p3.4}, defined on the open set of 
all the $\cc$-1-proximal elements. The corresponding  families 
of dominating planes $L(g)\subset\gg$ and the $s(g)$- complex structures 
on them (see the above definition) also depend analytically on $g$. 
We define the multiplication of vectors in $L(g)$ 
by complex numbers in the sense of the $s(g)$- complex structure. Thus, 

- a product of a vector in $L(g)$ and a complex number is again a vector in the plane 
$L(g)$;

- a product of a vector $\nu\in L(g)$ and a $\rr$- linear 
complex-valued 1- form $\sigma:\rr^n \to\cc$ 
is a linear operator $\nu\sigma:\rr^n\to L(g)$.

Set 
\begin{equation} \pc1=\{ \cc-1-\text{proximal elements} \ g\in G \ \text{with} \ 
|s(g)|<1\}.\label{pc1}\end{equation}

\begin{proposition} \label{c1pre}
Let $G$ be a semisimple linear algebraic group without proximal elements. 
The set $\pc1$ is open. It is dense in a neighborhood of the identity.
\end{proposition}
\begin{proof} Consider the following neighborhood of the identity in $G_0$:
$$W=\{ g\in G_0 \ | \ \text{ each complex eigenvalue of } Ad_g-Id \text{ lies in the 
unit disk}\}.$$
The set $W$ contains no 1- proximal element, since otherwise, the latter would 
be proximal (Proposition \ref{pr1pr}), - a contradiction to the no proximal element 
condition. The subset $U\subset G_0$ from Proposition \ref{p3.1} is 
open and dense. The intersection $W\cap U$ is dense in $W$ and consists 
of $\cc$-1- proximal elements (Proposition \ref{p3.2}). This proves Proposition 
\ref{c1pre}.
\end{proof}

\def\tt{\mathbb T}
\def\nn{\mathbb N}

\def\la{\lambda}
\def\La{\Lambda}

\begin{proof} {\bf of Proposition \ref{p3.1}.} Let $n=dim G$. Let $r$ denote the complex 
rank of $G$,  i.e., the minimal multiplicity of  zero eigenvalue of $ad_x$, $x\in\gg$. 
The characteristic polynomial $\chi_{\Ad_g}(\la)$ of the adjoint $Ad_g$, $g\in G$, 
has the type 
$$\chi_{Ad_g}(\la)=(1-\la)^rQ_g(\la), \ Q_g(\la)=\la^{n-r}+1+\sum_{j=1}^{n-r-1}
q_j(g)\la^j,$$
where $q_j(g)$ are polynomials in the matrix coefficients of  the operator $Ad_g$ 
(in some fixed basis of $\gg$). The polynomial $Q_g$ is real. It has unit 
constant term, since 
$\chi_{Ad_g}(0)=\det Ad_g=1$, by semisimplicity (Proposition \ref{ad+-}). 
Its highest coefficient equals $(-1)^{n-r}=1$, by definition and since the number 
$n-r$ is even. This follows from the same proposition. Let 
$$\la_1(g),\dots,\la_{n-r}(g) \text{ be the complex roots of the polynomial } Q_g(\la),$$
which are multivalued functions in $g$. Consider  the following auxiliary subsets  
in $G$: 
$$\La_0=\{ g\in G \ | \ \la_j(g)=1 \text{ for some } j\},$$
$$\La_1=\{ g\in G \ | \ \la_{j_1}(g)=\la_{j_2}(g) \text{ for some } j_1\neq j_2\},$$
$$\La_2=\{ g\in G \ | \ \la_{j_1}(g)-1=1-\la_{j_2}(g) \text{ for some } j_1\neq j_2\},$$
$$\La_3=\{ g\in G \ | \ (\la_{j_1}(g)-1)^2=(\la_{j_2}(g)-1)(\la_{j_3}(g)-1) 
\text{ for some distinct } j_1,j_2,j_3\},$$
$$\La_4=\{ g\in G \ | \ (\la_{j_1}(g)-1)(\la_{j_2}(g)-1)=(\la_{j_3}(g)-1)(\la_{j_4}(g)-1) 
\text{ for some distinct  } j_1,\dots,j_4\},$$
$$U=G\setminus\cup_{i=0}^4\La_i.$$
 
The set $U$ is Zariski open in $G$, since the sets $\La_i$ are Zariski closed in $G$.  
(A priori, a Zariski open subset may be empty.) Each $g\in U$ satisfies statements 
1) and 2) of Proposition \ref{p3.1}. Indeed, if a $g\in G$ does not satisfy statement 1), 
then $g\in\La_0\cup\La_1\subset G\setminus U$ by definition. If a $g\in G$ does not 
satisfy statement 2), then there exist eigenvalues $\nu_1,\nu_2$ of $Ad_g$ such that 
\begin{equation}|\nu_1-1|=|\nu_2-1| \text{ and } \nu_1\neq\nu_2,\bar\nu_2,
\label{nu12}
\end{equation}
which implies that $\nu_1,\nu_2\neq1$. Hence, 
\begin{equation}\nu_1=\la_{j_1}(g)\neq \nu_2=\la_{j_2}(g), \ j_1\neq j_2. 
\text{ We assume that } j_1=1, \ j_2=2,
\label{nul}\end{equation}
without loss of generality, choosing appropriate numeration of the $\la_j$' s. 

\medskip

{\bf Claim 1.} {\it One has $g\in\La_2\cup\La_3\cup\La_4\subset G\setminus U.$}
%\label{la14g}
\medskip

\begin{proof}
Case 1: $\la_1(g),\la_2(g)\in\rr$. Then $\la_1(g)-1=1-\la_2(g)$,  by 
(\ref{nu12}), (\ref{nul}). Thus,  $g\in\La_2.$

Case 2: $\la_1(g)\in\rr$, $\la_2(g)\notin\rr$. Set  
$\la_3(g)=\overline{\la_2(g)}$, which is also a root of  $Q_g(\la)$. By (\ref{nu12}), 
$$(\la_1(g)-1)^2=(\la_2(g)-1)(\la_3(g)-1), \text{ hence } g\in\La_3.$$

Case 3: $\la_1(g),\la_2(g)\notin\rr$. Set  
$\la_3(g)=\overline{\la_1(g)}$, $\la_4(g)=\overline{\la_2(g)}$.  By  (\ref{nu12}), 
(\ref{nul}), 

$(\la_1(g)-1)(\la_3(g)-1)=(\la_2(g)-1)(\la_4(g)-1)$. Hence $g\in\La_4.$
\end{proof}

Thus, the set $U$ is Zariski open in $G$ and consists of elements satisfying statements 
1) and 2). Now for the proof of Proposition \ref{p3.1} it suffices to show that 
the set $U$ 
is nonempty. Fix a Cartan subalgebra $\gh\subset\gg$ (see Subsection 2.5). Let  
$T=\exp\gh\subset G$ be the corresponding maximal connected 
commutative Lie subgroup. 

\medskip

{\bf Claim 2.} $U\cap T\neq\emptyset$.
%\label{uint}\end{equation}
\medskip

\begin{proof} The algebra $\gh_{\cc}$ has $n-r$ distinct roots.  
For every $g=\exp v\in T$, $v\in\gh$, one has 
\begin{equation}
\la_j(g)=\exp(\a_j(v)), \ 
\a_1,\dots,\a_{n-r}\in\gh_{\cc}^*\setminus0 \text{ are the roots of } \gh_{\cc}.
\label{laroot}\end{equation} 
Suppose that $U\cap T=\emptyset$. This means that  the sets 
$\La_0,\dots,\La_4$ (which are Zariski closed in $G$) cover $T$. 
One of the $\La_i$' s  contains $T$, since $T$ is a smooth connected analytic 
variety. 

Case 1: $T\subset\La_0$. Then by (\ref{laroot}), $\a_j\equiv0$ on $\gh$ for some $j$. 
Thus, $\a_j=0\in\gh_{\cc}^*$, - a contradiction.

Case 2: $T\subset\La_1$. Then $\a_{j_1}\equiv\a_{j_2}$ on $\gh$ for some distinct 
roots $\a_{j_1}$ and $\a_{j_2}$, - a contradiction. 

Case 3: $T\subset\La_2$. Then there are distinct roots, we number  them 
as $\a_1$, $\a_2$, such that 
$$e^{\a_1(v)}-1=1-e^{\a_2(v)} \text{ for every } v\in\gh.$$
Taking the derivative at $v=0$, one gets 
$\a_1\equiv-\a_2$. Thus, 
$$e^{\a_1(v)}-1\equiv1-e^{-\a_1(v)}\equiv e^{-\a_1(v)}(e^{\a_1(v)}-1).$$
Hence, $\a_1\equiv0$, - a contradiction. 

Case 4: $T\subset\La_3$. Then there exists a triple of distinct roots 
$\a_1$, $\a_2$, $\a_3$ such that 
$$(e^{\a_1(v)}-1)^2=(e^{\a_{2}(v)}-1)(e^{\a_{3}(v)}-1) \text{ for every } v\in\gh.$$
Taking the second derivative at $v=0$, we get 
\begin{equation}
\a_{1}^2(v)=\a_{2}(v)\a_{3}(v) \text{ for every } v\in\gh, \text{ and hence, 
for every } v\in\gh_{\cc},\label{newm}\end{equation}
by analyticity. Therefore, the zero hyperplanes in $\gh_{\cc}$ of $\a_{1}$, $\a_2$, 
$\a_3$ coincide,  thus, the latter roots are proportional and 
$\a_2,\a_3=\pm\a_1$ (statement d) from Subsection 2.5). Hence, the roots 
$\a_1,\a_2,\a_3$ are not distinct, -  a contradiction. 

Case 5: $T\subset\La_4$. Then there exist 4 distinct roots $\a_1,\dots,\a_4$, such that 
\begin{equation}
(e^{\a_1(v)}-1)(e^{\a_{2}(v)}-1)\equiv (e^{\a_{3}(v)}-1)(e^{\a_4(v)}-1) \ 
\text{ on } \gh.\label{eav}\end{equation}
Taking the second derivative at $v=0$, we get (as above) 
$\a_1(v)\a_2(v)\equiv\a_3(v)\a_4(v)$ on $\gh$, hence on $\gh_{\cc}$. 
Then as in (\ref{newm}), one of the roots $\a_3$, $\a_4$ (say, $\a_3$) is 
proportional to $\a_1$, and $\a_4$ is proportional to $\a_2$. Since the roots 
are distinct, one has  
\begin{equation}\a_3=-\a_1, \ \a_4=-\a_2, 
\label{a-a}\end{equation}
$$(e^{\a_1(v)}-1)(e^{\a_2(v)}-1)\equiv(e^{-\a_1(v)}-1)(e^{-\a_2(v)}-1)\equiv
e^{-(\a_1(v)+\a_2(v))}(e^{\a_1(v)}-1)(e^{\a_2(v)}-1),$$
by (\ref{eav}) and (\ref{a-a}). 
Therefore, $\a_1\equiv-\a_2$. This together with (\ref{a-a}) implies that 
$\a_1=\a_4$, and the roots $\a_1,\dots,\a_4$ are not distinct, - a contradiction.
This proves Claim 2.
\end{proof}

Thus, $U$ is a nonempty Zariski open subset in $G$. This proves Proposition \ref{p3.1}.
\end{proof}

\subsection{A preparatory linear algebra for $\cc$-1-proximal elements}

The results of this subsection (Propositions \ref{lem3} and \ref{lem4}) 
will be used in the proof of Theorem \ref{tconj} for  
groups without proximal elements. They deal with 

- a $n$- tuple of real-linear complex-valued forms $\sigma_1,\dots,\sigma_n:\rr^n\to\cc$ such that the linear operator $\Sigma=(\sigma_1,\dots,\sigma_n):\rr^n\to\cc^n$ is 
injective, 

- a $n$- dimensional real vector 
space $E$ and a $n$- tuple of two-dimensional subspaces $\Lambda_j\subset E$ 
equipped with structures of complex lines compatible with their $\rr$- linear structures, 

- a $n$- tuple $Y=(Y_1,\dots,Y_n)$, $Y_j\in\Lambda_j$, and the linear operator 
$\Omega_{Y,\Sigma}=\sum Y_j\sigma_j:\rr^n\to E$. 

The statements of both propositions say that under certain conditions the operator 
$\Omega_{Y,\Sigma}$ is invertible. We state Proposition \ref{lem3} and use both 
propositions in the case, 
when $E=\gg$ and $\Lambda_j$ are the dominating planes for 
appropriate $n$- tuple of $\cc$-1-proximal elements of $G$. 

\begin{remark} \label{complex} The products $Y_j\sigma_j$ in the above sum are linear operators $\rr^n\to \Lambda_j\subset E$, and their sum is a linear operator $\rr^n\to E$. 
\end{remark}

\begin{proposition} \label{lem3} Let $G$ be a real Lie group with irreducible 
adjoint representation, $n=dim G$. 
Let $A_1,\dots,A_n\in G$ be arbitrary collection of  $\cc$-1-proximal elements, 
$L(A_j)\subset\gg$ be the corresponding dominating planes equipped with the
 $s(A_j)$- complex structures (see Definition \ref{defstr}). 
Let $\Sigma=(\sigma_1,\dots,\sigma_n)$ be as above.  
Then there exist $r_1,\dots,r_n\in\{1,\dots,n\}$ (not necessarily distinct) 
and  elements  $H_1,\dots,H_n\in G$, so that for 
$$\wt A_j=H_jA_{r_j}H_j^{-1}, \ \wt\sigma_j=\sigma_{r_j}, \ 
\wt\Sigma=(\wt\sigma_1,\dots,\wt\sigma_n),$$
there exists a $n$- tuple $Y=(Y_1,\dots,Y_n)$,   $Y_j\in L(\wt A_j)$, 
for which the linear operator 
$\Omega_{Y,\wt\Sigma}=\sum_{j=1}^nY_j\wt\sigma_j:\rr^n\to\gg$ is invertible. 
\end{proposition}

\begin{proof} We construct elements $H_j$, indices $r_j$ and $Y_j\in L(\wt A_j)=
Ad_{H_j}L(A_{r_j})$ by induction on $j$  so that for every $j=1,\dots,n$ the linear operator 
\begin{equation}
\Omega_j=\sum_{i=1}^jY_i\sigma_{r_i}:\rr^n\to\gg \text{ has kernel } 
K_j=Ker\Omega_j \text{ of codimension at least } j.\label{ojcod}\end{equation}
By definition, $\Omega_{Y,\wt\Sigma}=\Omega_n$. This together with 
(\ref{ojcod})  proves the proposition. 

Induction base: $j=1$. Recall that $Ker\Sigma=0$. Take an arbitrary $r_1$ so that 
$\sigma_{r_1}\not\equiv0$  and arbitrary  $H_1\in G$,  
$Y_1\in L(\wt A_1)\setminus0$.  Then $\Omega_1=Y_1\sigma_{r_1}\not\equiv0$ 
and statement (\ref{ojcod}) is obvious.

Induction step: $1<j\leq n$. Let we have already chosen $r_i$, $H_i$, $Y_i$ for 
$i\leq j-1$ so that (\ref{ojcod}) holds with $j$ replaced by $j-1$.  
Let $\Omega_{j-1}$ be the corresponding operator, $K_{j-1}=Ker\Omega_{j-1}$. 
Let us construct $r_j$, $H_j$ and $Y_j$ for which (\ref{ojcod}) holds. 

By the induction hypothesis, 
$$codim K_{j-1}\geq j-1.$$ 
If $codim K_{j-1}\geq j$, then (\ref{ojcod}) holds with $Y_j=0$ and 
arbitrary $r_j$, $H_j$: in this case $K_j=K_{j-1}$. Thus, everywhere below we 
assume that
\begin{equation}
dim\Omega_{j-1}(\rr^n)=codim K_{j-1}=j-1<n, \text{ in particular, } K_{j-1}\neq0, \ 
\Omega_{j-1}(\rr^n)\neq\gg.\label{dimcodim}\end{equation}
Let us fix a $r_j\in\{1,\dots,n\}$ such that 
\begin{equation}\sigma_{r_j}|_{K_{j-1}}\not\equiv0.\label{sigmaneq}\end{equation}
It exists by (\ref{dimcodim}) and since $Ker\Sigma=0$. The irreducibility of the adjoint 
and (\ref{dimcodim}) imply that there exists a 
$H_j\in G$ such that for $\wt A_j=H_jA_{r_j}H_j^{-1}$ one has 
\begin{equation}L=L(\wt A_j)=Ad_{H_j}L(A_{r_j})\not\subset\Omega_{j-1}(\rr^n). 
\text{ Set  } \Lambda=\Omega_{j-1}(\rr^n)\cap L.\label{ola}\end{equation}
Fix this $H_j$. Then $\Lambda$ is a linear 
subspace in the plane $L$, and $\Lambda\neq L$. Thus, either $\Lambda=0$, 
or $\Lambda$ is a line in $L$. For every choice of  $Y_j\in L$ one has 
\begin{equation} 
\Omega_j=\Omega_{j-1}+Y_j\sigma_{r_j}, \ K_{j-1},K_j\subset P=\Omega_{j-1}^{-1}(\Lambda), \ codim P=codim K_{j-1}-dim\Lambda: 
\label{koj}\end{equation}
\begin{equation}
\Omega_j(v)=0 \text{ if and only if } \Omega_{j-1}(v)=-Y_j\sigma_{r_j}(v)\in 
\Omega_{j-1}(\rr^n)\cap L=\Lambda.\label{nstar}\end{equation}

Case 1: $\Lambda=0$. Take an arbitrary 
$Y_j\in L\setminus0$. Then the images of $\rr^n$ under the linear operators 
$\Omega_{j-1}$ and 
$-Y_j\sigma_{r_j}$ are linear subspaces in $\gg$ with zero intersection. This together with  (\ref{nstar}), implies that the kernel $K_j$ is the intersection of the kernels of the 
latter operators: $K_j=Ker(\sigma_{r_j}|_{K_{j-1}})$, $codim K_j\geq
codim K_{j-1}+1=j$, by (\ref{sigmaneq}). This proves (\ref{ojcod}). 

Case 2: $\La$ is a line. Recall that $K_{j-1}\subset P$ by (\ref{koj}). 
Hence, by  (\ref{sigmaneq}) and (\ref{koj}), 
\begin{equation}\sigma_{r_j}(P)\neq0, \ codim P=codim K_{j-1}-1=j-2.\label{AC}
\end{equation}
Consider the following subcases, where $\sigma_{r_j}(P)\subset\cc$ 
is either a real line, or the entire $\cc$. 

\def\oj{\Omega_j}
\def\ojl{\Omega_{j-1}}
\def\srj{\sigma_{r_j}}

Subcase 2.1: $\sigma_{r_j}(P)\subset\cc$ is a line. Take a $Y_j\in L\setminus0$ so that 
the lines $\hat l=Y_j\sigma_{r_j}(P)$ and $\La$ be distinct. These lines, which 
have zero intersection,  are the  images of the space $P$ under the operators from 
Case 1, and $K_{j-1},K_j\subset P$ by (\ref{koj}). This implies (\ref{ojcod}), as in Case 1.

 Subcase 2.2: $\srj(P)=\cc$. Fix an arbitrary $Y_j\in L\setminus0$. Then
 $Y_j\srj(P)=L\supset\Lambda$. Set 
$$P'=\{ v\in P \ | \ Y_j\srj(v)\in\La\}: \  
P' \text{ is a hyperplane in } P, \ codim P'=j-1,$$
 by (\ref{AC}). The kernel $K_j$ is the subspace in $P'$ defined by the linear equation
 \begin{equation} \Omega_{j-1}(v)=-Y_j\srj(v),\label{yjv}\end{equation}
by (\ref{nstar}). The left-hand side in (\ref{yjv}) is independent on $Y_j$. 
Its right-hand side does not 
 vanish identically on $P'$: $\srj|_{P'}\not\equiv0$, since $P'\subset P$ is 
 a hyperplane,  and $\srj(P)=\cc$. 
 The subspace $P'$ depends on $Y_j$, but does not change when 
 we multiply $Y_j$ by a real constant. Hence, multiplying $Y_j$ by a  
 constant, one can achieve  that 
  equation (\ref{yjv}) does not hold identically in $v\in P'$. Then its solution 
  space $K_j$ is smaller than $P'$, thus, $codim K_j\geq codim P'+1=j$. The 
 induction step is over. This proves (\ref{ojcod}) and Proposition \ref{lem3}.
 \end{proof}

\def\La{\Lambda}
\begin{proposition} \label{lem4} Let $n\in\mathbb N$, $E$, 
$\Lambda_1,\dots,\Lambda_n\subset E$ 
be the same, as at the beginning of the subsection. Let 
$\Sigma=(\sigma_1,\dots,\sigma_n):\rr^n\to\cc^n$ be a linear operator 
suct that there exist $Y_j\in\La_j$, for which the corresponding operator 
$\oys$ is invertible. Let 
$s_1,\dots,s_n\in\cc\setminus0$ be a $n$- tuple of complex numbers such that 
$0<|s_j|<1$ and $\a_j=\frac{\arg s_j}{\pi}\notin\cup_{q=1}^n\frac1q\mathbb Z$. 
Then for every other  
collection $\wt Y=(\wt Y_1,\dots,\wt Y_n)$, $\wt Y_j\in\La_j\setminus0,$ 
there exists an $l=(l_1,\dots,l_n)\in\zz^n$ such that for  
$$Y_j'=s^{l_j}_j\wt Y_j, \ Y'=(Y_1',\dots,Y_n')$$
the operator 
$\Omega_{Y',\Sigma}=\sum_{j=1}^nY_j'\sigma_j:\rr^n\to E$ 
is also invertible. 
\end{proposition}

\def\olj{\overline L_j}
\def\ol{\overline L}

\begin{proof} Consider the set 
$$Reg=\{Y=(Y_1,\dots,Y_n)\in  
\La_1\times\dots\times \La_n=\rr^{2n} \ | \ \text{the operator }  \oys 
\text{ is invertible}\}.$$
The set $Reg$ is nonempty, by the conditions of Proposition \ref{lem4}. 
Fix an arbitrary collection $\wt Y=(\wt Y_1,\dots,\wt Y_n)$, 
$\wt Y_j\in\La_j\setminus0$. Set 
\begin{equation}L_j=\{ s^q_j\wt Y_j \ | \ q\in\zz\}\subset\La_j, \ 
L=L_1\times\dots\times L_n\subset\rr^{2n}.\label{r2n}\end{equation}
The statement of the proposition is equivalent to the inequality 
$L\cap Reg\neq\emptyset$. 

\medskip

{\bf Claim 1.} {\it The complement $\rr^{2n}\setminus Reg$ is the zero set 
of a nonzero homogeneous polynomial of degree $n$.}

\begin{proof} 
Fix some bases in $\rr^n$, $E$ and in the planes $\Lambda_j$. 
The above  
complement is the zero set of the determinant $\det\oys$ of the matrix of the 
operator $\oys:\rr^n\to E$. The latter operator depends linearly on the vector 
parameter $Y$. Thus, $\det\oys$ is a homogeneous polynomial 
of degree $n$ in the components of $Y_j$. One has $\det\oys\not\equiv0$, 
since $Reg\neq\emptyset$. 
\end{proof}

{\bf Claim 2.} {\it The set $L$ is not contained in the zero 
set of a nonzero polynomial of degree $n$.}
 
\medskip

For the proof of Claim 2 we consider the real Zariski closures   
$\overline L_j\subset\La_j$, $\overline L\subset\rr^{2n}$ of the sets 
$L_j$  and $L$  and split the set $\{1,\dots,n\}$ into 
the two following auxiliary subsets $I$ and $J$:  
\begin{equation}
I=\{ j=1,\dots,n \ | \ \a_j\notin\mathbb Q\}, \ 
J=\{ j=1,\dots,n \ | \ \a_j=\frac{p_j}{q_j}\in\mathbb Q\}.\label{IJ}
\end{equation} 
Let $\La_{j,q}$ denote the line $\rr s^q_j\wt Y_j\subset \La_j$.  Then 
\begin{equation}\olj=\La_j \text{ if } j\in I; \ \olj=
\cup_{q=0}^{q_j-1}\La_{j,q} \text{ if } 
j\in J;  \ \ol=
\overline L_1\times\dots\times\overline L_n.\label{oljq}\end{equation}
Statement (\ref{oljq}) follows from definition and the convergence $s^q_j\to0$, as 
$q\to+\infty$ (recall that $0<|s_j|<1$). Recall 
that $q_j>n$ for every $j\in J$, by the conditions of Proposition \ref{lem4}. Thus, for 
$j\in J$, the set $\olj$ is a union of $q_j>n$ lines $\La_{j,q}$.

\begin{proof} {\bf of Claim 2.} We prove the claim by contradiction. 
Suppose that there exists a 
nonzero polynomial $P$ of degree $n$ that vanishes 
identically on $L$, and hence, on $\ol$. We show that $P\equiv0$ on 
$\rr^{2n}$, - a contradiction to the 
nontriviality of $P$. To do this, we introduce the auxiliary subspaces 
$L_{J',D_{J'}}\subset\rr^{2n}$ defined as follows. Given  
\begin{equation}
J'\subset J, \ D_{J'}=(d_j)_{j\in J'}, \ d_j\in\mathbb Z, \  
0\leq d_j\leq q_j-1, \label{j'j}\end{equation}
set
\begin{equation}L_{J',D_{J'}}=(\prod_{j\in I\cup(J\setminus J')}\La_j)\times
(\prod_{j\in J'}\La_{j,d_j}). 
\label{lla}\end{equation}
In particular, $L_{\emptyset}=\rr^{2n}$. Then we have  
$\ol=\cup_{D_J} L_{J,D_J}$, by (\ref{oljq}). We show that 
\begin{equation}
P|_{L_{J',D_{J'}}}\equiv0 \text{ for every } L_{J', D_{J'}}  
\text{ as in (\ref{lla}), including } L_{\emptyset}=\rr^{2n},\label{indst}
\end{equation}
by induction on the cardinality $|J\setminus J'|$.

Induction base: $|J\setminus J'|=0$. Then $J'=J$, $L_{J,D_J}\subset \ol$, 
see (\ref{lla}), hence, $P\equiv0$ on $L_{J,D_J}$. 

Induction step: $|J\setminus J'|=k\geq1$. We assume the induction hypothesis 
holds: 
$$P|_{L_{J'',D_{J''}}}\equiv0 \text{ for every } J'',D_{J''} \text{ as in (\ref{j'j}) with } 
|J\setminus J''|
<k.$$
Take some $D_{J'}=(d_j)_{j\in J'}$, as in (\ref{j'j}), and 
$s\in J\setminus J'$. Set  
$$J''=J'\cup\{ s\}.$$
For every $q=0,\dots,q_s-1$ set 
$$d_{j,q}=d_j \text{ for } j\in J', \ d_{s,q}=q, \ 
D_{J''}(q)=(d_{j,q})_{j\in J''}.$$
The linear space $L_{J',D_{J'}}$ contains $q_s$ distinct hyperplanes $L_{J'',D_{J''}(q)}$, 
$q=0,\dots,q_s-1$, on which $P\equiv0$ (the induction hypothesis). One has 
$n=deg P<q_s$. Hence, $P\equiv0$ on the ambient  
space $L_{J',D_{J'}}$. The induction step is over. Statement (\ref{indst}) is 
proved.

Statement (\ref{indst}) and its above proof work for the space 
$L_{\emptyset}=\rr^{2n}$, hence,   
$P|_{\rr^{2n}}\equiv0$, - a contradiction to the nontriviality of the 
polynomial $P$. This proves Claim 2.
\end{proof}
 
Proposition \ref{lem4} follows immediately from Claims 1 and 2.
\end{proof}

\subsection{Case of group $Aff_+(\rr)$.} 
Recall that $Aff_+(\rr)$ is the group of orientation-preserving affine transformations 
$\rr\to\rr$. For every $s>0$ and $u\in\rr$ set 
$$g_s:x\mapsto sx, \ t_u:x\mapsto x+u, \ g_s,t_u\in Aff_+(\rr).$$
Define the group $\G=<g,t\ |\ [[g,t],t]=1>$. For each $s\in\rr$ define the representation 
$\rho_s\in Hom(\G,Aff_+(\rr))$: $\rho_s(g)=g_s$, $\rho_s(t)=t_1$.   

\begin{proposition} \label{paff} 
For every $s_0>0$ there exist sequences $s_k>0$ and $\gamma_k\in\G$ such that  
$$s_k\to s_0, \text{ as } k\to\infty, \ \rho_{s_k}(\gamma_k)=1, \text{ but } 
\rho_s(\gamma_k)\not\equiv1 \text{ in } s\in\rr_+.$$ 
\end{proposition}

\begin{proof} It suffices to prove the statement of the proposition 
for open and dense subset of the values $s_0>0$ (afterwards we pass 
to the closure and diagonal sequences). Thus, without loss of 
generality we assume that $s_0\neq1$. We also assume that $0<s_0<1$, 
since $\rho_s(\G)=\rho_{s^{-1}}(\G)$. 

Each group $\rho_s(\G)$ contains the elements  
$ t_{ms^k}=\rho_s(g^kt^mg^{-k})$, $m\in\mathbb Z,\ k\in\mathbb N$. Set 
$$m_k=[s_0^{-k}], \ \gamma_k=g^kt^{m_k}g^{-k}t^{-1}.$$
The relation $\rho_s(\gamma_k)=1$ is equivalent to the equation $m_ks^k=1$. 
It obviously does not hold identically and has a solution $s_k\to s_0$, as $k\to\infty$. 
Indeed, consider the new parameter $\wt u=k(s-s_0)$ and mappings 
$\psi_k:\wt u\mapsto m_k(s_0+k^{-1}\wt u)^k$. One has $m_ks^k=\psi_k(\wt u)$, 
$m_ks^k_0\to1$, 
\begin{equation}\psi_k(\wt u)=m_ks_0^k(1+k^{-1}\frac{\wt u}{s_0})^k
\to 
\psi(\wt u)=e^{\frac{\wt u}{s_0}},\ \text{as}\ k\to\infty. 
\label{mskto}\end{equation}
The convergence is uniform with derivatives on compact sets. 
The limit is a diffeomorphism $\psi:\rr\to\rr_+$, and $\psi(0)=1$. Hence, 
the equations $\psi_k(\wt u)=1$ have solutions $\wt u_k\to0$, and the equations 
$\rho_s(\gamma_k)=1$ have solutions $s_k=s_0+k^{-1}\wt u_k\to s_0$. 
The proposition is proved. 
\end{proof}

The next  generalization of the limit (\ref{mskto}) will be used in the proof of 
Theorem \ref{tconj}.

\begin{proposition} \label{pann} Let $s:\rr^n\to\cc$ be a smooth function on a neighborhood of zero, 
$0<|s|<1$. Let $\var\in\cc$ and sequences $k_r, m_r\in\zz$ be such that 
$k_r\to+\infty$ and $m_rs^{k_r}(0)\to\var$. Then $m_rs^{k_r}(k_r^{-1}\wt u)\to \var e^{(d\ln s(0))\wt u}$, 
as $r\to\infty$, uniformly with derivatives on compact sets in $\rr^n$. 
\end{proposition}
\begin{remark} The logarithmic differential $d\ln s(0)$ is well-defined, since $s(0)\neq0$.
\end{remark}
\begin{proof} Analogously to (\ref{mskto}), one has 
$$s^{k_r}(k_r^{-1}\wt u)=(s(0)+k_r^{-1}(ds(0))\wt u+o(k_r^{-1}))^{k_r}=
s^{k_r}(0)(1+k_r^{-1}(d\ln s(0))\wt u+o(k_r^{-1}))^{k_r}.$$
The latter power of expression in the brackets tends to $e^{(d\ln s(0))\wt u}$, while 
$m_rs^{k_r}(0)\to\var$, by assumption. This proves the proposition.
\end{proof}

\subsection{Normal form of contracting dynamics with dominated splitting}

\begin{definition} A linear operator 
is {\it contracting}, if all its complex eigenvalues lie in the unit disk. 
\end{definition}

\begin{theorem} \label{dsplit}
 Let $\phi:(\cc^n,0)\to(\cc^n,0)$ be an analytic mapping  germ at its fixed point 
 $0$. Assume that its differential 
$F=d\phi(0)$ be contracting and have complementary invariant subspaces $E_{\pm}$: 
 $T_0\cc^n=E_+\oplus E_-$. Set 
$$F_{\pm}=F|_{E_{\pm}}, \  \pi_{\pm}:T_0\cc^n\to E_{\pm} \text{ (the projection 
along } E_{\mp}).$$
Assume that for each (complex) eigenvalue $\la_+$ of $F_+$ and any eigenvalue 
$\la$ of $F$ one has 
\begin{equation}\la_+\neq0, \ |\la_+^{-1}\la^2|<1.\label{la+la}\end{equation}
Then there exists a unique  germ of analytic mapping $x_+:(\cc^n ,0)\to(E_+,0)$ such that 
\begin{equation} x_+\circ\phi=F_+\circ x_+ \text{ and } dx_+(0)=\pi_+.
\label{vphi}\end{equation}
Thus, for any germ of analytic vector function $x_-:(\cc^n,0)\to (E_-,0)$ with 
$dx_-(0)=\pi_-$ one has local coordinates $x=(x_+,x_-)$ near 0, in which 
the mapping $\phi$ can be written as 
\begin{equation}\phi:(x_+,x_-)\mapsto(F_+x_+,F_-x_-+O(||x||^2)), \text{ as } x\to0.
\label{coord}\end{equation}
\end{theorem}

\def\bd{\mathcal B_{\delta}}
\def\cs{\mathcal S}
\def\la{\lambda}

\begin{remark}
In the case, when the differential $d\phi(0)$ is invertible, Theorem \ref{dsplit} 
follows from Poincar\'e-Dulac theorem (\cite{A}, chapter 5, section 25, subsection D). 
The author is pretty sure that the general statement of 
Theorem \ref{dsplit} is well-known to specialists, but he did not find it in literature. A 
proof of Theorem \ref{dsplit} is given below. 
\end{remark}

\begin{corollary} \label{cpoc} 
Let $\phi$, $E_{\pm}$, $F_+$, $x_+$ be the same, as in 
Theorem \ref{dsplit}, and let $T_0\cc^n=E_+\oplus E_-$ be a dominated splitting for 
$d\phi(0)$ (see Definition \ref{ddsplit}). Fix an arbitrary germ of biholomorphism $(T_0\cc^n,0)\to(\cc^n,0)$ tangent to the identity. Then 
\begin{equation}\phi^k(h)=F_+^kx_+(h)+o(||F_+^k||||h||), 
\text{ as } k\to\infty.\label{phicla}
\end{equation}
The vector $F_+^kx_+(h)\in E_+\subset T_0\cc^n$ represents a point in $\cc^n$ 
 via the above biholomorphism. 
The  ``o'' in (\ref{phicla}) is uniform with derivatives on  a neighborhood of zero 
(independent on $k$). 
\end{corollary}

{\bf Addendum to Theorem \ref{dsplit} and Corollary \ref{cpoc}.} {\it If in Theorem 
\ref{dsplit} the mapping germ $\phi$ and the splitting $T_0\cc^n=E_+\oplus E_-$ 
are complexifications of a real mapping germ $\phi:(\rr^n,0)\to(\rr^n,0)$ and 
a $d\phi(0)$- invariant splitting of $T_0\rr^n$, then the corresponding vector 
function $x_+$ is also real-valued on $\rr^n$. If $\phi$ and the above splitting depend analytically on some parameter, then $x_+$ also depends analytically on the same 
parameter. If in Corollary \ref{cpoc} the identification 
$(T_0\cc^n,0)\to(\cc^n,0)$ is fixed and the germ $\phi$ 
  depends analytically on a parameter, then 
the corresponding ``o'' in (\ref{phicla}) is uniform with derivatives on compact 
sets in the product of the phase and parameter spaces.}

The corollary and the addendum are proved below.

Let $G$ be a Lie group,  $\Pi_{\rr,1}$, $\Pi_{\cc,1}$, $s(g)$ be the same, as in 
(\ref{prox<1}) and (\ref{pc1}). Set
\begin{equation}
\hat\Pi=\{ g\in G \text{ is either 1- proximal, or } \cc-\text{1-proximal} \ | \ |s(g)|<1\}
=\Pi_{\rr,1}\sqcup\pc1.\label{hatpi}\end{equation}
\begin{remark} \label{hat=pi} Let $G$ be a semisimple linear algebraic group. 
Then the set $\hat\Pi$ is open and accumulates to the identity  
(Theorem \ref{1prpr} and Proposition \ref{c1pre}); $G$ has no proximal 
elements if and only if $\hat\Pi=\pc1$. 
\end{remark}
We apply Theorem \ref{dsplit}, 
Corollary \ref{cpoc} and the addendum to the following commutator mapping $\phi_g$ 
with $g\in\hat\Pi$:
\begin{equation}\phi_g:G\to G, \ h\mapsto ghg^{-1}h^{-1}: \ \  
\phi_g(1)=1, \ d\phi_g(1)=Ad_g-Id:\gg\to\gg.\label{commap}\end{equation}
In what follows, for every 1-proximal or $\cc$-1-proximal element $g\in\hat\Pi$ let 
$$L(g)\subset\gg\text{ denote the  dominating subspace (line or plane) of the operator 
 } Ad_g-Id.$$
%\begin{equation}s(g) \text{ denote the corresponding (complex) eigenvalue.}
%\label{sg}\end{equation}
In the case, when $g$ is $\cc$-1-proximal, see the convention of choice of the 
complex eigenvalue $s(g)$ in Subsection 2.7 (after Proposition \ref{p3.4}). 
\begin{remark} \label{sneq0} The  subspaces $L(g)$ and the eigenvalues $s(g)$ depend analytically on $g\in\hat\Pi$. The analyticity in the  1-proximal 
elements $g$ follows from the simplicity of the eigenvalue $s(g)$. 
The analyticity on $\pc1$ follows from Proposition \ref{p3.4} and the discussion afterwards.
\end{remark}

\begin{proposition} \label{pcomm} Let $G$ be a Lie group, $\hat\Pi$, $L(g)$  
be as above, $\hat\Pi\neq\emptyset$. There exist an open subset
\begin{equation}\Pi'\subset\hat\Pi\times G, \ \Pi'\supset\hat\Pi\times 1,\label{pi'}
\end{equation}
 and a $\gg$- valued vector function $v_g(y)$  analytic in $(g,y)\in\Pi'$,   
 $v_g(1)=0$, 
 such that  for every  $(g,y)\in\Pi'$ one has $v_g(y)\in L(g)$, 
\begin{equation} dv_g|_{L(g)}=
Id:L(g)\to L(g), 
\ \phi^k_g(y)=\exp(s^k(g)v_g(y)+o(|s^k(g)|dist(y,1))),\label{vg}\end{equation} 
 as   $k\to+\infty$.  
The latter ``$o$'' is uniform with derivatives in $(g,y)$ on compact subsets in $\Pi'$. 
\end{proposition}  
\begin{remark} \label{remdv} In the above proposition and in what follows 
$dv_g:\gg\to\gg$ is the differential of the mapping $v_g:G\to \gg$ 
at the point $1\in G$, $\phi_g^k$ is 
the $k$- th iteration of the mapping $\phi_g$. In the case, when $g$ is $\cc$-1-proximal, 
the multiplication of vectors in the plane $L(g)$ by complex numbers is defined in terms 
of the $s(g)$- complex structure, see Definition \ref{defstr}. 
\end{remark}
\begin{proof} {\bf of Proposition \ref{pcomm}.} For every $g\in\hat\Pi$ the complexified mapping $\phi_g$ satisfies 
the conditions of Theorem \ref{dsplit} and Corollary \ref{cpoc}. Indeed, 
the real mapping $\phi_g$ has contracting differential $F_{\rr}=Ad_g-Id$ at 1 with dominated splitting, 
$E_{+,\rr}=L(g)$, $F_{+,\rr}:L(g)\to L(g)$ being the multiplication by $\lambda_+=s(g)$. 
The latter operator on $L(g)$ has at most two distinct complex eigenvalues, and they 
have the same modulus $|s(g)|\neq0$. For 
each complex eigenvalue $\lambda$ of $F_{\rr}$ one has $|\lambda|\leq|s(g)|<1$, 
by definition and ($\cc$-)1- proximality. This implies (\ref{la+la}): 
$$|\lambda_+|=|s(g)|\neq0, \ |\lambda_+^{-1}\lambda^2|=|(s(g))^{-1}\lambda^2|\leq
|\lambda|\leq|s(g)|<1.$$
 Let $v_g(y)=x_+(y)$ be the corresponding vector function from (\ref{vphi}). Now 
the statements of Proposition \ref{pcomm} follow from (\ref{vphi}), 
Corollary \ref{cpoc} and its addendum applied to the complexified  local 
exponential identification $\exp:(\gg,0)\to(G,1)$. 
\end{proof}

\begin{corollary} \label{convexp}
Let $G$, $\Pi'$ be the same, as in Proposition \ref{pcomm}. Let 
$g(u)$, $h(u)$ be two families of elements of $G$ depending smoothly on a 
parameter $u\in\rr^n$. Let $(g(0),h(0))\in\Pi'$. Set
$$\chi_k(u)=\phi_{g(u)}^k(h(u))\in G,$$
$$s_g(u)=s(g(u)), \ \wt\nu_g(u)=v_{g(u)}(h(u)), \nu_g=\wt\nu_g(0)\in L(g(0))\subset\gg.$$
Let $\var\in\rr$ if $g(0)\in\Pi_{\rr,1}$, or $\var\in\cc$ if $g(0)\in\pc1$. Let $m_r,k_r\in\zz$ be two integer sequences such that $k_r\to+\infty$ 
and $m_rs_g^{k_r}(0)\to\var$, as $r\to\infty$. Then
\begin{equation}(\chi_{k_r}(k_r^{-1}\wt u))^{m_r}\to\exp(\var e^{(d\ln s_g(0))\wt u}\nu_g), 
\text{ as } r\to\infty\label{limnew}\end{equation}
uniformly with derivatives on compact sets in $\rr^n$.
\end{corollary}

\begin{proof} The functions $s_g(u)$ and $\wt\nu_g(u)$ are well-defined on a neighborhood of zero (openness of the set $\Pi'$). By (\ref{vg}), one has
$$(\chi_{k_r}(k_r^{-1}\wt u))^{m_r}=\exp(m_r(s_g^{k_r}(k_r^{-1}\wt u)\wt\nu_g(k_r^{-1}
\wt u)+o(|s_g^{k_r}(k_r^{-1}\wt u)|dist(h(k_r^{-1}\wt u),1)))), \text{ as } r\to\infty.$$
The expression under exponent tends to 
$\var e^{(d\ln s_g(0))\wt u}\nu_g$ uniformly with derivatives on compact sets, 
by Proposition \ref{pann} and since $\wt\nu_g(k_r^{-1}\wt u)\to\nu_g$, as $k_r\to\infty$. 
This proves (\ref{limnew}).
\end{proof}

In the proof of Theorem \ref{dsplit} and Corollary \ref{cpoc} we use the following notations. Let $\Lambda_{\pm}$ denote 
the collection of the (complex) eigenvalues of $F_{\pm}:E_{\pm}\to E_{\pm}$. Set
$$\la_{\pm,\min}=\min\{|\la| \ | \ \la\in\Lambda_{\pm}\}, \ \la_{\pm,\max}=
\max\{|\la| \ | \ \la\in\Lambda_{\pm}\}, \la_{\max}=\max\la_{\pm,\max}.$$
Let us fix $\la_1,\la_2,\la_3>0$ such that 
\begin{equation}0<\la_1<\la_{+,\min}\leq\la_{\max}<\la_2<
\la_3<1, \ 
\mu=\la_1^{-1}\la_3^2<1.\label{lala}\end{equation}
They exist by definition 
and inequality (\ref{la+la}). Fix some positive Hermitian scalar products on $E_{\pm}$, 
let $|| \ ||$ denote the norm on $T_0\cc^n$ corresponding to their direct sum, 
such that 
\begin{equation}\la_1||v||\leq 
||F_+v||\leq\la_2||v|| \text{ for every } v\in E_+, \   ||d\phi(0)||\leq\la_2.\label{preineq}\end{equation}
The existence of the scalar products satisfying (\ref{preineq}) follows by 
definition and (\ref{lala}). 
Set $D_{\delta}=\{||h||<\delta\}\subset\cc^n$. By (\ref{lala}) and (\ref{preineq}), 
there exists a $\delta>0$ (let us fix it)  such that 
\begin{equation} \phi \text{ is holomorphic on } 
\overline{D_{\delta}}, \ 
||\phi(h)||\leq\la_3||h|| \text{ for every } h\in\overline{D_{\delta}}, \text{ 
thus, } \phi(\overline{D_{\delta}})\subset D_{\delta}.\label{defde}\end{equation}

\begin{proof} {\bf of Theorem \ref{dsplit}.} 
We follow the classical scheme (cf. \cite{gorb} and \cite{ilyak} p. 65).  
Let us identify $\cc^n$ with $T_0\cc^n$ in the canonical way. 
The subspaces in $\cc^n$ corresponding to $E_{\pm}\subset T_0\cc^n$ will be denoted 
by the same symbols $E_{\pm}$. The norm on $\cc^n$ is the one induced by 
that from (\ref{preineq}). Given $\phi$, define $Q_{\phi}:(\cc^n,0)\to(\cc^n,0)$ by 
$$\phi(h)=F_+\pi_+h+F_-\pi_-h+Q_{\phi}(h): \ Q_{\phi}(h)=O(||h||^2), \text{ as } 
h\to0.$$
The  function $x_+$ we are looking for has the following type:
$$x_+(h)=\pi_+h+Q(h), \  Q(h)\in E_+, \ Q(h)=O(||h||^2), \text{ as } h\to0.$$
 Substituting the latter to (\ref{vphi}), cancelling linear terms and 
applying the inverse $F_+^{-1}$  
translates (\ref{vphi}) into  the fixed point equation $Q=\Phi(Q)$ 
for the following transformation $\Phi$. Set 
\begin{equation}
R(h)=F_+^{-1}\pi_+Q_{\phi}(h), AQ=F_+^{-1}\circ Q\circ\phi, \ \Phi(Q)=R+AQ.
\label{ra}\end{equation}
(The operator $F_+:E_+\to E_+$ is invertible, since $0\notin\Lambda_+$, see 
(\ref{la+la}).) 
Recall that the mapping $\phi$, and hence, $Q_{\phi}$, is holomorphic on 
$\overline{D_{\delta}}$, see (\ref{defde}). 
We show that the transformation $\Phi$ is contracting, and hence, has a 
unique fixed point,  in the Banach space 
$$\bd=\{Q:\overline{D_{\delta}}\to E_+ \ | \ Q \text{ is holomorphic on } 
D_{\delta} \text{ and continuous on } \overline{D_{\delta}}, \  dQ(0)=0\}$$
$$\text{ with the norm } ||Q||_{\delta}=\sup_{\overline{D_{\delta}}\setminus0}
\frac{||Q(h)||}{||h||^2}.$$
 Note that $R$ is a constant linear non-homogeneous operator 
 (independent on $Q$), thus, the 
 contractibility of $AQ+R$ is equivalent to that of $AQ$. 
 
 \medskip
 
 {\bf Claim.} {\it  Let $A$ be the operator from (\ref{ra}). For every 
 $Q\in\bd$ one has $AQ\in\bd$, 
 $||AQ||_{\delta}\leq\mu||Q||_{\delta}$ , where $\mu=\la_1^{-1}\la_3^2<1$, 
 see (\ref{lala}).}
 
 \begin{proof} Let $Q\in\bd$. 
 By definition, $AQ=F^{-1}_+\circ Q\circ\phi$. By (\ref{defde}), one has 
 $Q\circ\phi\in\bd$,  
 $$||Q\circ\phi||_{\delta}\leq\sup_{\overline{D_{\delta}}\setminus\{\phi=0\}}
 \frac{||Q(\phi(h))||}{||\phi(h)||^2}
 \sup_{\overline{D_{\delta}}\setminus0}(\frac{||\phi(h)||}{||h||})^2
 \leq||Q||_{\delta}\la_3^2, \ 
 ||F_+^{-1}||\leq\la_1^{-1},$$
 by (\ref{defde}) and (\ref{preineq}). The two above inequalities together  imply the claim.
 \end{proof}

 Every analytic germ $Q(h)=O(||h||^2)$ with values in $E_+$ 
 represents an element of the space $\bd$ with  a  $\delta$ satisfying (\ref{defde}). 
The transformation $\Phi$ is a contracting self-mapping $\bd\to\bd$, 
by the claim. Hence, it has a unique 
fixed point there. This proves Theorem \ref{dsplit}. 
\end{proof}

\begin{proof} {\bf of Corollary \ref{cpoc}.} Consider the coordinates $x=(x_+,x_-)$, 
see (\ref{coord}). It suffices to prove (\ref{phicla}) for the canonical identification 
$T_0\cc^n=\cc^n$ induced by the coordinates $x$: the asymptotics (\ref{phicla}) 
survives after passing to another identification germ $(T_0\cc^n,0)\to(\cc^n,0)$ 
tangent to the identity. Let $\la_0,\la_1,\la_3,\delta>0$ and the norms on the spaces 
$E_{\pm}$ and $T_0\cc^n$ be the same, as in (\ref{lala})-(\ref{defde}), 
such  that in addition, 
\begin{equation}\la_{-,\max}<\la_0<\la_1<\la_{+,\min}, \ ||F_-||<\la_0.\label{newla}
\end{equation}
These $\la_0$, $\la_1$, $\la_3$, $\delta$ exist by our assumption on the 
dominated splitting. There exists a $c>0$ such that 
$$||x_-\circ\phi(h)-F_-x_-(h)||\leq c||h||^2 \text{ for every } h\in\overline{D_{\delta}},$$
by (\ref{coord}) and (\ref{defde}). Hence, by (\ref{defde}) and (\ref{newla}), 
for each $k>0$, $0\leq l\leq k-1$, and $h\in\overline{D_{\delta}}$ 
$$||F^l_-x_-\circ\phi^{k-l}(h)-F_-^{l+1}x_-\circ\phi^{k-l-1}(h)||\leq 
c||F_-^l||||\phi^{k-l-1}(h)||^2\leq c\lambda_0^l\lambda_3^{2(k-l-1)}||h||^2.$$
Summing these inequalities for all $l$ and taking the maximum of the 
right-hand sides yields
\begin{equation}||x_-\circ\phi^k(h)-F_-^kx_-(h)||\leq ck
\max_{0\leq l<k}(\lambda_0^l\lambda_3^{2(k-l-1)})||h||^2=
ck(\max\{\lambda_0,\lambda_3^2\})^{k-1}||h||^2.\label{x-h}\end{equation}
Recall that $\la_0,\la_3^2<\la_1<1$ by (\ref{lala}) and  (\ref{newla}), 
$||F_+^k||\geq\la_1^k$  
by (\ref{preineq}), $||F_-^k||\leq\la_0^k$ by (\ref{newla}).  Hence, 
$||F_-^k||, k\la_0^{k-1}, k\la_3^{2(k-1)}=o(\la_1^k)=o(||F_+^k||)$, as $k\to\infty$. 
This together with 
(\ref{x-h}) implies that $x_-\circ\phi^k(h)=o(||F_+^k||||h||)$. This 
 implies (\ref{phicla}), by (\ref{coord}). 
\end{proof}

\begin{proof} {\bf of Addendum.} Let the mapping $\phi$ and the  
subspaces $E_{\pm}$ from Theorem 
\ref{dsplit} be complexifications of real ones. Then the contracting 
transformation $\Phi$ from (\ref{ra}) preserves 
the subspace in $\bd$ of the real-valued vector functions on $\rr^n$. 
Hence, its unique fixed 
point $x_+$ also belongs to this subspace. If $\phi$ and $E_{\pm}$ depend analytically 
on a parameter, then the mapping $\Phi$, and hence  $x_+$, both 
depend analytically on the same parameter. The uniformity of the ``o'' from (\ref{phicla}) 
with derivatives follows from the proof of Corollary \ref{cpoc}. 
\end{proof}

\section{Special cases of Theorems \ref{thth} and \ref{tconj}}

\subsection{A simple proof of Theorem \ref{thth} for $G=PSL_2(\rr)$}   
Let $\G$ have rank $M$, thus, $R(\G,G)=G^M$. 
Without loss of generality we assume that $\overline{\rho(\G)}=G$. 
Otherwise, $\rho(\G)$ would be dense in a Lie subgroup of 
dimension at most two, which is solvable, and hence, $\G=\rho(\G)$ cannot be free.  
The group $G=PSL_2(\rr)$ acts conformally on the unit disk $D_1$. 
There is an open subset $U\subset G\setminus1$ formed by 
elliptic transformations, which are conformally conjugate to 
nontrivial rotations. The {\it rotation number} (the rotation angle 
divided by $2\pi$) is a nonvanishing    
analytic function $\R:U\to\rr\slash\zz$. A transformation $f\in U$ has order 
$k$, if and only if $\R(f)=\frac mk(mod\zz)\in\mathbb Q\slash\zz$. 
Let us lift the representation $\rho$ to $SL_2(\rr)$, which is possible 
since $\G$ is free. For every $w\in\G$ the trace $tr(\phi(w))$ is 
an analytic function of $\phi\in R(\G,SL_2(\rr))$. Ellipticity of $g\in SL_2(\rr)$ 
is equivalent to the inequality 
$$-2<tr(g)<2.$$
Fix a $w\in\G$ with $\rho(w)\in U$: then $\rho(w)\neq1$, hence, 
$w\neq1$. 
Let $\lambda^{\pm1}(g)$ denote the eigenvalues of a matrix $g\in SL_2(\rr)$. 
Since $\lambda(g)+\lambda^{-1}(g)=tr(g)$ and, for elliptic elements, $\lambda(g)=
e^{\pm \pi i\R(g)}$ (up to multiplication by $-1$), it follows that either the trace of 
$\phi(w)$ is constant on the 
representation variety to $SL_2(\rr)$, or $\R(\psi(w))\in\mathbb Q\slash\zz$ for 
$\psi$ arbitrarily close to any given $\phi$. In the latter case there exists a sequence 
$\rho_n\to\rho$ such 
that $\R(\rho_n(w))=\frac{m_n}{k_n}(mod\zz)\in\mathbb Q\slash\zz$. Thus, 
$\rho_n(w^{2k_n})=1$  and the representations 
$\rho_n$ are non-injective. The 
former case (the constance of the trace) is impossible by deforming $\phi$ to the 
trivial representation $\phi_0$ for which $tr(\phi_0(w))=2$.

\subsection{On flexibility of representations of arbitrary finitely-generated 
group}
Everywhere in this subsection $\G$ is a finitely-generated (not necessarily free) group, 
$G$ is a real semisimple linear algebraic group. First we briefly  
overview results related to Theorem \ref{tconj} on structural stability of 
representations $\G\to G$. Then we state an open question. 

\begin{definition} An injective representation $\rho\in\rgg$ is called {\it flexible}, if 
it is a limit of non-injective representations. Otherwise, the 
 representation is called {\it structurally-stable} (\cite{sul}, page 243). 
\end{definition}

\begin{example} \label{ex} If $\rho\in\rgg$ is a dense representation  and 
$[\rho]\in\xgg$ is an isolated point, 
then $\rho$ is  structurally-stable: each close representation is conjugate to 
$\rho$, by Proposition \ref{densecenter}. 
Other examples of structurally-stable 
representations include those discrete representations that are not limits 
of nondiscrete ones. The next theorem of D.Sullivan describes the 
structurally-stable finitely-generated torsion-free subgroups in $PSL_2(\cc)$. 
N.Hitchin studied the space of those representations of a surface group 
to $G=PSL_n(\rr)$, $n>2$, whose compositions with the adjoint representation are  
completely 
reducible. He proved that this representation space has three components if $n$ is 
odd, and six components if $n$ is even (\cite{hit}, p.450, theorem B). The 
{\it Teichm\"uller component} is the one 
containing the injective discrete cocompact Fuchsian representations to a given 
subgroup $PSL_2(\rr)\subset G$.  F.Labourie proved that all the representations in the 
Teichm\"uller component  are discrete and injective (as the above ones) 
and structurally stable in this 
component (\cite{lab}, page 53, theorem 1.5).   
\end{example}

\begin{theorem} \label{sull} (\cite{sul}, page 243) An injective representation 
$\rho:\G\to G=PSL_2(\cc)$ of a finitely-generated torsion-free  group $\G$ 
is structurally-stable, if and only if either $[\rho]\in\xgg$ is an isolated 
point, or $\rho(\G)$ is discrete geometrically finite without parabolics.
\end{theorem}

\begin{example} \label{exref} (suggested by the referee). 
There exist injective dense representations that correspond to 
an isolated point in the character variety. For example, let $\G$ be the group of 
$\zz[\sqrt2]$- integer automorphisms of the quadratic form 
$x_1^2+\dots+x_n^2-\sqrt2x_{n+1}^2$, $n\geq4$. The Galois conjugation 
$\sigma:a+b\sqrt2\mapsto a-b\sqrt2$ maps $\G$ onto a dense subgroup of the compact Lie group $O(n+1)$ preserving the quadratic form 
$x_1^2+\dots+x_n^2+\sqrt2x_{n+1}^2$. The conjugation $\sigma$ induces an 
%$\mathbb Q[\sqrt2]$- 
isomorphism (also denoted $\sigma$) 
between the Lie algebras of the Lie groups preserving the above quadratic forms. 
The representation $\sigma:\G\to O(n+1)$ 
is injective and locally rigid: the corresponding point of the character variety is 
isolated, since  $H^1(\G,Ad\circ\sigma)=0$. (See Subsection 2.3 for the definition of 
the first cohomology group.) The latter holds true, since $H^1(\G,Ad)=0$ 
(Garland-Raghunathan rigidity theorem, \cite{garl}, p.284, theorem 0.10) and 
$\sigma$ is induced by an underlying field 
isomorphism, and hence, defines an isomorphism between the above cohomology 
groups. 
\end{example}

\begin{remark} Theorem \ref{tconj} for $G=PSL_2(\cc)$ and $\Gamma$ torsion-free 
 follows from Theorem \ref{sull}. 
\end{remark}
Theorem \ref{tconj} implies the following corollary.

\begin{corollary} Suppose that $\G$ has deficiency at least 2, i.e., it admits 
a finite presentation with $k$ generators and $k-2$ relators. Then every 
injective dense representation of $\G$ to a complex simple linear algebraic  
group $G$  is flexible.
\end{corollary}

\begin{proof} The representation variety $\rgg$ is given by $(k-2)dim G$ 
equations on $k(dim G)$ variables. (Here all the dimensions are over $\cc$.) 
Therefore, $dim\rgg\geq2dim G$ at every point. 
% (similarly to the proof of proposition 3.2.1 in \cite{cs}
Every dense representation $\rho:\G\to G$ lies in the set $S$ from Proposition \ref{densecenter}. 
Hence, $dim_{[\rho]}\xgg\geq dim G$, by Proposition \ref{densecenter}. 
Every complex simple Lie group is simple as a real Lie group.   
Now, the assertion of the corollary follows from Theorem \ref{tconj}.
\end{proof} 

\begin{remark} (due to the referee). Theorem \ref{tconj} fails for dense 
representations to general semisimple 
linear algebraic groups (not necessarily without connected normal Lie subgroups). 
Namely, let $\G$ be an arithmetic lattice in $SO(n,1)$, $n\geq4$, from 
Example \ref{exref}. Then for every $j\geq2$ the group $\G$ contains a finite index 
subgroup which admits an epimorphism to 
$F_j$ (the free group with $j$ generators). This follows from 
the proof of theorem 3.1 in \cite{lub}, p.75.  Take $j=3$, and let $\G$ denote 
the latter subgroup. Let $\rho_1=\sigma:\G\to O(n+1)$ denote the dense representation 
from Example \ref{exref}. Let 
$\rho_2:\G\to F_3\to SO(n+2)$ be a dense representation. 
Set $\rho_0=(\rho_1,\rho_2):
\G\to G=O(n+1)\times SO(n+2)$. The representation $\rho_0$ is dense, since 
so are $\rho_1$ and $\rho_2$ and 
there are no nontrivial closed Lie subgroups in $G$ projected epimorphically onto 
each factor of $G$: the groups $O(n+1)$ and $SO(n+2)$ are non-isomorphic 
simple groups. The representation $\rho_0$ admits nontrivial deformations  
coming from deforming the homomorphism $F_3\to SO(n+2)$. Therefore, 
$dim_{[\rho_0]}\xgg\geq 2 dim(SO(n+2))>dim G$, since this is the case for 
dense representations of $F_3$ to $SO(n+2)$. However the representation $\rho_0$ is 
structurally stable, as is $\rho_1$. Thus, all the representations close to $\rho_0$ are 
injective. 
\end{remark}

{\bf Question 3.} {\it Let $G$ be a simple linear algebraic group. 
Do there exist a finitely-generated group $\G$ and a structurally-stable dense 
 representation 
$\rho:\G\to G$ such that $[\rho]\in\xgg$ is not an isolated point?}

For $G=PSL_2(\cc)$ and $\G$ torsion-free, Theorem \ref{sull} gives a negative answer.

\section{Proof of Theorem \ref{tconj}: plan and Main Technical Lemma}
Here we sketch the proof of Theorem \ref{tconj} in the general case and finish it 
for groups with proximal elements. We complete it for groups without proximal elements 
 in Section 5.

\subsection{Motivation and the plan of the proof of Theorem \ref{tconj}} 
Each representation 
satisfying the conditions of Theorem \ref{tconj} is a limit of representations 
$\rho\in S_{reg}\subset R_{reg}(\G,G)$, see Proposition \ref{densecenter}, such that  
\begin{equation}\overline{\rho(\G)}\supset G_0, \ [\rho]\in X_{reg}(\G,G), \  \operatorname{rank}(\pi)= dim_{[\rho]}\xgg\geq dim G \text{ at } \rho,
\label{zv}\end{equation}
recall that $\pi$ denotes the projection $\rgg\to\xgg$. 
It suffices to prove Theorem \ref{tconj} for every injective 
representation $\rho$ satisfying (\ref{zv}). Its statement for every limit of 
injective representations (\ref{zv}) then follows by diagonal sequence argument. 
(If the approximants (\ref{zv}) are already non-injective, then the flexibility 
of $\rho$ follows immediately.) 

\begin{definition} \label{dimmers} Let $\G$ be a finitely-generated group, 
$G$ be a semisimple linear algebraic group, $n=dim G$.                                                                                              A family of 
representations $\rho_u:\G\to G$ depending smoothly on the parameter $u\in\rr^n$ 
is called $(G,X)$- {\it immersive}, if $\xgg$ is smooth at the points 
$[\rho_u]$, and 
\begin{equation}\text{the mapping } \rr^n\to\xgg: \ u\mapsto[\rho_u], 
\text{ is an immersion to } X_{reg}(\G,G). 
\label{immers}\end{equation}
\end{definition}

In this section, whenever the contrary is not specified, we 
make the following assumption.

\begin{assumption}\label{*} 
The group $\G$ is finitely-generated, $G$ is a semisimple linear 
algebraic group,  $n=dim G$, $\rho=\rho_0:\G\to G$ is an 
injective  representation satisfying (\ref{zv}), the Lie subgroup
$\overline{\rho(\G)}$  has no nontrivial 
connected normal Lie subgroups,  $\rho_u$ is a  
$(G,X)$- immersive deformation of $\rho$,  which exists by (\ref{zv}). 
\end{assumption}

We construct sequences $w_r\in\G$ and $u_r\in\rr^n$ such that 
\begin{equation}u_r\to0, \text{ as } r\to\infty, \ \rho_{u_r}(w_r)=1, \ 
\rho_u(w_r)\not\equiv1.\label{wklook}
\end{equation}
Then $w_r\neq1$ (by the last inequality), the 
representations $\rho_{u_r}$ are non-injective, and hence, 
$\rho=\lim_{r\to\infty}\rho_{u_r}$ is flexible.  
This proves Theorem \ref{tconj}. 

First let us motivate the construction of $w_r$ and $u_r$. One would naturally 
look for $w_r$ among those sequences for which $\rho(w_r)\to1$, as $r\to\infty$. 
This is the case, e.g, if we 
take two elements $h_1,h_2\in\G$ with $\rho(h_1)$, $\rho(h_2)$ 
close to the identity and set $w_r=[\dots[h_1,h_2],\dots,h_2]$: the $r$- th iterated 
commutator. To show that $\rho_{u_r}(w_r)=1$ for some $u_r\to0$, one has to 
prove a lower bound of derivative implying, e.g., that for every $v\in T_0\rr^n$ 
with $|v|=1$ the derivative 
$\frac{d\rho_u(w_r)}{dv}$ asymptotically dominates $dist(\rho(w_r),1)$, 
as $r\to\infty$. But for the above  $w_r$ both latter quantities tend to zero 
exponentially fast. We modify $w_r$ using the following observation, which 
holds under some genericity assumptions on the families of elements  $\rho_u(h_1)$ and 
$\rho_u(h_2)$. Fix a $\Delta>0$ small enough and a sequence 
$m_r\in\nn$ such that $\Delta<dist(\rho(w_r^{m_r}),1)<2\Delta$ for all $r$. 
Then for  a generic $v$, the derivative 
$\frac{d\rho_u(w_r^{m_r})}{dv}$  grows linearly in $r$, as $r\to\infty$.  

In what follows we construct appropriate 

- elements $g_1,\dots,g_n,h,w\in\G$
and define recurrently the iterated commutators 
\begin{equation}w_{i,0}=h,\ w_{i,r}=g_iw_{i,r-1}g_i^{-1}w_{i,r-1}^{-1},
\label{defwik}\end{equation}

- vector $l=(l_1,\dots,l_n)\in\zz^n$, 

- sequences $k_r\in\mathbb N$, $k_r\to+\infty$, as $r\to\infty$, 
and $m_{jr}\in\zz$, $j=1,\dots,n$, and set 
\begin{equation}\omega_r=w_{1,k_r+l_1}^{m_{1r}}\dots w_{n,k_r+l_n}^{m_{nr}}, \ 
w_r=w^{-1}\omega_r.\label{owr}\end{equation}
We show that
\begin{equation}\rho_{k_r^{-1}\wt u}(\omega_r)\to\Psi(\wt u), \text{ as } r\to\infty; \ 
\Psi:\rr^n\to G_0 \ \text{is a local diffeomorphism at} \ 0,\label{convpsi}\end{equation}
the above convergence is uniform with derivatives on compact subsets in $\rr^n$. 
Theorem \ref{tconj} will be easily deduced from (\ref{convpsi})  
at the end of the subsection.  

First we make a preliminary construction of the above $g_j$, $h$, $k_r$ and 
$m_{jr}=m_{jr}(\var)$, the latter 
depending on an additional small parameter $\var>0$,  so that the 
mapping sequence (\ref{convpsi}) converges for every $l$. 
The limit mapping $\Psi$ depends only on $g_j$, $h$, $l$ and $\var$. 
Afterwards we show (and this is the main technical part of the proof of 
Theorem \ref{tconj}) that one can adjust the latter so that $\Psi$ is 
a local diffeomorphism at 0 (Main Technical Lemma \ref{lsi}, 
 Corollary \ref{sind} and Lemma \ref{lem1}). Lemma \ref{lsi}  is proved 
 in Subsection 4.3. Lemma \ref{lem1} is proved in Subsection 4.2 for groups 
with proximal elements and in Section 5 for groups without proximal elements.  

Let $\hat\Pi=\Pi_{\rr,1}\sqcup\Pi_{\cc,1}\subset G$,  $\Pi'\subset\hat\Pi\times G$, $s(g)$ and 
$v_g$ be the same, as in (\ref{hatpi}) and Proposition \ref{pcomm}. We construct 
elements $g_j$ and $h$ so that 
\begin{equation}(\rho_0(g_j), \rho_0(h))\in\Pi'\text{ for all } j.\label{AA}\end{equation}
Set 
\begin{equation}s_j(u)=s(\rho_u(g_j))\in\cc, \ \wt\nu_j(u)=v_{\rho_u(g_j)}
(\rho_u(h))\in\gg, \ \nu_j=\wt\nu_j(0).\label{nui}
\end{equation}
For a natural sequence $k_r\to+\infty$ and an 
$\var>0$, for every $r\in\mathbb N$ and $j=1,\dots,n$ set 
\begin{equation} m_{jr}=[\var|s_j|^{-k_r}(0)]\in\mathbb N.\label{mjr}\end{equation}

{\bf Preliminary construction of $\omega_r$: case when $G$ has proximal elements.} 
We set 
$$k_r=2r, \ l_j=0.$$
We choose $g_j$ and $h$ so that (\ref{AA}) holds and each 
$\rho_0(g_j)$ is 1-proximal and take $m_{jr}$  as in (\ref{mjr}).

\medskip

{\bf Claim 0.} {\it In the above assumptions,  
\begin{equation}\rho_{k_r^{-1}\wt u}(\omega_r)\to
\Psi(\wt u)=
\exp(\var e^{(d\ln s_1(0))\wt u}\nu_1)\dots
\exp(\var e^{(d\ln s_n(0))\wt u}\nu_n), \text{ as } r\to\infty,
\label{rholim}\end{equation}
 uniformly with derivatives on compact sets in $\rr^n$.}

\medskip

The claim and a more general Proposition \ref{propconvnew} are proved below. 

{\bf Preliminary construction of $\omega_r$: case when $G$ has no proximal 
elements.} We 
 choose $g_j$ and $h$ so that (\ref{AA}) holds and each $\rho_0(g_j)$ is 
 $\cc$-1-proximal. Set 
\begin{equation} \zeta_j=\arg s_j(0).\label{nunu}\end{equation}

\begin{proposition} \label{pkr} For every real vector 
$\zeta=(\zeta_1,\dots,\zeta_n)\in\rr^n$ 
there exists a sequence of numbers $k_r\in\mathbb N$, $k_r\to+\infty$, as $r\to\infty$, 
such that 
\begin{equation}k_r\zeta_j\to0 (mod2\pi), \ \text{as} \ r\to\infty, \ \text{for every} \ 
j=1,\dots,n.\label{zetaj}\end{equation}
\end{proposition}
\begin{proof} Consider $\zeta$ as an element of the torus $\tt^n=\rr^n\slash2\pi\zz^n$. 
The subgroup $<\zeta>\subset\tt^n$ either is discrete, or accumulates to 0. In both 
cases there exists a sequence  $k_r\in\nn$, $k_r\to\infty$, such that 
$k_r\zeta\to0$ in $\tt^n$ (which is equivalent to (\ref{zetaj})). In the second 
case this follows from definition. In the first case the group $<\zeta>$ is finite cyclic 
by compactness. Let  $m$ denote its order, $k_r=rm$. Then $k_r\zeta=0$ for all 
$r\in\nn$. This proves Proposition \ref{pkr}.
\end{proof} 

We take $k_r$ the same, as in (\ref{zetaj}), with $\zeta_j$ as in  (\ref{nunu}), 
and $m_{jr}$ the same, as in (\ref{mjr}), 
and an arbitrary  vector $l\in\zz^n$. The preliminary construction of 
$\omega_r$ is complete. 

\begin{proposition} \label{propconvnew}
Let $\G$ be a finitely-generated group, $G$ be a Lie group, 
$n=dim G$,
$\rho_u:\G\to G$ be a family of representations depending on a parameter $u\in\rr^n$. 
Let $g_1,\dots,g_n,h\in\G$ be such that $(\rho_0(g_j),\rho_0(h))\in\Pi'$ for all 
$j$. (Here any $\rho_0(g_j)$ may be either 1-proximal, or $\cc$-1-proximal.)  
Let $s_j$, $\nu_j$ and $\zeta_j$ be the same, as in (\ref{nui}) and 
(\ref{nunu}). Let $k_r\to+\infty$ be a natural sequence satisfying (\ref{zetaj}). Let 
$\var>0$, and let $m_{jr}$ be the same, as in (\ref{mjr}). Let 
$l=(l_1,\dots,l_n)\in\zz^n$, and let 
$\omega_r\in\G$ be the corresponding sequence of elements (\ref{owr}). Then 
\begin{equation}\rho_{k_r^{-1}\wt u}(\omega_r)\to\Psi(\wt u)=
\exp(\var s_1^{l_1}(0)e^{(d\ln s_1(0))\wt u}\nu_1)\dots
\exp(\var s_n^{l_n}(0)e^{(d\ln s_n(0))\wt u}\nu_n),\label{defpsi}\end{equation}
as $r\to\infty$, uniformly with derivatives on compact sets in $\rr^n$. 
In the case, when the element $\rho_0(g_j)$ is $\cc$-1-proximal, the multiplication of the 
vector $\nu_j$ in the dominating plane  $L(\rho_0(g_j))$ by complex numbers is defined 
in terms of the $s(\rho_0(g_j))$- complex structure. 
\end{proposition}

\begin{proof} One has $m_{jr}s_j^{k_r+l_j}(0)\to\var s_j^{l_j}(0)$, as $r\to\infty$, 
by (\ref{mjr}) and (\ref{zetaj}). Therefore,
$$\rho_{k_r^{-1}\wt u}(w_{j,k_r+l_j}^{m_{jr}})\to\exp(\var s_j^{l_j}(0)
e^{(d\ln s_j(0))\wt u}\nu_j),$$
 by Corollary \ref{convexp} applied to $g(u)=\rho_u(g_j)$, $h(u)=\rho_u(h)$, 
 $m_r=m_{jr}$, $k_r$ being replaced by $k_r+l_j=k_r(1+o(1))$ and $\var$ being replaced by 
 $\var s_j^{l_j}(0)$. This implies the proposition.
 \end{proof}
 
 \begin{proof} {\bf of Claim 0.} One has $\zeta_j=\arg s_j(0)\in\pi\zz$: $s_j(0)\in\rr$, 
 since $\rho_0(g_j)$ are 1-proximal. Thus, $k_r\zeta_j=2r\zeta_j=0(mod2\pi)$.  
 Now Proposition \ref{propconvnew} with $l_j=0$ implies the claim. 
 \end{proof}

\begin{lemma} \label{lsi} {\bf (Main Technical Lemma).} 
Let $\G$, $G$, $n$, $\rho_u$ be the same, as in the above Assumption \ref{*}. 
 Let $U\subset \overline{\rho(\G)}$ be an arbitrary open subset,  
and let $\sigma:U\to\mathbb R$ be a smooth locally nonconstant function. 
Then there exist $n$ elements $g_1,\dots,g_n\in\G$ with $\rho_0(g_j)\in U$ so that the 
function  $u\mapsto(s_1(u),\dots,s_n(u))$, defined by 
$s_j(u)=\sigma(\rho_u(g_j))$, is a local diffeomorphism at 0. 
Moreover, for any given $A_1,\dots,A_n\in U$ 
one can chose $g_1,\dots,g_n$ so that in addition, the elements 
$\rho_0(g_j)\in U$ are  arbitrarily close to $A_j$.
\end{lemma}

\begin{corollary} \label{sind} Let $\G$, $G$, $n$, $\rho_u$ be  as in 
Assumption \ref{*},  $\hat\Pi$  be the same, as in (\ref{hatpi}). Set
$$\sigma(g)=s(g) \text{ if } g\in\Pi_{\rr,1}, \ \sigma(g)=\arg s(g) \text{ if } g\in\pc1.$$
Let $\Pi^0\subset G_0$ denote the union of those connected components of $\hat\Pi$ 
that intersect at least one one-parametric subgroup in $G$. 
 Then  there exist $g_1,\dots,g_n\in\G$ such that 
$\rho_0(g_j)\in\Pi^0$ and the function 
$u\mapsto(s_1(u),\dots,s_n(u))$ defined by 
$s_j(u)=\sigma(\rho_u(g_j))$ is a local diffeomorphism at 0.
Moreover, for any given $A_1,\dots,A_n\in\Pi^0$ one can choose the above 
$g_1,\dots,g_n$ so that  $\rho_0(g_j)$ are arbitrarily close to $A_j$. 
\end{corollary}

\begin{proof}  The function 
$\sigma(g)$ is analytic and single-valued on $\hat\Pi$ (Proposition \ref{p3.4} 
and Remark \ref{sneq0}). It is locally 
nonconstant on $\Pi^0$, since it is nonconstant along each one-parametric 
subgroup $H$ intersecting $\hat\Pi$. 
Indeed,  the  eigenvalue $s(g)+1$ of the adjoint $Ad_g$ yields a 
nontrivial homomorphism $H\to\cc^*$, since $s(g)\neq0$ for $g\in\hat\Pi$. 
Hence, $s\not\equiv const$ along $H$ and $\sigma=s\not\equiv const$ along the 
intersection $H\cap\Pi_{\rr,1}$ (if the latter is non-empty). In the case, when  
$H\cap\Pi_{\cc,1}\neq\emptyset$,  
one has $\sigma(g)=\arg s(g)$ and $s(g)\notin\rr$ on this intersection 
($\cc$-1-proximality). Therefore, 
in this case the image of the above homomorphism is not contained in the real line, 
thus $\sigma\not\equiv const$ along $H$. 
This together with Lemma \ref{lsi} proves Corollary \ref{sind}.  
\end{proof}

\begin{lemma} \label{lem1} Let $\G$, $G$, $n$, $\rho_u$ be as in 
Assumption \ref{*}. Then there exist elements $g_1,\dots,g_n,h\in\G$ satisfying (\ref{AA}) 
and a vector $l\in\zz^n$ such that  for every $\var>0$ small enough, the 
corresponding mapping $\Psi:\rr^n\to G_0$ from (\ref{defpsi}) is a local diffeomorphism 
at 0. In the case of group $G$ with proximal elements one can take $l=0$ for 
appropriate $g_j$ and $h$. 
\end{lemma}

\def\odd{\overline D_{\delta}}

\begin{proof} {\bf of Theorem \ref{tconj} modulo  Lemmas \ref{lsi} and \ref{lem1}.} 
Let $g_1,\dots,g_n$, $h$, $l$, $\var$ be 
as in Lemma \ref{lem1}, $s_j(u)=s(\rho_u(g_j))$, $\zeta_j=\arg s_j(0)$. 
Let $k_r\to+\infty$ be a natural sequence satisfying (\ref{zetaj}), 
$m_{jr}$ be the 
numbers from (\ref{mjr}). Let $\omega_r$ be the corresponding iterated commutator 
power product (\ref{owr}), $\Psi$ be the mapping from (\ref{defpsi}). Fix a 
$\delta>0$ and an element $w\in\G$ such that $\Psi:\overline D_{\delta}\to
\Psi(\overline D_{\delta})$ is a diffeomorphism and $\rho_0(w)
\in\Psi(D_{\delta})$. Set $w_r=w^{-1}\omega_r$. Then  
$\rho_{k_r^{-1}\wt u}(w_r)\to\psi(\wt u)=\rho_0(w^{-1})\Psi(\wt u)$, as 
$r\to\infty$. Furthermore, 
$\psi:\overline D_{\delta}\to\psi(\overline D_{\delta})\subset G$ is a 
diffeomorphism, $1\in\psi(D_{\delta})$, and the above 
convergence is uniform with derivatives on compact subsets in $\rr^n$, by construction 
and Proposition \ref{propconvnew}.  Thus, for each $r$ large enough  
  $\rho_{k_r^{-1}\wt u_r}(w_r)=1$ for some $\wt u_r\in D_{\delta}$. 
 Set $u_r=k_r^{-1}\wt u_r$.  Then $\lim_{r\to\infty}u_r=0$ and $\rho_{u_r}(w_r)=1$. 
 One has $\rho_u(w_r)\not\equiv1$, since the mappings $\wt u\mapsto
 \rho_{k_r^{-1}\wt u}(w_r)$ are  diffeomorphisms on $\overline D_{\delta}$ for large $r$, 
as is their limit $\psi(\wt u)$. Thus, the elements 
$w_r$ satisfy (\ref{wklook}).  This proves Theorem \ref{tconj}.  
\end{proof}

The above discussion implies the next theorem, which summarizes the technical 
results of Sections 4 and 5. 
It will be used in the proof of Theorems \ref{cdiap0} and \ref{diap} (Section 7).

\begin{theorem} \label{tuple} Let $\G$, $G$, $n$, $\rho_u$ be the same, as in 
Assumption \ref{*}. Then 
there exist sequences of elements  $w_r\in\G$, numbers $k_r\in\nn$, $k_r\to+\infty$, as 
$r\to\infty$,  a real-analytic 
mapping $\psi:\rr^n\to G_0$  and a $\delta>0$ such that  
\begin{equation} \rho_{k_r^{-1}\wt u}(w_r)\to\psi(\wt u) \text{ uniformly with 
derivatives on 
compact subsets in } \rr^n, \text{ as } r\to\infty,\label{wrto}\end{equation}
\begin{equation}\psi:\odd\to\psi(\odd)\subset G_0 \text{ is a diffeomorphism, and } 
1\in\psi(D_{\delta}).\label{psidif}\end{equation}
\end{theorem} 
\begin{definition} \label{dtuple} Under the assumptions of Theorem \ref{tuple} the 
tuple 
$(\{ w_r\}, \{ k_r\}, \psi,\delta)$ is called the {\it converging tuple} 
associated to the given $(G,X)$- immersive representation family $\rho_u$. 
\end{definition}

\def\psie{\Psi_{\var}}
\def\psije{\Psi_{j,\var}}
\def\wtpsie{\wt\Psi_{\var}}
\def\wtpsije{\wt\Psi_{j,\var}}

\subsection{The diffeomorphicity: case of proximal elements}
In this section we assume that $G$ contains proximal elements and prove 
Lemma \ref{lem1} modulo the Main Technical Lemma \ref{lsi}. 
To do that, we consider the mapping $\Psi=\Psi_{\var}(\wt u)$ corresponding to fixed 
$g_j,h\in\G$ and variable  $\var$. Recall that in the case under consideration $l=0$ and 
$s_j(u)$ are real-valued functions on a neighborhood of zero, $0<|s_j(u)|<1$. Set 
$$\sigma_j=d\ln s_j(0):T_0\rr^n=\rr^n\to\rr, \  
\psije(\wt u)=\exp(\var e^{\sigma_j(\wt u)}\nu_j): \ 
\psie(\wt u)=\Psi_{1,\var}(\wt u)\dots\Psi_{n,\var}(\wt u);$$ 
\begin{equation} \wtpsie(\wt u)=\psie(\wt u)(\psie(0))^{-1}.\label{sifi}
\end{equation}

\begin{proposition} \label{asfor} Let $\sigma_1,\dots,\sigma_n:T_0\rr^n\to\rr$ 
be arbitrary linear forms, $\nu_1,\dots,\nu_n\in\gg$. Let $\psie$, $\wtpsie$ 
be the corresponding mappings (\ref{sifi}). Then 
\begin{equation}\wtpsie'(0)=\var\sum_{j=1}^n\nu_j\sigma_j+O(\var^2), \text{ as }
\var\to0.\label{wtpsi}\end{equation}
Here $\wtpsie'$ is the derivative of the mapping $\wtpsie$. 
\end{proposition}

\begin{proposition} \label{nond} Under Assumption \ref{*}  
there exist $g_1,\dots,g_n,h\in\G$ satisfying (\ref{AA}) and such that the 
elements $\rho_0(g_j)\in G$ are 1-proximal and the linear operator 
$\sum_{j=1}^n\nu_j\sigma_j:T_0\rr^n\to\gg$ 
is invertible. Here $\sigma_j=d\ln s_j(0)$, and $\nu_j$ are the same, as in 
(\ref{nui}). 
\end{proposition}

\begin{proof} {\bf of Lemma \ref{lem1}.} The $g_j$ and $h$ from Proposition 
\ref{nond} and $l=0$ satisfy the statement of the lemma: if $\var$ is small 
enough, then the derivative 
$\wtpsie'(0)$ (and hence, $\psie'(0)$) is invertible, as is the operator 
$\sum_{j=1}^n\nu_j\sigma_j$ in (\ref{wtpsi}), 
by (\ref{wtpsi}) and Proposition \ref{nond}. 
\end{proof}

\begin{proof} {\bf of Proposition \ref{asfor}.} Let $F$ be the free group with 
generators $x_1,\dots,x_n$. Consider the family of representations 
$\phi_{\wt u,\var}:F\to G$ defined by 
$$\phi_{\wt u,\var}(x_j)=\psije(\wt u).$$
Then 
$$\phi_{\wt u,\var}(x_1\dots x_n)=\psie(\wt u).$$
For every $\var>0$ and $v\in T_0\rr^n$ the vector function $c_{v,\var}:F\to\gg$:  
$c_{v,\var}(x)=\frac{d\phi_{\wt u,\var}(x)}{dv}
(\phi_{0,\var}(x))^{-1}$ is a cocycle with respect to the homomorphism  
$Ad\circ\phi_{0,\var}$, see Subsection 2.3. Set 
$\wtpsije(\wt u)=\psije(\wt u)(\psije(0))^{-1}$. Then 
$c_{v,\var}(x_j)=\frac{d\wtpsije(\wt u)}{dv}$, $c_{v,\var}(x_1\dots x_n)=
\frac{d\wtpsie(\wt u)}{dv}$. Therefore, 
\begin{equation}\wtpsie'(0)=\wt\Psi_{1,\var}'(0)+
Ad_{\Psi_{1,\var}(0)}\wt\Psi_{2,\var}'(0)+\dots+
Ad_{\Psi_{1,\var}(0)\dots\Psi_{n-1,\var}(0)}\wt\Psi_{n,\var}'(0), 
\label{cocpr}\end{equation}
by the generalized cocycle identity (\ref{cocn}) applied to $c_{v,\var}$ for all $v$. 
One has 
$\psije(0)=\exp(\var\nu_j)$, 
$$\wtpsije(\wt u)=\exp(\var e^{\sigma_j(\wt u)}\nu_j)
\exp(-\var\nu_j)=\exp(\var\nu_j(e^{\sigma_j(\wt u)}-1)).$$ 
Hence,  
\begin{equation}\wtpsije'(0)=\var\nu_j\sigma_j,\label{derj}\end{equation}
\begin{equation}
Ad_{\psije(0)}=Id+O(\var), \ Ad_{\Psi_{1,\var}(0)\dots\Psi_{k,\var}(0)}=
Id+O(\var), \text{ as } \var\to0,\label{adps}\end{equation} 
since $dist(\psije(0),1)=O(\var)$. Substituting (\ref{derj}) and (\ref{adps}) 
to (\ref{cocpr}) yields (\ref{wtpsi}).
\end{proof}

\begin{proof} {\bf of Proposition \ref{nond}.} Let $\Pi^0$ be the same, as 
in Corollary \ref{sind}. The intersection $\Pi^0\cap\Pi_{\rr,1}$ is open, nonempty 
(Theorem \ref{1prpr}) and invariant under conjugations. 
Fix an $A_1\in\Pi^0\cap\Pi_{\rr,1}$. Fix some 
$H_2,\dots,H_n\in \overline{\rho(\G)}$, set $H_1=Id$, $A_j=H_jA_1H_j^{-1}$, such that the lines 
$L(A_j)=Ad_{H_j}L(A_1)\subset\gg$ are linearly independent (the irreducibility 
of the adjoint of the Lie group $\overline{\rho(\G)}$, which has no nontrivial connected 
normal Lie subgroups). Then $A_j\in\Pi^0\cap\Pi_{\rr,1}$. Fix elements $g_j\in\G$ from 
Corollary \ref{sind} (already proved modulo Lemma \ref{lsi}) so that 
$A_j'=\rho_0(g_j)$ are close enough to $A_j$: 
the lines $L(A_j')$ should be also linearly independent. Fix a $h\in\G$ so that 
$(A_j',\rho_0(h))\in\Pi'$ and $\nu_j=v_{A_j'}(\rho_0(h))\neq0$ for all $j$. 
This $h$ exists, since $dv_{A_j'}\neq0$, see 
(\ref{vg}), and $\rho_0(\G)$ is dense in $G_0$. The vectors $\nu_j$ are linearly 
independent, as are the ambient lines $L(A_j')$. The forms 
$\sigma_j=d\ln s_j(0)$ are also independent, by Corollary \ref{sind}. Thus, the operator 
$\sum_{j=1}^n\nu_j\sigma_j$ is invertible. 
Proposition \ref{nond} is proved.  Lemma \ref{lem1} and Theorem 
\ref{tconj} are proved  for Lie groups with proximal elements 
modulo the Main Technical Lemma \ref{lsi}. 
\end{proof}

\def\la{\lambda}

\def\ton{T_0\rr^n\setminus0}

\subsection{Independent eigenvalues. Proof of Main Technical Lemma \ref{lsi}} 
Let $\G\subset G$ be a  subgroup whose closure 
$\overline\G\subset G$ is an open Lie subgroup. In what follows we denote 
$G$ the latter open Lie subgroup. Let $\rho_0:\G\to\G\subset \overline\G$
 be  the identity representation. It suffices to prove Lemma \ref{lsi} for deformations of 
the identity representation $\rho_0$. In this subsection we fix a right-invariant 
Riemannian metric on $G$. For every $v\in T_0\rr^n\setminus0$ we set 
$$c_v:\G\to\gg \text{ the 1- cocycle associated to the derivative } \frac{d\rho_u}{dv}
\in T_{\rho_0}\rgg, \text{ see (\ref{corcoc})}.$$

\medskip

{\bf Claim 1.} {\it For every $v\in T_0\rr^n\setminus0$ the cocycle $c_v$ is 
not Lipschitz at 1 in the metric on $\G=\rho_0(\G)\subset G$ induced from $G$.}

\medskip

\begin{proof} Fix a $v\in T_0\rr^n\setminus0$. Assume to the contrary that $c_v$ is 
Lipschitz at 1. Then $c_v$ extends to a locally Lipschitz 
cocycle on $G$ with coefficients in $Ad$ (Proposition \ref{p}). 
Therefore, $c_v$ is a coboundary 
(Corollary \ref{corp}). Hence, the derivative $\frac{d\rho_u}{dv}$ is tangent to the 
$G$- orbit of $\rho_0$ (Remark \ref{cohom}). Thus, 
$\frac{d[\rho_u]}{dv}=0$, - a contradiction to the 
$(G,X)$- immersivity, see (\ref{immers}). 
\end{proof}

\begin{proposition} \label{plimline} 
Let $G$ be a Lie group with irreducible adjoint representation, $\G\subset G$ be a dense 
subgroup. Let $c\in Z^1(\G,Ad)$ be a cocycle  
that is not Lipschitz at the identity in the metric on $\G$ induced from $G$. 
Then for every one-dimensional linear subspace 
$\Lambda\subset\gg$ there exists a sequence of elements $b_k\in\G$ such that 
\begin{equation} c(b_k)\neq0 \text{ for all } k; \ 
 b_k\to1, \ \frac{c(b_k)}{dist(b_k,1)}\to\infty, \ \rr c(b_k)\to\Lambda,
\text{ as } k\to\infty.\label{limline}\end{equation}
Here the lines $\rr c(b_k)$ converge to $\Lambda$ as points of the projectivization of 
$\gg$. 
\end{proposition}

Below we prove Proposition \ref{plimline} and deduce the following proposition from it. 

\begin{proposition} \label{climline} Let $G$, $\G$, $c$ be as in Proposition \ref{plimline}. 
Consider the vector field $\nu$ on $\G$, 
$\nu(\gamma)=c(\gamma)\gamma\in T_{\gamma}G$, $\gamma\in\G$,  corresponding 
to the cocycle $c$ (denoted $v_a$ in (\ref{corcoc})). Let $U\subset G$ be an 
open subset, and let 
$\sigma:U\to\rr$ be a $C^1$- smooth locally-nonconstant function. Then each $g\in U$ is 
 a limit of a sequence $h_k\in\G$ such that $\frac{d\sigma}{d\nu(h_k)}\neq0$ for all $k$. 
\end{proposition}

\begin{proof} {\bf of Main Technical Lemma \ref{lsi} modulo 
 Proposition \ref{climline}.} For every $v\in T_0\rr^n\setminus0$ the vector field 
 $$\nu_v(\gamma)=\frac{d\rho_u(\gamma)}{dv}=c_v(\gamma)\gamma\in T_{\gamma}G, \ 
  \gamma\in\G,$$
 satisfies the conditions of Proposition \ref{climline} with 
 $G=\overline{\rho_0(\Gamma)}$,  since the cocycle $c_v$ is not 
 Lipschitz at 1 (Claim 1) and $Ad_G$ is irreducible 
 (Remark \ref{ad-irr}). 
 Fix a collection $A_1,\dots,A_n\in U$,  a $\var>0$ and a vector 
 $v^1\in T_0\rr^n\setminus0$. There exists an element $g_1\in U\cap\G$ 
$\var$- close to $A_1$ such that $\frac{d\sigma}{d\nu_{v^1}(g_1)}\neq0$. 
This follows from Proposition \ref{climline} 
applied to $g=A_1$: $g_1$ is  an $\var$- approximant $h_k$ of $g$.  One has 
$$s_1(u)=\sigma(\rho_u(g_1)), \ 
\frac{ds_1}{dv^1}=\frac{d\sigma}{d\nu_{v^1}(g_1)}\neq0.$$ 
Hence, $codim Ker(ds_1(0))=1$. 
Fix a  $v^2\in Ker(ds_1(0))\setminus0$. We construct similarly a  
$g_2\in U\cap\G$ $\var$- close to $A_2$, such that 
$\frac{ds_2}{dv^2}=\frac{d\sigma}{d\nu_{v^2}(g_2)}\neq0$.  One has:  
$codim (Ker(ds_1(0))\cap Ker(ds_2(0)))=2$. 
Now fix a $v^3\in Ker(ds_1(0),ds_2(0))\setminus0$ and construct $g_3$ etc. 
The elements $g_1,\dots,g_n\in\G$ thus constructed are the ones we are looking for: 
$codim Ker(ds_1(0),\dots,ds_n(0))=n$.  Thus, this kernel vanishes,  
$(s_1(u),\dots,s_n(u))$ is a local diffeomorphism at 0. Lemma \ref{lsi} is 
proved.
\end{proof}

\def\clr{\mathcal L_r}
 
 \begin{proof} {\bf of Proposition \ref{plimline}.} For every one-dimensional 
 linear subspace $\Lambda\subset\gg$ define the following set of sequences and 
 union of lines:
 $$\mathcal P_{\Lambda}=\{\{b_k\}_{k\in\nn}\subset\G \ | \text{ statement 
 (\ref{limline}) holds}\},$$
 \begin{equation}\hat L=\text{ the union of the lines } \Lambda\subset\gg \ 
 \text{for which} \ 
 \mathcal P_{\Lambda}\neq\emptyset.\label{defhl}\end{equation}
 The set $\hat L$ is nonempty. Indeed, there exists a sequence $\{ b_k\}_{k\in
 \mathbb N}\subset\G\subset G$ such that $b_k\to1$ and $dist(b_k,1)=o(c(b_k))$, 
 as $k\to\infty$, since the cocycle is not Lipschitz at 1. Passing to a 
 subsequence, one can achieve that $c(b_k)\neq0$ for all $k$ and the lines 
 $\rr c(b_k)$ converge to a line $\Lambda$. Then $\{ b_k\}_{k\in
 \mathbb N}\in\mathcal P_{\Lambda}$ and $\Lambda\subset\hat L$. The set $\hat L$ 
  is closed by definition and the diagonal sequence argument. We claim that this set is 
  dense in $\gg$ (this will prove Proposition \ref{plimline}). To prove density we fix  
some $\Lambda$ with $\mathcal P_{\Lambda}\neq\emptyset$, a 
$\{b_k\}_{k\in\mathbb N}\in\mathcal P_{\Lambda}$, a $\xi\in\Lambda\setminus0$  
and $g_1=id, g_2,\dots,g_n\in \G$ such that 
\begin{equation}
\text{ the vectors } \xi=Ad_{g_1}\xi, Ad_{g_2}\xi,\dots,Ad_{g_n}\xi \text{ form a basis in }
\gg.\label{adbase}\end{equation}
These elements $g_i$ exist by the irreducibility of  $Ad$ and the density of 
$\G$.  For any $r=(r_1,\dots,r_n)\in\zz^n\setminus0$ 
we consider the following auxiliary linear operator, lines and sequence of elements of 
$\G$:
\begin{equation}\clr=\sum_{j=1}^nr_jAd_{g_j}:\gg\to\gg, \ \Lambda_r=\clr(\Lambda), \ 
b_{r,k}=g_1b_k^{r_1}g_1^{-1}\dots g_nb_k^{r_n}g_n^{-1}.\label{defclr}
\end{equation}
Note that for any $r\in\zz^n\setminus0$ the operator $\clr:\Lambda\to\Lambda_r$ is 
invertible, by (\ref{adbase}), and hence, the space $\Lambda_r$ is a 
line. The union of all the lines $\Lambda_r$, $r\in\zz^n\setminus0$, is dense in $\gg$, 
by (\ref{adbase}).

\medskip

{\bf Claim 2.} {\it For every $r\in\zz^n\setminus0$ one has 
$\{b_{r,k}\}_{k\in\mathbb N}\in\mathcal P_{\Lambda_r}$: thus, $\Lambda_r
\subset\hat L$.}

\medskip

The claim (proved below) together with the previous statement implies that 
the closed set $\hat L$ contains a dense subset $\cup_r\Lambda_r\subset\gg$. 
Hence, $\hat L=\gg$. This proves Proposition \ref{plimline}.

In the proof of Claim 2 we use the following asymptotic formula: 
\begin{equation} c(b_{r,k})=\clr c(b_k)+o(c(b_k)), \text{ as } k\to\infty.
\label{asbk}
\end{equation}
\begin{proof} {\bf of (\ref{asbk}).} Let us calculate the cocycle value $c(b_{r,k})$. Set 
$$\chi_{j,k}=g_jb_k^{r_j}g_j^{-1}.$$ 
One has $b_{r,k}=\chi_{1,k}\dots\chi_{n,k}$.  Hence, by (\ref{cocn}), 
\begin{equation}c(b_{r,k})=c(\chi_{1,k})+Ad_{\chi_{1,k}}c(\chi_{2,k})+
\dots+Ad_{\chi_{1,k}\dots\chi_{n-1,k}}c(\chi_{n,k}).\label{prodchi}
\end{equation}

{\bf Claim 3.} $ c(\chi_{j,k})=r_jAd_{g_j}c(b_k)+o(c(b_k)),$ {\it as} $k\to\infty$.

\medskip
\begin{proof} For simplicity we prove the claim assuming that $r_j>0$; for  
$r_j\leq0$ the proof is analogous and uses (\ref{cocin}). 
Writing $\chi_{j,k}=g_jb_k\dots b_kg_j^{-1}$ and applying (\ref{cocn}), 
(\ref{cocin}) yields  
$$c(\chi_{j,k})=c(g_j)+(Ad_{g_j}+Ad_{g_jb_k}+\dots+Ad_{g_jb_k^{r_j-1}})c(b_k)+
Ad_{g_jb_k^{r_j}}c(g_j^{-1})$$
$$=c(g_j)+Ad_{g_j}(Id+Ad_{b_k}+\dots+Ad_{b_k}^{r_j-1})c(b_k)- 
Ad_{g_jb_k^{r_j}}Ad_{g_j^{-1}}c(g_j)$$
\begin{equation}
=(Id-Ad_{\chi_{j,k}})c(g_j)+r_jAd_{g_j}c(b_k)+
Ad_{g_j}\sum_{i=1}^{r_j-1}(Ad_{b_k}^i-Id)c(b_k).
\label{chij}\end{equation}
Recall that (since $g_j$'s are fixed), as $k\to\infty$,  
\begin{equation}
dist(b_k,1)=o(1)=o(c(b_k)),\  dist(\chi_{j,k},1)=O(dist(b_k,1))=o(1)=o(c(b_k)), 
\label{auxas}\end{equation}
by (\ref{limline}) and the definition of $\chi_{j,k}$. Hence, 
the adjoints 
$Ad_{b_k}$, $Ad_{\chi_{j,k}}$ are $O(dist(b_k,1))=o(1)=o(c(b_k))$- close to the identity. 
This implies that the first and 
third terms in the right-hand side of (\ref{chij}) are $o(c(b_k))$. This together with 
(\ref{chij}) proves the claim.
\end{proof}

Claim 3 together with (\ref{prodchi}) and the latter  asymptotics of adjoints implies 
(\ref{asbk}). 
\end{proof}

\begin{proof} {\bf of Claim 2.} One has 
$$\lim_{k\to\infty}b_{r,k}=1,$$ 
since $b_k\to1$ and $g_j$'s are fixed. Moreover, $c(b_{r,k})\neq0$ for large $k$ and 
$\rr c(b_{r,k})\to\Lambda_r$, by  (\ref{asbk}), the invertibility of the operator 
$\clr:\Lambda\to\Lambda_r$ 
and since $c(b_k)\neq0$ and $\rr c(b_k)\to\Lambda$, see
(\ref{limline}). Now for the proof of Claim 2 it suffices to show that 
\begin{equation} dist(b_{r,k},1)=o(c(b_{r,k})), \text{ as } k\to\infty.\label{obk}
\end{equation}
Indeed, one has 
\begin{equation} dist(b_{r,k},1)=O(dist(b_k,1))=o(c(b_k)),\label{dbk}
\end{equation}
by definition and (\ref{limline}). On the other hand, 
\begin{equation}c(b_k)=O(c(b_{r,k})).\label{ocbk}\end{equation}
This follows from the invertibility of the operator $\clr:\Lambda\to\Lambda_r$, 
(\ref{limline}) and (\ref{asbk}). 
Statements (\ref{dbk}) and (\ref{ocbk}) imply (\ref{obk}). This proves Claim 2. 
\end{proof}

Claim 2 implies Proposition \ref{plimline}, as was shown above. 
\end{proof}
 
 \begin{proof} {\bf of Proposition \ref{climline}.} 
We say that an element $g\in U$ is {\it $\sigma$- regular}, if 
$d\sigma(g)\neq0$. Set $\mathcal S=\{h\in U\cap\G \ | \ h 
\text{ is } \sigma- \text{ regular}\}$. Each 
$g\in U$ is a limit of  elements of $\mathcal S$ 
(by the local nonconstance of $\sigma$ and the  density of  
$\Gamma$). Hence, it suffices to prove Proposition \ref{climline} for every 
$g\in\mathcal S$. Fix a $g\in\mathcal S$. If $\frac{d\sigma}{d\nu(g)}\neq0$, we set 
$h_k=g$, and we are done. Thus, in what follows, we assume that 
\begin{equation}\frac{d\sigma}{d\nu(g)}=0.\label{dng}\end{equation}
Fix a one-dimensional linear subspace $\Lambda'\subset T_gG$ transversal to the 
 hypersurface $\sigma=\sigma(g)$, which is smooth by the $\sigma$- regularity 
 of $g$. Set 
$$\Lambda=\Lambda'g^{-1}\subset\gg.$$
By Proposition \ref{plimline}, there exists a sequence $b_k\in\G$ satisfying 
(\ref{limline}). Set 
$$ h_k=b_kg.$$
Then $\lim_{k\to\infty}h_k=g$. Now it suffices to prove the inequality  
$\frac{d\sigma}{d\nu(h_k)}\neq0$ for large $k$. To do 
this, we decompose $\nu(h_k)$ via the cocycle identity (\ref{coc}): 
$c(h_k)=c(b_k)+Ad_{b_k}c(g)$. Thus, 
\begin{equation}
\nu(h_k)=c(h_k)h_k=v_k+\zeta_k, \ v_k,\zeta_k\in T_{h_k}G, \ 
v_k=c(b_k)h_k, \ \zeta_k=(Ad_{b_k}c(g))h_k.\label{vwk}\end{equation}
We show next that the derivative $\frac{d\sigma}{dv_k}$ dominates 
$\frac{d\sigma}{d\zeta_k}$ (Claims 4 and 5 below). 

{\bf Claim 4.} {\it There exists a $p>0$ such that for every $k$ large enough}
\begin{equation}|\frac{d\sigma}{dv_k}|>p|c(b_k)|.
\label{ds>}\end{equation}

\begin{proof} The sequence of lines $\rr v_k$ tends to $\Lambda g=\Lambda'$ (in 
the projectivization of the tangent bundle over $G$), as $k\to\infty$, 
by (\ref{vwk}) and since $\rr c(b_k)\to\Lambda$, see (\ref{limline}), and $h_k\to g$. 
The limit line $\Lambda'$ is tranversal to the level 
hypersurface $\sigma=\sigma(g)$ at the $\sigma$- regular point $g$. 
Therefore, there exists a $p>0$ such that for every $k$ large enough 
$|\frac{d\sigma}{dv_k}|>p|v_k|$. This together with the equality 
$|v_k|=|c(b_k)|$ (the right invariance of the metric) implies (\ref{ds>}).
\end{proof}

Let $G_g$ be the component of $G$ containing $g$. 
We measure the distances between the vectors in $TG_g$ 
 (that may belong to distinct tangent spaces) in the  
 metric induced from $G$. 
 
 \medskip
 
 {\bf Claim 5.} $dist(\zeta_k,\nu(g))=O(dist(b_k,1))$, as $k\to\infty$. 
 
 \medskip
 
 \begin{proof} One has $\zeta_k=(Ad_{b_k}c(g))h_k$, $Ad_{b_k}=Id+
 O(dist(b_k,1))$, as $k\to\infty$. Hence, 
 \begin{equation}\zeta_k=
 c(g)h_k+O(dist(b_k,1));\label{zetak}\end{equation}
 \begin{equation}
 dist(\zeta_k,\nu(g))\leq dist(\zeta_k,c(g)h_k)+dist(c(g)h_k,\nu(g)).
 \label{zety}\end{equation}
 The first term in the latter right-hand side is $O(dist(b_k,1))$, by 
 (\ref{zetak}). The second term is also $O(dist(b_k,1))$, since 
 $\nu(g)=c(g)g$ and $dist(h_k,g)=dist(b_k,1)$ 
 (the right invariance of the metric). 
 This together with (\ref{zety}) proves the claim.
 \end{proof}
 
 Claim 5 together with (\ref{dng}) implies that 
 $\frac{d\sigma}{d\zeta_k}=O(dist(b_k,1))$, as $k\to\infty$. 
 Hence, $|\frac{d\sigma}{d v_k}|>|\frac{d\sigma}{d\zeta_k}|$ for every $k$ 
 large enough, by Claim 4 and (\ref{limline}). This together with 
 the equality $\nu(h_k)=v_k+\zeta_k$  implies that 
 $\frac{d\sigma}{d\nu(h_k)}\neq0$. Proposition \ref{climline} is proved.
  \end{proof}

\section{Groups  without proximal elements: proof of Lemma \ref{lem1}}

\def\psijhe{\Psi_{j,\wt h,\var}}
\def\psihe{\Psi_{\wt h,\var}}
\def\wtpsijhe{\wt\Psi_{j,\wt h,\var}}
\def\wtpsihe{\wt\Psi_{\wt h,\var}}

In this section we assume that the group $G$ contains no proximal elements and prove 
Lemma \ref{lem1} for this case. Together with the discussion at the 
end of Subsection 4.1, this proves Theorem \ref{tconj}. 

Set $\rho=\rho_0$. In what follows we will be using elements $g_1,\dots,g_n,h\in\G$ 
such that  $A_j=\rho(g_j)\in\pc1$ and $\rho(h)$ is sufficiently close to $1\in G$, so that 
  $\wt h=\log(\rho(h))\in\gg$ is well-defined and $(A_j,\rho(h))\in\Pi'$ for all $j$. Here 
  $\log$ is a branch of the inverse to $\exp:\gg\to G$, so that $\log(1)=0$.  
  Then $\nu_j=v_{A_j}(\exp(\wt h))$, by (\ref{nui}). Set $A=(A_1,\dots,A_n)$.  We 
consider the mapping $\Psi=\psihe$ corresponding to fixed $g_j$, $l_j$ and 
variable $\wt h$, $\var$. We first compute the asymptotics of its derivative at 0, as 
$\wt h,\var\to0$ (Proposition \ref{propas}, which is an analogue of Proposition 
\ref{asfor}). Then we show that the dominant term 
in its asymptotic expression is an injective linear map for appropriate 
$g_j$, $l_j$ and $\wt h$ (Proposition \ref{lem2}). 

Recall that in the case under consideration $s_j(u)$ are complex-valued functions, 
$0<|s_j(u)|<1$, and $L(A_j)$ are planes equipped with complex structures. 
We write the mapping $\psihe$ and its derivative using the following notations. 
Let us introduce  linear 1- forms $\sigma_j$, linear operator $\Sigma$ 
and auxiliary functions $\Lambda_{j,\wt h,\var}(\wt u)$, $\psijhe(\wt u)$: 
\begin{equation}
\sigma_j=d\ln s_j(0): T_0\rr^n\to\cc, \ 
 \Sigma=(\sigma_1,\dots,\sigma_n):T_0\rr^n\to\cc^n,\label{ajrh}\end{equation}
 \begin{equation} \Lambda_{j,\wt h,\var}(\wt u)=
 \var s^{l_j}(A_j)e^{\sigma_j(\wt u)}
v_{A_j}(\exp(\wt h))\in L(A_j), \ 
 \psijhe(\wt u)=\exp\circ\Lambda_{j,\wt h,\var}(\wt u).\label{arr} \end{equation}
 Then
 \begin{equation}
 \Psi(\wt u)=\Psi_{\wt h,\var}(\wt u)=
 \Psi_{1,\wt h,\var}(\wt u)\dots\Psi_{n,\wt h,\var}(\wt u). \label{rappel}\end{equation}
 Set
 $$\wtpsihe(\wt u)=\psihe(\wt u)(\psihe(0))^{-1},$$
\begin{equation}
 Y_j=Y_j(A_j,\wt h,l_j)=s^{l_j}(A_j)dv_{A_j}(\wt h)\in L(A_j), 
 \ Y=Y(A,\wt h,l)=(Y_1,\dots,Y_n).
 \label{yjlh}\end{equation}
 
 {\it Sketch of the proof of Lemma \ref{lem1}.} 
 Recall that $dv_{A_j}$ denotes the differential at 1 of the vector function $v_{A_j}$, 
 see Remark \ref{remdv}. 
 Proposition \ref{propas} below provides an asymptotic formula for the derivative  
 $\wtpsihe'(0)$, as $\wt h,\var\to0$. The main asymptotic 
 term there is $\var$ times the linear operator 
\begin{equation}
\Omega_{Y,\Sigma}=\sum_{j=1}^nY_j\sigma_j:T_0\rr^n\to\gg,\label{omys}\end{equation}
which is defined as in Remark \ref{complex}. 
The principal part of the proof of Lemma \ref{lem1} is Proposition \ref{lem2}, 
which says that the operator $\oys$ is invertible for appropriately chosen 
$g_j$, $l$ and $\wt h$. To prove that, we first choose (using the density of $\rho(\G)$ 
in $G_0$  and Corollary \ref{sind}) some 
``preliminary'' $g_j\in\G$ so that  $A_j\in\pc1$, 
\begin{equation}
 \text{the operator } \Sigma:T_0\rr^n\to\cc^n \text{ is injective,  \ and } 
\a(A_j)=\frac{\arg s(A_j)}{\pi}\notin\cup_{q=1}^n\frac1q\zz.\label{condg}\end{equation}
Proposition  \ref{lem3} shows that after replacing these $g_j$ 
by appropriate conjugates $\wt g_j=h_jg_{r_j}h_j^{-1}$ 
(which consists of replacing $A_j$ by  a conjugate $\wt A_j$ of $A_{r_j}$ and  
$\sigma_j$ by $\sigma_{r_j}$) one can 
choose some $Y_j\in L(\wt A_j)$ such that the corresponding operator 
$\oys$ from (\ref{omys}) is invertible. Here we use only the 
injectivity  condition from  (\ref{condg}) and the irreducibility of the adjoint of the group 
$\overline{\rho(\G)}$.  
For the new $\wt g_j$, Proposition 
\ref{lem4}  implies that the latter $Y_j$ can be realized as $Y_j(\wt A_j,\wt h,l_j)$ 
for some $\wt h\in\gg$ and $l_j\in\zz$. Here we use the 
argument condition from (\ref{condg}), which is invariant under conjugation.  
Lemma \ref{lem1} is 
deduced from Propositions \ref{propas} and \ref{lem2} at the end of the subsection. 

We now proceed with the proofs. 

\begin{proposition} \label{propas} Let $G$ be a real Lie group. Let  
$\Pi_{\cc,1}$ be the same, as 
in (\ref{pc1}), $A_1,\dots,A_n\in\Pi_{\cc,1}$, $l\in\zz^n$, 
$\wt h\in\gg$. Let $\Sigma=(\sigma_1,\dots,\sigma_n)$ 
be a collection of $\rr$- linear complex-valued 1- forms $\sigma_j:T_0\rr^n
=\rr^n\to\cc$. 
Let $v_{A_j}(y)$ be the vector function from Proposition \ref{pcomm}. Let  
$\psijhe$, $\psihe$, $\wtpsihe$, $Y=Y(A,\wt h,l)$,  $\oys$ be as in  
(\ref{arr})-(\ref{omys}). Then 
\begin{equation}\wtpsihe'(0)=\var(\oys+O(|\wt h|^2)), \text{ as } 
\wt h,\var\to0.\label{psi'to}
\end{equation}
\end{proposition}

\begin{proof} Set $\wtpsijhe(\wt u)=\psijhe(\wt u)(\psijhe(0))^{-1}$. 
\medskip

{\bf Claim.} $\wtpsijhe'(0)=\var(Y_j+
O(|\wt h|^2))\sigma_j,$ {\it as} $\var,\wt h\to0$.
%\label{psihe}\end{equation}

\begin{proof} Let 
$\Lambda_{j,\wt h,\var}(\wt u)=\log\psijhe(\wt u)\in L(A_j)\subset\gg$, see 
(\ref{arr}), $T_{j,\wt h,\var}:G\to G$ be the right translation    
$T_{j,\wt h,\var}(g)=g(\psijhe(0))^{-1}$. One has 
\begin{equation}\wtpsijhe=T_{j,\wt h,\var}\circ\exp\circ\Lambda_{j,\wt h,\var}. 
\label{compos}\end{equation}
Let us calculate asymptotically the derivative at 0 of the latter composition. 
Recall that $v_{A_j}(1)=0$, see Proposition \ref{pcomm}, $\exp'(0)=Id$. Hence, 
\begin{equation}v_{A_j}(\exp(\wt h))=dv_{A_j}(\wt h)+O(|\wt h|^2)=O(|\wt h|), 
\text{ as } \wt h\to0.\label{vath}\end{equation}
Substituting (\ref{vath}) and $e^{\sigma_j(0)}=1$ to the expressions (\ref{arr}) for 
$\Lambda_{j,\wt h,\var}$ and $\psijhe$ yields  
\begin{equation}
\Lambda_{j,\wt h,\var}(0)=\var s^{l_j}(A_j)(dv_{A_j}(\wt h)+O(|\wt h|^2))
 =O(\var|\wt h|), \  
 dist(\psijhe(0),1)=O(\var|\wt h|).\label{disthe}\end{equation}
Recall that $Y_j=Y_j(A_j,\wt h,l_j)=s^{l_j}(A_j)dv_{A_j}(\wt h)=O(|\wt h|)$. 
Hence, by (\ref{arr}) and (\ref{disthe}), 
\begin{equation}
\Lambda_{j,\wt h,\var}'(0)=\Lambda_{j,\wt h, \var}(0)\sigma_j=
\var (Y_j+O(|\wt h|^2))\sigma_j.\label{lahe}
\end{equation}

Writing  $\wtpsijhe$ and $T_{j,\wt h,\var}$ in the exponential chart 
as mappings $\rr^n\to\gg$ and $\gg\to\gg$ yields 
$\wtpsijhe=T_{j,\wt h,\var}\circ\Lambda_{j,\wt h,\var}$, $T_{j,\wt h,\var}'=Id+
O(\var|\wt h|)$ 
 uniformly on compact sets in $\gg$, by (\ref{compos}) and (\ref{disthe}). 
 This together with (\ref{lahe})  and  the asymptotics $Y_j=O(|\wt h|)$ implies the claim.
 \end{proof}

We can now finish the proof of Proposition \ref{propas}. Similarly to (\ref{cocpr}), 
one has 
\begin{equation}
\wtpsihe'(0)=\wt\Psi_{1,\wt h,\var}'(0)+Ad_{\Psi_{1,\wt h,\var}(0)}
\wt\Psi_{2,\wt h,\var}'(0)+\dots+
Ad_{\Psi_{1,\wt h,\var}(0)\dots\Psi_{n-1,\wt h,\var}(0)} 
 \wt\Psi'_{n,\wt h,\var}(0).\label{cocneww}\end{equation}
Each adjoint in (\ref{cocneww}) is $Id+O(\var|\wt h|)$, by (\ref{disthe}), as in the 
proof of Proposition \ref{asfor}. Substituting 
this and the formula from the claim to (\ref{cocneww}) yields (\ref{psi'to}), as in the 
proof of Proposition \ref{asfor}.  
 \end{proof}

\begin{proposition} \label{lem2} Let $\G$, $G$, $n$, $\rho_u$ satisfy Assumption 
\ref{*} at the beginning of Section 4, and $G$ have no proximal elements. 
Let $\pc1\subset G$ be 
the subset from (\ref{pc1}). Then there exist elements 
$g_1,\dots,g_n\in\G$  with $A_j=\rho_0(g_j)\in\pc1$ that satisfy 
the following statements. Let $s_j(u)=s(\rho_u(g_j))$, $\Sigma$ be the same, 
as in (\ref{ajrh}), $v_{A_j}(y)$ be the vector function from Proposition \ref{pcomm}, 
$dv_{A_j}$ be its differential at 1. Then for every 
$h'\in\gg$ such that $dv_{A_j}(h')\neq0$ for all $j$ there exists an 
 $l\in\zz^n$ such that the  operator 
$\Omega_{Y(A,h',l),\Sigma}$, see (\ref{yjlh}), (\ref{omys}), is invertible.
\end{proposition}
\begin{proof} Step 1:  ``preliminary'' elements $g_j$ satisfying 
(\ref{condg}). Let $\Pi^0\subset\hat\Pi=\pc1$ be the open set from Corollary 
\ref{sind}: $\Pi^0\neq\emptyset$ and 
$\arg s(g)|_{\Pi^0}$ is locally nonconstant (see Remark \ref{hat=pi}  
and the proof of Corollary \ref{sind}). Fix arbitrary 
$A_1',\dots,A_n'\in \Pi^0$ 
such that $\a(A_j')=\frac{\arg s(A_j')}{\pi}\notin \mathbb Q$. Fix  
a $\delta>0$ so that the $\delta$- neighborhood of each $A_j'$ 
consists of elements $g\in\Pi^0$ with $\a(g)\notin\cup_{q=1}^n\frac1q\zz$. 
The elements $g_j\in\G$ given by Corollary \ref{sind}  
with $A_j=\rho_0(g_j)$ being $\delta$- close to $A_j'$ and $\Sigma=(\sigma_1,\dots,
\sigma_n)$ with $\sigma_j=d\ln s_j(0)$ satisfy (\ref{condg}). Indeed, the inequality on 
arguments from (\ref{condg}) follows by construction. The 
operator $\Sigma$ is injective, since its imaginary part 
$(d(\arg s_1)(0),\dots,d(\arg s_n)(0))$ is injective (Corollary \ref{sind}). 

Step 2: making the operator $\oys$ invertible for some $Y$. Fix elements 
$g_1,\dots,g_n$ constructed on Step 1. The Lie group $\overline{\rho_0(\G)}$ 
and the above $A_j$ and 
$\sigma_j$ satisfy the conditions of Proposition \ref{lem3}, by construction. 
By Proposition \ref{lem3}, there exist $H_1,\dots,H_n\in \overline{\rho_0(\G)}$,  
$r_1,\dots,r_n\in\{1,\dots,n\}$  and $Y_j\in Ad_{H_j}L(A_{r_j})$, $j=1,\dots,n$, 
such that the operator $\sum_jY_j\sigma_{r_j}:T_0\rr^n\to\gg$ is 
invertible. Without loss of generality we can assume that $H_j=\rho_0(h_j)$ 
for some $h_j\in\G$ (density and persistence 
of invertibility under  perturbations). Fix these 
$h_j$ and set 
$$\wt g_j=h_jg_{r_j}h_j^{-1}, \ 
\wt A_j=\rho_0(\wt g_j), \ \wt s_j(u)=s(\rho_u(\wt g_j)),\ \wt\sigma_j=
\sigma_{r_j}=d\ln\wt s_j(0),  \ 
 \wt\Sigma=(\wt\sigma_1,\dots,\wt\sigma_n).$$
Step 3: Let $\wt g_j$, $\wt A_j$, $\wt\Sigma$ be as above, 
$\wt A=(\wt A_1,\dots,\wt A_n)$. Fix an arbitrary $h'\in\gg$ such that 
$$\wt Y_j=dv_{\wt A_j}(h')\neq0 \text{ for each } j.$$
Set $\wt Y=(\wt Y_1,\dots,\wt Y_n)$. Recall that the operator $\Omega_{Y,\wt\Sigma}
=\sum_jY_j\wt\sigma_j$ 
is invertible for a certain $Y=(Y_1,\dots,Y_n)$, $Y_j\in L(\wt A_j)$, and 
$Y_j(\wt A_j,h',l_j)=s^{l_j}(\wt A_j)\wt Y_j$, see (\ref{yjlh}). Therefore, 
there exists a vector $l=(l_1,\dots,l_n)\in\zz^n$ such 
that the operator $\Omega_{Y(\wt A,h',l),\wt\Sigma}$ is invertible, 
by (\ref{condg}) and Proposition \ref{lem4} (applied to $E=\gg$, 
$\La_j=L(\wt A_j)$, $s_j=s(\wt A_j)=s(A_{r_j})$ and $\wt\Sigma$). 
This proves Proposition \ref{lem2} for the above $\wt g_j$ and $\wt\Sigma$.
\end{proof} 
 
\begin{proof} {\bf of Lemma \ref{lem1}.} Fix $g_1,\dots.g_n\in\G$, $h'\in\gg$ 
and $l\in\zz^n$ satisfying the statements of Proposition \ref{lem2}: the operator 
$\Omega_{Y(A,h',l),\Sigma}$, see (\ref{ajrh}), (\ref{yjlh}), (\ref{omys}), is invertible. 
Let us show that there exists an element $h\in\G$ 
such that for every $\var>0$ small enough the derivative $\Psi'(0)$ is 
invertible. To do this, we fix a $c>0$ and an open cone $K\subset\gg$, 
$h'\in K$, such that for every $h''\in K$ and each $v\in T_0\rr^n$ one has 
\begin{equation}||\Omega_{Y(A,h'',l),\Sigma}(v)||\geq c||h''||||v||.\label{lowest}
\end{equation}
These $c$ and $K$ exist by the linear dependence of the operator 
$\Omega_{Y(A,h'',l),\Sigma}$ on the parameter $h''$ and the persistence of its 
invertibility under small perturbations. By Proposition \ref{propas}, 
\begin{equation}\wtpsihe'(0)=\var(\Omega_{Y(A,\wt h,l),\Sigma}+o(\wt h)), 
\text{ as }\wt h,\var\to0.\label{psito}
\end{equation} 
Fix a $\var_0>0$ and a neighborhood $W$ of zero in $\gg$ such that for every 
$\var<\var_0$ and each $\wt h\in W$
\begin{equation} ||o(\wt h)||\leq \frac c2||\wt h||.\label{owth}\end{equation}
(These $\var_0$ and $W$ exist by definition.) Fix a $h\in\G$ such that 
$\wt h=\log(\rho_0(h))\in W\cap K$ 
(it exists by the density of $\rho_0(\G)$ in $G_0$). Then for every $\var<\var_0$ 
and each $v\in T_0\rr^n$ one has 
$$||\wtpsihe'(0)v||\geq\frac c2\var||\wt h||||v||,$$
by (\ref{lowest})-(\ref{owth}). Thus, the derivative 
$\wtpsihe'(0)$, and hence, $\Psi'(0)$ is invertible for every $\var<\var_0$. 
The proof of Lemma \ref{lem1} is complete.
\end{proof}

\def\oh{\omega_h}

\section{Proof of Theorem \ref{thth} for arbitrary Lie group} 

Let $G$ be a real Lie group, $\hat G_{ss}$ denote its   
quotient by the radical of its identity component. The group $\hat G_{ss}$ 
is semisimple, hence $G_{ss}=Ad(\hat G_{ss})$ is an open subgroup in 
the semisimple linear algebraic group $Aut(\hat\gg_{ss})=Aut(\gg_{ss})$ 
(see \cite{VO}, p.214, theorem 1 and problem 3). 

\begin{definition} \label{sss} The above  group $G_{ss}$ will be called the {\it 
canonical semisimple part of} $G$. Let $p:G\to G_{ss}$ denote the quotient 
projection. Given a representation $\rho:\G\to G$, the induced representation 
$\rho_{ss}=p\circ\rho:\G\to G_{ss}\subset Aut(\gg_{ss})$ is called the canonical 
semisimple part of $\rho$. 
\end{definition}

Without loss of generality we assume that the representation 
$\rho:\G\to G$ is dense, passing to the Lie subgroup $\overline{\rho(\G)}\subset G$. 
Recall that $\G$ is a free group of rank $M\geq2$.

{\bf Case 1): $G\simeq G_{ss} \subset Aut(\gg)$ and $G$ 
has no nontrivial connected normal Lie subgroups.} 
One has $R(\G,Aut(\gg))=(Aut(\gg))^M$. Therefore, 
$dim_{[\rho]}X(\G,Aut(\gg))=(M-1)dim G\geq dim G$ (Proposition \ref{densecenter} applied to 
the representation $\rho:\G\to Aut(\gg)$, see Example \ref{exdense}). 
Hence, $\rho$ is a limit of non-injective representations (Theorem \ref{tconj}). 

{\bf Case 2): $G\simeq G_{ss}$, general case.} We show that $\rho$  is almost a 
 product of  representations for which Theorem \ref{thth} is already proved 
 (Proposition \ref{pprod}). We then 
deduce Theorem \ref{thth} for the representation $\rho$ using  purely algebraic 
Propositions \ref{itercom} and \ref{pnontriv} stated below. 

\begin{proposition} \label{pprod} Let $G$ be a semisimple Lie group  
 with trivial centralizer of the identity component. Then there exists a 
 collection of semisimple Lie groups $H_1,\dots,H_s$ with trivial centralizers of their 
 identity components and without 
 nontrivial connected normal Lie subgroups and an injective homomorphism 
$$\pi:G\to H_1\times\dots\times H_s$$
that is a local diffeomorphism, whose image $\pi(G)$ projects surjectively onto each 
factor $H_j$. 
\end{proposition}
\begin{proof} Each $Ad_G$- invariant linear subspace is an ideal in the 
semisimple Lie algebra $\gg$.  The only ideals are direct sums of 
simple factors of $\gg$. Therefore, the adjoint representation of $G$ is a 
direct sum of irreducible representations $Ad_G:\gh_j\to\gh_j$, $\gh_j$ 
are semisimple Lie algebras, and  
$\gg=\gh_1\oplus\dots\oplus\gh_s$. Let $G_j\subset G$ be 
the connected normal Lie subgroup with the Lie algebra $\gh_j$. Then take  
$H_i=G\slash\prod_{j\neq i}G_j$. Recall that $G$ has trivial centralizer of 
the identity component, which is equivalent to the injectivity of $Ad_G$. 
The adjoint $Ad_{H_i}=Ad_G|_{\gh_i}$ 
of each group $H_i$ is injective (as is that of $G$) and irreducible.  
Therefore, each group $H_i$ is semisimple with trivial centralizer of the identity 
component and has no nontrivial connected normal Lie subgroups. 
There is an obvious homomorphism from $G$ to 
$\prod H_i$ satisfying the required conditions. Proposition \ref{pprod} is 
proved.
\end{proof}
 
\begin{proposition} \label{itercom} Let $G=H_1\times\dots\times H_s$ be a direct 
product of groups $H_1,\dots, H_s$, $\pi_j:G\to H_j$ be the corresponding 
projections. Let $g_1,\dots,g_s\in G$ be such that 
$\pi_j(g_j)=1$ for all $j$. Then
\begin{equation}[\dots[[g_1,g_2],g_3],\dots,g_s]=1.\label{itcom1}\end{equation}
\end{proposition}
\begin{proof} One has $\pi_1([g_1,g_2])=\pi_2([g_1,g_2])=1$, since 
$\pi_1(g_1)=\pi_2(g_2)=1$. Analogously $\pi_j([[g_1,g_2],g_3])=1$ for every $j=1,2,3$, 
etc. This proves (\ref{itcom1}).  
\end{proof}

\def\glm{(\gamma_1,\dots,\gamma_M)}
\def\ga{\gamma}
\def\al{\alpha}
\def\gjk{\ga_{jk}}
\def\rjk{\rho_{jk}}

Let $\pi:G\to H_1\times\dots\times H_s$, $\pi_j:G\to H_j$ 
be respectively the  homomorphism from Proposition \ref{pprod} and its 
compositions with the projections to $H_j$. Each group $H_j$ is isomorphic to an open 
Lie subgroup of the linear algebraic group $Aut(\gh_j)$ (triviality of the centralizer of the 
identity component and semisimplicity).  Each representation 
$\pi_j\circ\rho:\G\to H_j$ is dense, as is $\rho$. Thus, Theorem 
\ref{thth} is already proved for the representations $\pi_j\circ\rho:\G\to H_j\subset 
Aut(\gh_j)$: for every $j=1,\dots,s$ there exist sequences 
$\rjk\in R(\G,H_j)$ and  $\gjk\in\G\setminus1$ such that 
\begin{equation}\rjk\to\pi_j\circ\rho, \text{ as } k\to\infty, \text{ and } 
\rjk(\gjk)=1.\label{rjkgjk}\end{equation}
(The representation $\pi_j\circ\rho$ may be not injective; in this case we  
choose $\rjk=\pi_j\circ\rho$ and $\gjk$ independent on $k$.) For every $k$ the 
collection  $\rho_{1k},\dots,\rho_{sk}$ lifts to a 
representation 
$$\rho_k:\G\to G \text{ such that } \pi_j\circ\rho_k=\rjk \text{ and } 
\rho_k\to\rho, \text{ as } k\to\infty$$
(the local diffeomorphicity of the mapping $\pi$ and the freeness of 
$\G$). For every $k\in\nn$ set 
\begin{equation}w_k=[\dots[[\ga_{1k},\ga_{2k}],\ga_{3k}],\dots,\ga_{sk}]: \ \ 
\rho_k(w_k)=1 \text{ for every } k,\label{wk=1} 
\end{equation}
by (\ref{rjkgjk}) and Proposition \ref{itercom} with  
$g_j=\pi\circ\rho_k(\ga_{jk})$. The next proposition shows 
that replacing some $\gjk$ by their conjugates by appropriate generators of 
$\G$, one can achieve that $w_k\neq1$. 

\begin{proposition} \label{pnontriv} Let $\G$ be a noncyclic free group. For every 
$s\in\nn$ and every $\gamma_1,\dots,\gamma_s\in\G\setminus1$ 
there exist elements $a_1,\dots,a_s\in\G$ (that may be chosen to 
be either trivial or some generators of  $\G$) such that for $\hat\ga_j=a_j\ga_ja_j^{-1}$  
we have $\hat w_s=[\dots[[\hat\ga_1,\hat\ga_2],\hat\ga_3],\dots,\hat\ga_s]\neq1$. 
\end{proposition}
\begin{proof} Induction on $s$.

Induction base: $s=1$. Set $a_1=1$. The element $\hat\ga_{1}=\ga_{1}$ 
is nontrivial by assumption.

Induction step for $s=2$. If $[\ga_{1},\ga_{2}]\neq1$, set $a_2=1$. 
Otherwise, if $[\ga_1,\ga_2]=1$, fix a generator $a_2$ of $\G$ that does not commute 
with $\ga_{2}$. Then by construction, $\hat w_{2}=[\ga_{1},\hat\ga_{2}]\neq1$. 

Induction step for arbitrary $s$. Let we have already constructed $a_j$ 
for $j\leq s-1$ so that $\hat w_{s-1}\neq1$. Then $a_{s}$ is constructed 
by applying the above induction step to $\ga_{1}$ replaced by $\hat w_{s-1}$ 
and $\ga_{2}$ replaced by $\ga_{s}$. By construction, $\hat w_{s}\neq1$. 
Proposition \ref{pnontriv} is proved. 
\end{proof}

We can assume that $w_k\neq1$ by replacing 
$\gamma_{jk}$ by their appropriate conjugates $a_{jk}\ga_{jk}a_{jk}^{-1}$, 
$|a_{jk}|\leq1$, by Proposition \ref{pnontriv}. (The latter inequalities will be used in the 
next section.) Equalities $\rho_{jk}(\gamma_{jk})=1$ and hence,   $\rho_k(w_k)=1$ 
remain valid. Thus, the representations $\rho_k$ are non-injective, 
and $\rho_k\to\rho$. This proves Theorem \ref{thth}.

{\bf Case 3): $G\not\simeq G_{ss}$.} Let $H\subset G$ denote the kernel of the quotient 
projection $p:G\to G_{ss}$. Its identity component $H_0$ is the radical of 
$G_0$. Let $\G$ be a noncyclic free group. 
\begin{proposition} \label{psolv} Let $\G$, $G$, $H$, $H_0$, $G_{ss}$ be as above, 
and let $s$ be the derived length of $H_0$:
$H^{(0)}_0=H_0$, $H^{(1)}_0=[H_0,H_0]$, $H^{(2)}_0=[H^{(1)}_0,H^{(1)}_0], \dots, 
H^{(s)}_0=1$;  $H^{(s-1)}_0\neq1$. 
For each $\gamma\in\G\setminus1$ there exists a $\wt \gamma\in\G\setminus1$ such that for every nondiscrete representation 
$\rho\in R(\G,G)$  with $\rho(\gamma)\in H$ one has  
\begin{equation}\rho(\wt\gamma)=1; \ |\wt \gamma|\leq 4^{s+2}(|\gamma|+l(\rho)), 
\ l(\rho)=\min\{ |a| \ | \ a\in\G\setminus1, \ \rho(a)\in G_0\}+2.
\label{solv}\end{equation}
\end{proposition}
\begin{proof} First observe that for every $g\in G_0$ and $h\in H$ one has $[g,h]\in H_0$. 
Indeed, since $H$ is normal in 
$G$, we have $[G,h]\subset H$. Moreover, since the commutator map is continuous 
and $G_0$ is the identity component of $G$, it follows that for every $h\in H$ the set 
$[G_0,h]$ is contained in $H_0$.  Fix an element $a\in\G\setminus1$  such that 
$\rho(a)\in G_0$. One can achieve that $a$ does not commute with $\gamma$ 
and $|a|\leq l(\rho)$, choosing first $a$ with the minimal possible value $|a|$ and 
then applying conjugation by appropriate generator of $\G$. (This conjugation may 
increase $|a|$ by at most 2.) 
Fix a generator $b$ of $\G$ that does not commute with $[a,\gamma]$. Set 
inductively 
\begin{equation}\sigma_{0,0}=[a,\gamma], \ \sigma_{0,1}=[b,\sigma_{0,0}],\ 
\sigma_{i+1,0}=[\sigma_{i,0},\sigma_{i,1}], \ \sigma_{i+1,1}=[\sigma_{i,0},
\sigma_{i,1}^{-1}]; \ \wt\gamma=\sigma_{s,1}.\label{sigmas}\end{equation}
The element $\wt\gamma\in\G$ is nontrivial and $\rho(\wt\gamma)=1$ by construction:  
$\rho(\sigma_{i,j})\in H^{(i)}_0$ for all $i$, $j$ by definition. Inequality 
(\ref{solv}) follows from definition. Proposition \ref{psolv} is proved.
\end{proof}

Let $p:G\to G_{ss}$ be the quotient projection, 
$\rho\in R(\G,G)$ be an injective dense representation. 
The representation $\rho_{ss}=p\circ\rho:\G\to G_{ss}$ is also dense and 
injective. Indeed, if 
$\rho_{ss}(\gamma)=1$ for some $\gamma\neq1$, then $\rho$ is not injective: 
$\rho(\gamma)\in H$, hence $\rho(\wt\gamma)=1$ for the corresponding $\wt\gamma\neq1$ from 
(\ref{solv}), - a conradiction to the injectivity of $\rho$. Thus, Theorem 
\ref{thth} is already proved for the representation $\rho_{ss}$. Hence,  
$\rho_{ss}$ is a limit of non-injective representations $\rho_{ss,k}:\G\to G_{ss}$. 
The latter can be lifted to representations $\rho_k:\G\to G$ such that 
$p\circ\rho_k=\rho_{ss,k}$ and $\rho_k\to\rho$,  as $k\to\infty$
(by the freeness of $\G$ and the submersivity of the projection). The 
representations $\rho_k$ are non-injective, as are $\rho_{ss,k}$, see 
the above argument for $\rho$ and $\rho_{ss}$. Hence, $\rho$ is 
a limit of non-injective representations $\rho_k$. The proof of Theorem 
\ref{thth} is complete.

 \medskip
 
 \begin{proof} {\bf of the Addendum to Theorem \ref{thth}.} 
Let $M$ denote the rank of $\G$. We prove the addendum by induction on 
$M$. Without loss of generality we assume that the representation $\rho$ is dense, 
replacing $G$ by its Lie subgroup $\overline{\rho(\G)}$.  
Fix an arbitrary sequence 
of non-injective representations $\rho_r:\G\to G$ converging to $\rho$. 
These representations are dense, whenever $r$ is large enough (Proposition \ref{open}). 
Set $F_r=\rho_r(\G)\subset G$. 

Induction base: $M=2$. The group $F_r$ is not cyclic for every $r$ large enough, 
since abelian groups cannot approximate nonabelian ones. Suppose that $F_r$ is 
free. Then it has to be of rank 2, and thus, isomorphic to $\G$. In other words, 
one obtains a non-injective epimorphism $\rho_r:\G\to F_r=\G$. This contradicts the 
fact that finitely-generated free groups are Hopfian.  

Induction step: $M>2$. Suppose the statement of the addendum is proved for the 
free groups of ranks less than $M$. Let us prove it for the given $\G$. Suppose 
$F_r$ are free for all $r$ large enough. Fix an arbitrary  $r$, for which 
$F_r$ is free and dense. Let $k$ denote the rank of $F_r$, and let $\G'$ denote 
the abstract free group of rank $k$. One has $k<M$, since $\rho_r:\G\to F_r\simeq\G'$ is 
a non-injective epimorphism and the free groups are Hopfian, 
as in the above argument.  The canonical identification 
$I:\G'\stackrel{\simeq}{\to}F_r$ is a dense injective
representation of $\G'$ to $G$. The group $\G'$ satisfies the statement of the 
addendum, by the induction hypothesis. Thus, the representation $I$ is a limit 
of non-injective representations $\phi_i:\G'\to G$ whose images are not free. The 
compositions $\phi_i\circ\rho_r:\G\to G$ are representations with the same 
non-free images that can be chosen arbitrarily close to $\rho$, taking $r$ and $i$ 
large enough. This together with a diagonal sequence argument yields 
a sequence of representations $\G\to G$ converging to $\rho$ with non-free images. 
The induction step is over. The addendum is proved.
\end{proof}

 \def\krl{k_r^{-1}}
 \def\rkr{\rho_{k_r^{-1}\wt u}}

\section{Approximations by non-injective representations. 
Proof of Theorem \ref{cdiap0}}

In the proof of Theorem \ref{cdiap0} we use results of approximations of 
elements of a semisimple Lie group by elements of its finitely-generated subgroups.
These results are stated in Subsection 7.1 and proved in the same subsection and 
Section 8. In Subsection 7.2 we prove Theorem \ref{cdiap0} as a corollary of  more 
general Theorem \ref{diap} stated at the same place. Theorem \ref{diap} 
is proved in Subsections 7.3 
(case when $G\simeq G_{ss}$ and $G$ has no nontrivial connected normal 
Lie subgroups), 7.4  (general case $G\simeq G_{ss}$) and 7.5 (case 
$G\not\simeq G_{ss}$).

\subsection{Approximations by elements of dense subgroups.} 
\begin{definition} \label{dnet} Given a metric space $E$, a  subset 
$K\subset E$ and a $\delta>0$, we say that a  subset $S$ in $E$ 
is a $\delta$- net on $K$, if for the $\delta$- neighborhood $S_{\delta}$ (respectively, 
$K_{\delta}$) of $S$ (respectively, $K$) we have $K\subset S_{\delta}$ and 
$S\subset K_{\delta}$. 
\end{definition}

In what follows (unless we say otherwise), given a point $a$ in a metric space (which 
will be either $\rr^n$, or a Lie group equipped with a 
Riemannian metric)  and $r>0$, we denote  
$$D_r(a) \ \text{the ball centered at} \ a \ \text{of radius} \ r, \ D_r=D_r(0)\subset\rr^n 
\ \text{(respectively,} \ D_r= D_r(1)\subset G_0).$$ 
Usually we measure 
the distance on a connected component of a Lie group in a given left-invariant 
Riemannian metric. We use the following property of a  
left-invariant distance. 

\begin{proposition} \label{pnet} Let $\delta_1,\delta_2>0$, $G$ be a connected 
 Lie group equipped with a left-invariant Riemannian metric, $K\subset G$ be an arbitrary 
 subset. Let 
$\Omega,\Omega'\subset G$ be two subsets such that $\Omega$ contains 
a $\delta_1$- net on $K$ and $\Omega'$ contains a $\delta_2$- net on the 
$\delta_1$- ball $D_{\delta_1}\subset G$. Then the product 
$\Omega\Omega'\subset G$ contains a $\delta_2$- net on $K$. 
\end{proposition}
\begin{proof} Take an arbitrary $x\in K$ and some its $\delta_1$- approximant 
$\omega\in\Omega$. Then $x'=\omega^{-1}x\in D_{\delta_1}$ 
(the left invariance of the metric). Take a $\delta_2$- approximant 
$\omega'\in\Omega'$ of $x'$. Then $\omega\omega'$ is a $\delta_2$- 
approximant of $x$: 
$dist(\omega\omega',x)=dist(\omega',x')<\delta_2.$ 
This proves the proposition.
\end{proof}

\begin{definition} A decreasing function $\var:\rr_+\to\rr_+$ is said to have 
{\it superlinear decay}, if there exists a $X>0$ so that for every $c>1$ and  $x\geq X$ 
we have 
\begin{equation} \var(cx)<c^{-1}\var(x).\label{condvar}\end{equation}
\end{definition}

\begin{example} For every $\kappa>0$ the function $\var(x)=e^{-x^{\kappa}}$ 
has superlinear decay. 
\end{example}

Given a free group $\G$ of rank $M$ and a dense representation $\rho:\G\to G$, 
we will approximate 
the elements of $G$ by elements $\rho(w)$, $w\in\G$, so that the approximation rate is controlled by a certain superlinear decay function 
$\var(|w|)$, as described below. 

\def\nn{\mathbb N}
\def\omk{\Omega_{K,m}}
\def\omr{\Omega_{D_R,m}}
\def\oml{\Omega_{D_1,m}}
\def\lmr{l_m(\rho,D_R)}
\def\lml{l_m(D_1)}

For every $w\in\G$ set 
\begin{equation}R_w:\rgg=G^M\to G: \ R_w(\phi)=\phi(w).
\label{rw}\end{equation}
Let $\rho:\G\to G$ be a dense representation, $\var(x)$ be a superlinear decay function.  
\begin{definition} \label{grapp} Let $\G$, $M$, $G$, $\rho$, $\var(x)$ be as above, 
$K\subset G_0$ be a bounded set. We say that $G$ is 
{\it $\var(x)$- approximable
on $K$ by $\rho(\G)$}, if there exist a $c=c(\rho,K)>0$, a sequence of numbers 
$l_m=l_m(\rho,K)\in\nn$ (called {\it length majorants}), 
$l_m\to\infty$, as $m\to\infty$, and a sequence of subsets 
$\Omega_{K,m}=\Omega_{K,\rho,m}\subset\G$  such that for every $m\in\mathbb N$ 
\begin{equation}|w|\leq l_m \text{ for each } w\in\Omega_{K,m} \text{ and}
\label{w<l}\end{equation}
\begin{equation}\text{the subset } \rho(\Omega_{K,m}) \text{ lies in } G_0 
\text{ and contains an }  
\var(cl_m)- \text{ net on } K.\label{onet}\end{equation}
We say that $G$ is $\var(x)$- approximable on $K$ by $\rho(\G)$ 
{\it with bounded 
derivatives}, if $c$, $l_m$ and $\omk$ satisfying (\ref{w<l}) and (\ref{onet}) 
may be chosen so that the subset $\cup_m\rho(\omk)\subset G_0$ is 
bounded and there exists a neighborhood $V\subset\rgg=G^M$ of 
$\rho$ such that  the mappings (\ref{rw}) corresponding to $w\in\cup_m\omk$ 
have  uniformly bounded derivatives on $V$.
\end{definition}

The next proposition shows that the  $\var(x)$- approximability is 
metric-independent. 

\begin{proposition}\label{pro2} Let $\G$, $G$, $\rho$, $\var(x)$ be as above, 
$K\subset G_0$ be a bounded set. Let $g_1$, $g_2$ be two complete 
Riemannian metrics on $G$. Let the group $G$ 
equipped with the metric $g_1$ be $\var(x)$- approximable on $K$ by 
$\rho(\G)$ (with bounded derivatives), $\omk\subset\G$, $l_m=l_m(\rho,K)$, 
$c_1=c(\rho,K)$ be the corresponding sequences and constant from  
(\ref{w<l}) and (\ref{onet}). Let 
$$\delta=\max_m\var(c_1l_m), \ K_{\delta} \text{ be the closed } 
\delta- \text{ neighborhood of } K 
\text{ in the metric } g_1.$$
Then the group $G$ equipped with the metric $g_2$ is 
also $\var(x)$- approximable on $K$ by $\rho(\G)$ (with bounded 
derivatives), with respect to the same sequences $\omk$, $l_m$ and the new constant 
$$c_2=c_2(\rho,K)=r^{-1}c_1, \ r=\max\{\sup_{x,y\in K_{\delta}, x\neq y}
\frac{d_{g_2}(x,y)}{d_{g_1}(x,y)}, 1\},$$
where $d_{g_j}$ is the distance in the metric $g_j$, $j=1,2$. 
\end{proposition}
\begin{proof} Each set $\rho(\omk)$ contains an $\var(c_1l_m)$- net on $K$ in the 
metric $g_1$. By definition, the latter net is contained in $K_{\delta}$ and is 
an $r\var(c_1l_m)$- net on $K$ in the metric $g_2$. One has    
$$r\var(c_1l_m)\leq\var(r^{-1}c_1l_m)=\var(c_2l_m), \text{ whenever } 
m \text{ is large enough,}$$
by  (\ref{condvar}). This proves the $\var(x)$- approximability 
in the metric $g_2$. Now let the group $G$ equipped with the 
metric $g_1$ be $\var(x)$- approximable on $K$ with bounded derivatives: the 
set $\cup_m\rho(\omk)$ be bounded and the mappings $R_w$, see 
(\ref{rw}), corresponding to 
$w\in\cup_m\omk$ have uniformly bounded derivatives on a bounded    
neighborhood $V$ of $\rho$  (in the metric $g_1$). Then 
the set $\wt V=\cup_{w\in\cup_m\omk}R_w(V)$ is bounded and hence, 
$\sup_{x,y\in\wt V, x\neq y}\frac{d_{g_2}(x,y)}{d_{g_1}(x,y)}<+\infty$. The latter 
inequality together with the above uniform boundedness of the derivatives 
on $V$ in the metric $g_1$ implies their uniform boundedness on $V$ in the 
metric $g_2$. This proves the proposition.
\end{proof}

The following proposition shows that the $\var(x)$- approximability on the unit ball centered at 1 implies the $\var$- approximability on arbitrary bounded subset 
in the identity component.  

\begin{proposition} \label{pro1} Let $\G$, $G$, $\rho$, $\var(x)$ 
be as above, and 
let the Riemannian metric on $G$ be left-invariant. Let $G$ be $\var(x)$- approximable 
by $\rho(\G)$ (with bounded derivatives) on the unit ball $D_1\subset G_0$. Let  
$c(\rho,D_1)$, $l_m(D_1)=l_m(\rho,D_1)$, $\Omega_{D_1,m}$ be the corresponding 
constant and sequences from  (\ref{w<l}) 
and (\ref{onet}). Let $R>1$, $\Omega_R\subset\G$ be a finite subset such that 
the set $\rho(\Omega_R)$ is contained in $G_0$ and forms a 1- net on 
$D_R\subset G_0$, 
$$l(R)=\max_{w\in\Omega_R}|w|.$$
Then $G$ is $\var(x)$- approximable on $D_R$ by $\rho(\G)$ 
(with bounded derivatives), where  
\begin{equation}\Omega_{D_R,m}=\Omega_R\Omega_{D_1,m}, \ 
l_m(D_R)=l_m(\rho,D_R)=l(R)+l_m(D_1), \ 
c(\rho,D_R)=\frac{c(\rho,D_1)}{l(R)+1}.\label{lr}\end{equation}
\end{proposition}

\begin{proof} Let $\Omega_{D_R,m}$, $l_m(D_R)$ be the finite sets and 
numbers given by (\ref{lr}). For every $m\in\nn$ the set $\rho(\omr)$ contains a 
$\delta(m)$- net on $D_R$,
$$\delta(m)=\var(c_1l_m(D_1)), \ c_1=c(\rho,D_1), $$
by Proposition \ref{pnet} applied to $K=D_R$,  $\Omega=\rho(\Omega_R)$, 
$\delta_1=1$,   $\Omega'=\rho(\oml)$, $\delta_2=\delta(m)$. (The latter satisfy 
the conditions of the proposition by  $\var(x)$- approximability on $D_1$.) 
One has 
$$|w|\leq l_m(D_R) \text{ for every } w\in\omr,$$
$$\delta(m)\leq\var(c_1(\inf_m\frac{l_m(D_1)}{l_m(D_R)})l_m(D_R))\leq
\var(c(\rho,D_R)l_m(D_R)).$$
This follows by definition, 
 (\ref{lr}), the inequality $\frac{\lml}{l_m(D_R)}\geq\frac 1{l(R)+1}$ and the decreasing of 
the function $\var(x)$. 
If in addition, the subset $\cup_m\rho(\oml)\subset G_0$ is bounded and 
the mappings (\ref{rw}) with $w\in\cup_m\oml$ have uniformly 
bounded derivatives on a neighborhood $V$ of $\rho$, then the same holds 
with $\oml$ replaced by $\omr$, since the  
collection $\Omega_R$ is finite. This proves the $\var(x)$- approximability 
on $D_R$ (with bounded derivatives) and Proposition \ref{pro1}.
\end{proof}

\begin{corollary} Any Lie group $\var(x)$- approximable by some $\rho(\G)$ on the 
unit ball in its identity component (with bounded derivatives) is $\var(x)$- 
approximable by $\rho(\G)$ on each bounded subset of its identity component 
(with bounded derivatives).
\end{corollary}

\begin{definition} \label{grapp10} Let $\G$ be a finitely-generated free group, $G$ 
be a Lie group. Let $\rho:\G\to G$ be a dense representation. 
We say that $G$ is {\it $\var(x)$- approximable 
(with bounded derivatives)} by  $\rho(\G)$, if it is $\var(x)$- approximable on each
bounded subset of the identity component $G_0$. We say that $G$ is 
{\it $\var(x)$- approximable (with bounded derivatives)}, 
if it is $\var(x)$-  approximable (with bounded derivatives) for every dense 
representation $\rho:\G\to G$. 
\end{definition}

Let $M$ denote the rank of $\G$. The following  well-known question is open.  
It was stated in \cite{sk3}, p.624 (without bounds of derivatives) for the groups $SU(n)$. 

\medskip

{\bf Question 4.} {\it Is it true that  each semisimple  Lie group having at least 
one irrational $M$- tuple of elements is always  
$\var(x)$- approximable 
with $\var(x)=e^{-x}$? If yes, does the same hold   with bounded 
derivatives?}

\medskip

\begin{theorem} \label{gapprox} Let $G$ be an arbitrary semisimple Lie group, 
$\G$ be a finitely-generated free group, $\rho:\G\to G$ be a dense representation. 
Then the group 
$G$ is $\var(x)$- approximable by $\rho(\G)$ with bounded derivatives, where 
\begin{equation}\var(x)=e^{-x^{\kappa}},\ \kappa=\frac{\ln1.5}{\ln9},\label{enka}
\end{equation} 
\begin{equation}\text{ the majorants } l_m=l_m(\rho,D_1) 
\text{ form a geometric progression } l_{m+1}=9l_m.\label{lkgensk}\end{equation}
\end{theorem}

Theorem \ref{gapprox} follows from Lemma \ref{lwsk} and Theorem 
\ref{gapprox1} (both formulated below). 

We will show that for many Lie groups the above approximation rate can be 
slightly improved. To state the corresponding result, let us introduce the following 

\begin{definition} \label{defsk} Let $\gg$ be  a Lie algebra with a fixed a positive definite 
scalar product on it. We say that $\gg$ {\it has surjective commutator}, if  
for each $z\in\gg$ there exist $x,y\in\gg$ such that 
\begin{equation}[x,y]=z.\label{commsurj}\end{equation}
We say that $\gg$ satisfies {\it the Solovay-Kitaev inequality}, 
if there exists a  $c>0$ such that for every 
$z\in\gg$ there exist $x,y\in\gg$ satisfying (\ref{commsurj}) and such that 
\begin{equation} |x|=|y|\leq c\sqrt{|z|}.\label{insk}
\end{equation}
\end{definition}

\begin{theorem} \label{thbr} (G.Brown, \cite{br}). Each complex semisimple Lie 
algebra and each real semisimple split Lie algebra 
(see \cite{VO}, p.288) have surjective commutator. 
\end{theorem}

It is also known that in every connected compact semisimple Lie group 
and in every  connected complex semisimple Lie group each element is 
a commutator \cite{goto, pw, ree}. The Lie algebras $\su_n$ satisfy the Solovay-Kitaev 
inequality \cite{sk1, sk2, sk3}. As it is shown in Subsection 8.1, Brown's  proof of  
Theorem \ref{thbr} together with simple estimates imply the following 

\begin{theorem} \label{thsk} 
Each complex semisimple Lie algebra and each real semisimple split Lie algebra 
satisfy the Solovay-Kitaev inequality.
\end{theorem}

{\bf Question 5.} {\it Is it true that each real semisimple Lie algebra has 
surjective commutator? If yes, is it true that it satisfies the Solovay-Kitaev 
inequality?}

The following theorem is implicitly contained in \cite{sk1, sk2, sk3}. 

\begin{theorem} \label{tsk} (R.Solovay, A.Kitaev).  Let $G$ be a 
Lie group whose Lie algebra satisfies the Solovay-Kitaev inequality, $\G$ be a 
finitely-generated free group, and $\rho:\G\to G$ be a dense representation. 
Then  $G$  is $\var'(x)$- approximable by $\rho(\G)$, where 
\begin{equation}\var'(x)=e^{-x^{\kappa'}},\ \kappa'=\frac{\ln1.5}{\ln5}, 
\label{var'}\end{equation} 
\begin{equation} \text{ the majorants } l_m=l_m(\rho,D_1) \text{ form a geometric 
progression } 
l_{m+1}=5l_m. 
\label{eta}\end{equation}
\end{theorem}

\begin{remark} \label{remsk} In fact, in  Theorem \ref{tsk} the Lie group is 
$\var'(x)$- 
approximable with bounded derivatives, with majorants 
$l_m(\rho,D_1)$ satisfying (\ref{eta}). 
This can be easily derived from 
Kitaev's proof \cite{sk1, sk2, sk3}. A proof of this statement  is outlined at the end of 
Subsection 8.2.  
\end{remark}

\begin{definition} \label{weaksk} 
Let $\gg$ be a Lie algebra with a fixed positive definite 
scalar product. We say that $\gg$ satisfies {\it the weak Solovay-Kitaev 
inequality}, if there exists a constant $c>0$ such that for every 
$z\in\gg$, there exist $x_j,y_j\in\gg$, $j=1,2$, such that 
\begin{equation}z=[x_1,y_1]+[x_2,y_2], \ |x_j|=|y_j|\leq c\sqrt{|z|}.
\label{wsk}\end{equation}
\end{definition}

\begin{remark} \label{rksk} 
The condition that a Lie algebra satisfies a (weak or strong) Solovay-Kitaev 
inequality is independent on the choice of the scalar product, while the corresponding 
constant $c$ depends on the scalar product: any two positive scalar products define 
equivalent norms. A Lie algebra satisfying 
the strong Solovay-Kitaev inequality obviously satisfies the weak one. 
\end{remark}

\begin{lemma} \label{lwsk} Each semisimple Lie algebra satisfies the weak 
Solovay-Kitaev inequality.
\end{lemma}

\begin{theorem} \label{gapprox1} Let a Lie group $G$ have a Lie algebra 
satisfying the weak Solovay-Kitaev inequality. Let $\G$ be a finitely-generated 
free group, $\rho:\G\to G$ be a dense representation. 
Then the group $G$ is $\var(x)$- approximable with 
bounded derivatives, where $\var(x)$, $l_m=l_m(\rho,D_1)$ 
are the same as in (\ref{enka}) and (\ref{lkgensk}) respectively.  
\end{theorem}

As it is shown in Subsection 8.1, Lemma \ref{lwsk} easily follows 
from Theorem \ref{thsk}. Theorem \ref{gapprox1} is proved in 
Subsection 8.2 analogously to the proof from \cite{sk1, sk2, sk3} of Theorem \ref{tsk}. 
Together, they imply Theorem \ref{gapprox}. 

\subsection{Approximations by non-injective representations} 
Let $G$ be a Lie group, $dim G>0$. We fix a Riemannian metric on $G$. We define 
the distance between any two points of $G$ to be the Riemannian distance, if they 
lie in the same connected component of $G$, and the infinity otherwise. Let $\G$ be 
a finitely-generated noncyclic free group. The above distance induces a distance 
 on the space of representations $\G\to G$ as a subset in $G^M$, $M=rank\G$. Let 
 $\rho:\G\to G$ be a nondiscrete injective representation. Theorem \ref{cdiap0} says 
that $\rho$ is $e^{-x^{\kappa}}$- approximable by non-injective representations, 
as in the following definition. 

\begin{definition} \label{aprel} 
Let $G$, $\G$, $\rho$ be as above. Let $\var(x)$ be a superlinear decay function. We say that 
$\rho$ is $\var(x)$- {\it approximable by non-injective representations}, if there exist 
a  $c=c(\rho)>0$ and  sequences of numbers $l_r\in\nn$ 
(called {\it the relation length majorants}), $l_r\to\infty$, as $r\to\infty$, elements 
 $w_r\in\G\setminus1$  and representations $\rho_r:\G\to G$ such that  
\begin{equation}\rho_r(w_r)=1, \ |w_r|\leq l_r \ \text{and} \ dist(\rho_r,\rho)<\var(cl_r) \ 
\text{for every} \ r\in\nn.\label{eapprox}\end{equation}
\end{definition}

\begin{remark} \label{aprmetr} 
The above definition and the corresponding  sequence $w_r$ are independent on the choice of the 
Riemannian metric on $G$ (while the constant $c$ depends on the metric). 
The proof of this statement is analogous to the proof of 
Proposition \ref{pro2}. 
\end{remark}
Let $G_{ss}$ and $\rho_{ss}:\G\to G_{ss}$ be the semisimple parts of $G$ and $\rho$ 
respectively,  
see Definition \ref{sss}. Note that $dim G_{ss}>0$. Indeed, otherwise  the group 
$G_0$, and hence, $\rho(\G)\cap G_0$ are solvable. But the latter is free noncyclic, 
being a nontrivial normal subgroup in the noncyclic free group $\rho(\G)$, - a 
 contradiction. 
 
 We will prove Theorem \ref{cdiap0} as a corollary of the following theorem. 
\begin{theorem} \label{diap} 
Let $\G$, $G$, $G_{ss}$, $\rho$, $\rho_{ss}$ be as above, and the representation  
$\rho_{ss}$ be dense. Let $\var(x)$ be a superlinear decay function, and the 
group $G_{ss}$ be $\var(x)$- approximable with bounded derivatives by 
$\rho_{ss}(\G)$ (see Definition \ref{grapp10}).  
Then $\rho$ is $\var(x)$- approximable by non-injective representations. 
\end{theorem}

{\bf Addendum to Theorem \ref{diap}.} {\it Under the assumptions of Theorem 
\ref{diap} let $l_m=l_m(\rho,D_1)$ be the  length majorants from 
(\ref{w<l}) for the $\var(x)$- approximations of $G_{ss}$ by $\rho_{ss}(\G)$ 
on the unit ball $D_1\subset G_{ss}$. 
Then there exist constants $d,q\in\nn$ depending only on $\rho$ such that 
$\rho$ is $\var(x)$- approximable by non-injective representations 
with relation length majorants} 
\begin{equation}l_m'=dl_m, \ m\geq q.\label{c''l}\end{equation}

\begin{proof} {\bf of Theorem \ref{cdiap0} modulo Lemma \ref{lwsk} and Theorems 
\ref{gapprox1}, \ref{diap}.} Without loss of generality we assume that 
the representation $\rho$ is dense, replacing $G$ by its Lie subgroup 
$\overline{\rho(\G)}$. Then  $\rho_{ss}$ is also dense. The function 
$\var(x)=e^{-x^{\kappa}}$ with $\kappa=\frac{\ln1.5}{\ln9}$ satisfies the conditions 
of Theorem \ref{diap} and its addendum with a majorant sequence 
$l_m$ such that $l_{m+1}=9l_m$ (Theorem \ref{gapprox} applied to the 
representation $\rho_{ss}$). This, together with Theorem \ref{diap} and its addendum, 
see (\ref{c''l}), implies the $\var(x)$- approximability of $\rho$ 
by non-injective representations with majorants $l_m'=dl_m$. 
This proves Theorem \ref{cdiap0}. 
\end{proof} 

\begin{corollary} \label{cdiap} Let $\G$, $G$, $G_{ss}$, $\rho$, $\rho_{ss}$ be as 
above. Let the semisimple part $\gg_{ss}$ satisfy the Solovay-Kitaev inequality, and 
the representation $\rho_{ss}$ be dense. Then $\rho$ is 
$\var'(x)=e^{-x^{\kappa'}}$- approximable  by non-injective representations, where 
$\kappa'=\frac{\ln1.5}{\ln5}$, see  (\ref{var'}). The corresponding relation length 
majorant  sequence $l_r$ can be chosen so that $l_{r+1}=5l_r$. 
\end{corollary}

Corollary \ref{cdiap} follows from Theorem \ref{diap} 
(with the addendum), Theorem \ref{tsk} and 
Remark \ref{remsk}, analogously to the above proof of Theorem \ref{cdiap0}.  Theorem 
\ref{diap} together with its addendum are proved in the next three subsections. 

\medskip

{\bf Question 6.} {\it Is it true that any 
injective dense representation of a finitely-generated noncyclic 
free group to a 
Lie group is always $e^{-x}$- approximable by non-injective representations?}

\medskip

By Theorem \ref{diap}, a positive  solution of Question 4 with bounded 
derivatives (see Subsection 7.1) would imply a positive answer to Question 6.

\subsection{Proof of Theorem \ref{diap}: case when $G\simeq G_{ss}$ and $G$ 
has no nontrivial connected normal Lie subgroups} 
The representation 
$\rho:\G\to G\subset Aut(\gg)$ satisfies (\ref{zv}) (with $G$ 
replaced by $Aut(\gg)$), analogously to the discussion at the beginning of Section 6. 
Fix a left-invariant Riemannian metric on $G$ and 
an $(Aut(\gg),X)$- immersive deformation $\rho_u$ of $\rho$,  see Definition \ref{dimmers}:  
\begin{equation}\rho_u:\G\to G, \ u\in\rr^n, \ \rho_0=\rho, \ 
||\frac{d\rho_u}{du}||\leq1. 
\label{rhou}\end{equation}
The inequality in (\ref{rhou}) may be achieved by a reparametrization of the 
given family $\rho_u$. 
 We show that $\rho$ is $\var(x)$- approximable by appropriate 
non-injective representations 
$\rho_{v_m}$. To do this, we fix a converging tuple 
$(\{ w_r\}, \{ k_r\}, \psi,\delta)$, see Definition \ref{dtuple}: 
$k_r\to+\infty$, as  $r\to\infty$,
\begin{equation} \rho_{k_r^{-1}\wt u}(w_r)\to\psi(\wt u),   \ 
\psi:\odd\to\psi(\odd)\subset G_0 \text{ is a diffeomorphism,} 
\label{krpsi}\end{equation}
the latter convergence being uniform with derivatives on $\odd$. Fix an $R>0$ such that 
\begin{equation}\psi(\odd)\cup(\psi(\odd))^{-1}\Subset D_R\subset G_0.\label{inclr}\end{equation}
Recall that $G$ is $\var(x)$- approximable on $D_R$ by $\rho_0(\G)$ with 
bounded derivatives. Let 
\begin{equation}\Omega_m=\omr\subset\G, \ \wt l_m=l_m(\rho_0,D_R), \ 
c=c(\rho_0,D_R)>0 
\label{omlc}
\end{equation}
be the corresponding subset and majorant sequences and constant
from Definition \ref{grapp}. One has 
\begin{equation}\rho_{k_r^{-1}\wt u}(ww_r)=\rho_{k_r^{-1}\wt u}(w)\rkr(w_r)\to
\chi_{w}(\wt u)=\rho_0(w)\psi(\wt u), \text{ as } r\to\infty,\label{rhokr}\end{equation}
uniformly in $\wt u\in\odd$, $w\in\cup_m\Omega_m$, with derivatives in 
$\wt u$. Indeed, $\rho_{k_r^{-1}\wt u}(w)\to\rho_0(w)$ in the latter sense: 
the derivatives $\frac{d\rho_u(w)}{du}$ are uniformly bounded on a fixed  
neighborhood $V=V(0)\subset\rr^n$ (Definition \ref{grapp} and (\ref{rhou}));  
$k_r\to\infty$ (in particular, $k_r^{-1}\odd\subset V$, whenever $r$ is large enough). 
This, together with (\ref{krpsi}), implies (\ref{rhokr}). 

Let $C>0$ be the Lipschitz constant of the diffeomorphism 
$\psi^{-1}|_{\psi(\odd)}$. The limit mappings $\chi_{w}(\wt u)|_{\odd}$, see 
(\ref{rhokr}), are also diffeomorphisms, and their inverses have the same 
Lipschitz constant on the images of $\odd$ (the left invariance of the metric). Fix a 
$r>0$ such that for every $w\in\cup_m\Omega_m$ the 
mapping $\wt u\mapsto\rho_{k_r^{-1}\wt u}(ww_r)$ is a diffeomorphism of 
$\odd$ onto its image,  its inverse is $2C$- Lipschitz there, 
and $\rho_0(w_r^{\pm1})\in D_R$. This $r$ exists by (\ref{rhokr}) and the convergence 
$\rho_0(w_r^{\pm1})\to(\psi(0))^{\pm1}\in D_R$, see (\ref{krpsi}) 
and (\ref{inclr}). For each $m\in\nn$ fix an element 
\begin{equation}
\omega_m\in\Omega_m \text{ such that } dist(\rho_0(\omega_m),\rho_0(w_r^{-1}))<
\var(c\wt l_m).
\label{omm}\end{equation}
It exists by Definition \ref{grapp} and the inclusion 
$\rho_0(w_r^{\pm1})\in D_R$. Set
$$\wt w_m=\omega_mw_r, \ \psi_m(\wt u)=\rho_{k_r^{-1}\wt u}(\wt w_m).$$ 
Then $\psi_m:\odd\to\psi_m(\odd)$ are diffeomorphisms,  
\begin{equation}\psi_m^{-1}:\psi_m(\odd)\to\odd \text{ are } 2C- \text{ Lipschitz 
diffeomorphisms,}\label{2lip}\end{equation} 
by the choice of $r$.  We assume that $\var(c\wt l_m)<R$, restricting ourselves 
to sufficiently large $m$'s. Thus, $\rho_0(\omega_m)\in D_{2R}$, by (\ref{omm}). 
Consider the right translations $G\to G$: $x\mapsto xg$. 
Let $K$ denote the maximum of norms of their derivatives in $x$ at 
$x\in\overline D_{2R}$ for all $g\in\overline D_{2R}$. 
One has $\psi_m(0)=\rho_0(\wt w_m)=\rho_0(\omega_m)\rho_0(w_r)$, 
\begin{equation}dist(\psi_m(0),1)\leq K
dist(\rho_0(\omega_m),\rho_0(w_r^{-1}))< K\var(c\wt l_m),\label{psimo}
\end{equation}
by the definition of $K$, the inclusions 
$\rho_0(w_r^{\pm1}),\rho_0(\omega_m)\in D_{2R}$ and 
(\ref{omm}). By (\ref{2lip}) and (\ref{psimo}), for every $m$ large enough 
(so that $K\var(c\wt l_m)<\frac{\delta}{2C}$) there exists a $\wt u_m\in D_{\delta}$ such that 
\begin{equation}\psi_m(\wt u_m)=1, \ |\wt u_m|\leq2Cdist(\psi_m(0),1)<
2CK\var(c\wt l_m)<\delta. \label{vmum}\end{equation}
Let $v_m=k_r^{-1}\wt u_m$, $l_m=l_m(D_1)$ be the length majorants for the $\var(x)$- approximations on $D_1$ by $\rho(\G)$:  
\begin{equation}\wt l_m=l(R)+l_m\geq l_m, \ l(R) \text{ is the same, as in 
(\ref{lr}).}\label{lrm}\end{equation}
By (\ref{vmum}), (\ref{lrm}) and elementary inequalities, we have:
$$\rho_{v_m}(\wt w_m)=\psi_m(\wt u_m)=1, \ |\wt w_m|\leq|\omega_m|+|w_r|\leq
\wt l_m+|w_r|\leq l_m'=dl_m, \ d=l(R)+|w_r|+1,$$
 $$dist(\rho_{v_m},\rho_0)\leq|v_m|\leq|\wt u_m|<2CK\var(c\wt l_m)\leq
 2CK\var(cl_m)\leq\var(c'l_m'), \ c'=\frac c{2CKd}$$
 for every $m$ large enough, 
 by the inequality in (\ref{rhou}), (\ref{vmum}) and  (\ref{condvar}). 
 The elements $\wt w_m\in\G$ are nontrivial, since 
 $\psi_m(\wt u)=\rho_{k_r^{-1}\wt u}(\wt w_m)\not\equiv1$, by the diffeomorphicity 
 of $\psi_m$. The three latter statements together imply that the injective 
 representation $\rho=\rho_0$ is $\var(x)$- approximable by non-injective 
 representations $\rho_{v_m}$ with relation length majorants 
 $l_m'=dl_m$. This proves Theorem \ref{diap} and its addendum.

\subsection{Proof of Theorem \ref{diap}: general case $G\simeq G_{ss}$} 
Let $\pi:G\to H_1\times\dots\times H_s$  
be the homomorphism from Proposition \ref{pprod}, $\pi_j:G\to H_j$ be the 
compositions of $\pi$ with the product projections. Theorem \ref{diap} and its 
addendum are already proved for each group $H_j$. 
We fix some left-invariant Riemannian metrics on the groups $H_j$. The $\pi$- pullback 
of their product yields a left-invariant metric on $G$. Each $\pi_j$ is surjective and maps 
 $D_1\subset G_0$ 
onto $D_1\subset H_j$ without increasing distances, by the choice of the metrics. 
Set $\rho_j=\pi_j\circ\rho$. 
Then, each representation $\rho_j$ is dense, as is $\rho$,  the 
group $H_j$ is $\var(x)$- approximable by $\rho_j(\G)$ on $D_1\subset H_j$ with 
bounded derivatives, as is $G$. The $\var(x)$- approximating 
subsets $\Omega_m=\Omega_{D_1,m}$, the majorants $l_m=l_m(D_1)$  and the 
constant $c=c(\rho,D_1)>0$ corresponding to $G$ and $\rho$ are the same for 
$H_j$ and $\rho_j$. Each  $\rho_j$ is a limit of a representation sequence 
$\rho_{jm}:\G\to H_j$ such that there exist a $d_j\in\mathbb N$ and a sequence  
$w_{jm}\in\G\setminus1$ for which
\begin{equation}\rho_{jm}(w_{jm})=1, \ |w_{jm}|\leq l_{jm}=d_jl_m, \ 
dist(\rho_{jm},\rho_j)<\var(cl_{jm})\leq\var(cl_m).  
\label{rhojm}\end{equation}
This follows from Theorem \ref{diap} and its addendum, applied to each 
representation $\rho_j:\G\to H_j$, and (\ref{condvar}). (If $\rho_j$ is not injective, then 
we set $\rho_{jm}=\rho_j$ and choose $w_{jm}$ independent on $m$, and set 
$d_j=|w_{jm}|$.) The representations $\rho_{jm}$ may be lifted to representations  
\begin{equation}
\rho_m:\G\to G, \ \pi_j\circ\rho_m=\rho_{jm}, \ 
dist(\rho_m,\rho)\leq s\max_jdist(\rho_{jm},\rho_j)<s\var(cl_m), 
\label{rhomj}
\end{equation}
 by the freeness of $\G$ and (\ref{rhojm}). By construction, $\rho_m\to\rho$, 
 as $m\to\infty$.  Set
 \begin{equation}w_m=[\dots[[w_{1m},w_{2m}],w_{3m}],\dots,w_{sm}]: \ 
\rho_m(w_m)=1,\label{wmcom}\end{equation}
 by (\ref{rhojm}) and Proposition \ref{itercom} applied to 
 $g_j=\pi\circ\rho_m(w_{jm})$. One can achieve that $w_m\neq1$  
replacing $w_{jm}$ by their conjugates with appropriate generators of 
 $\G$ (Proposition \ref{pnontriv}). The latter operation may increase $|w_{jm}|$ 
 no more than three times. 
 The representations $\rho_m$ and the elements $w_m\in\G\setminus1$ 
 thus constructed satisfy the following statements:
 $$\rho_m(w_m)=1, \ |w_m|\leq3\times 4^{s-1}\max_j|w_{jm}|\leq \wt l_m, 
 \ \wt l_m=d l_m, \ d=3\times 4^{s-1}\times\max_jd_j,$$ 
$$dist(\rho_m,\rho)<s\var(cl_m)\leq\var(\wt c\wt l_m), \ \wt c=\frac c{ds},$$
for large $m$, 
by (\ref{rhomj}). Thus, $\rho$ is $\var(x)$- approximable by non-injective representations 
$\rho_m$ with relation length majorants $\wt l_m=dl_m$. 
This proves Theorem \ref{diap} and its addendum.

\subsection{Proof of Theorem \ref{diap}: case $G\not\simeq G_{ss}$}
Let $H\subset G$ be the kernel of the quotient projection $p:G\to G_{ss}$.  Fix arbitrary 
left-invariant Riemannian metrics on $G$ and $G_{ss}$. The representation  
$\rho_{ss}=p\circ\rho$ is injective, as at the end of Section 6. 
Recall that by the conditions of Theorem \ref{diap}, 
the representation $\rho_{ss}$ is dense, and the group $G_{ss}$ is  
$\var(x)$- approximable by $\rho_{ss}(\G)$ on $D_1$ with bounded derivatives. Let 
$l_m$ be the 
corresponding length majorants. By Theorem \ref{diap} and its addendum, already 
proved for $\rho_{ss}$, there  exist   $c>0$, $d\in\mathbb N$ and  
representations $\rho_{ss,m}:\G\to G_{ss}$ and 
$w_m\in\G\setminus1$ such that for all sufficiently large $m$ 
\begin{equation}\rho_{ss,m}(w_m)=1, \ |w_m|\leq l_m'=dl_m, \ 
dist(\rho_{ss,m},\rho_{ss})<\var(cl_m').\label{rhoss}\end{equation}
The representations $\rho_{ss,m}$ may be lifted to representations $\rho_m:\G\to G$ 
such that 
\begin{equation}\rho_{ss,m}=p\circ\rho_m,
\ dist(\rho_m,\rho)\leq Kdist(\rho_{ss,m},\rho_{ss})<K\var(cl_m'),\label{distrm}
\end{equation}
where $K>0$ is a constant depending on the metrics on $G$ and $G_{ss}$ and 
$\rho$. One has $\rho_m(w_m)\in H$, by (\ref{rhoss}). 
Let $\wt w_m\in\G\setminus1$ be the elements  from (\ref{solv}) corresponding to $w_m$, $l(\rho_m)$ be the corresponding constant from (\ref{solv}). Note that 
$l(\rho_m)=l(\rho)$ for all $m$ large enough, since $\rho_m\to\rho$. Let $s$ 
be the derived length of $H_0$. Then  
$$\rho_m(\wt w_m)=1, \ 
|\wt w_m|\leq 4^{s+2}(|w_m|+l(\rho))\leq\wt l_m=\wt dl_m, \wt d=4^{s+2}(d+
l(\rho)),$$
$$dist(\rho_m,\rho)<K\var(cl_m')\leq K\var(cl_m)<\var(\wt c\wt l_m), \ \wt c=\frac c{K\wt d},$$
whenever $m$ is large enough, 
by (\ref{rhoss}), (\ref{distrm}) and (\ref{condvar}). Thus, $\rho$ is $\var(x)$- 
approximable by non-injective representations $\rho_m$ 
with relation length majorants 
$\wt l_m=\wt d l_m$. The proof of Theorem \ref{diap} and its addendum is complete.

\section{Approximability of semisimple Lie groups} 
In this section we first prove Theorem \ref{thsk} and Lemma \ref{lwsk}, 
and then Theorem \ref{gapprox1}.

\subsection{The Solovay-Kitaev inequality:  Theorem \ref{thsk} and Lemma \ref{lwsk}}

\begin{proof} {\bf of Theorem \ref{thsk}.} The proof repeats the arguments 
from \cite{br} (pp.765-767) with additional estimates. (The result and the arguments 
from the paper \cite{br} were independently rediscovered and explained to the author by 
J.-F.Quint.)  We prove Theorem \ref{thsk} for 
the real semisimple split Lie algebras $\gg$; the proof for the complex semisimple Lie algebras 
is analogous. 
Let $G$ be  the identity component of $Aut(\gg)$.
Fix a subalgebra $\gh\subset\gg$ corresponding to a maximal connected 
$\rr$- split torus in $Aut(\gg)$. 
Its adjoint action $\ad_{\gh}:\gg\to\gg$ is $\rr$- diagonalizable. The non-identically zero 
eigenvalues of the latter action are the values of the roots of $\gh$: the 
restrictions to $\gh$ of the 
roots of the complex Cartan subalgebra $\gh_{\cc}\subset\gg_{\cc}$. Set 
$$E=\text{ the sum of the root eigenlines of } \ad_{\gh}: \ \gg=\gh\oplus E; \ \ad_{\gh}(E)=E.$$
\begin{remark} \label{remsk1} A regular element $\wt x\in\gh$ induces a linear automorphism 
$\ad_{\wt x}:E\to E$. Therefore, each $z\in E$ is a Lie bracket: 
$[\wt x,(\ad_{\wt x}|_E)^{-1}z]=z$. 
\end{remark}

It was shown in \cite{br} that each element of $\gg$ has a conjugate in $E$. As it is 
shown below, Theorem \ref{thsk} is easily implied by the following more precise lemma.

\begin{lemma} \label{lemdec} There exists a compact subset $K\subset G$ such that 
each element of $\gg$ is conjugate to an element of $E$ by an element of $K$. 
\end{lemma}

\begin{proof} Without loss of generality 
we assume that the Lie algebra $\gg$ is simple. Let $\Phi\subset\gh^*$ 
denote the collection of the roots of $\gh$. Let $\{ e_{\beta}\}_{\beta\in\Phi}$ be 
the standard $\ad_{\gh}$- eigenbase of $E$, and let $h_{\beta}$ be the corresponding 
elements of $\gh$:

\begin{equation}[h_{\beta},e_{\pm\beta}]=\pm2e_{\pm\beta}, \ 
[e_{\beta},e_{-\beta}]=h_{\beta}, \ E=\oplus_{\beta\in\Phi}\rr e_{\beta}.
\label{sl2eh}\end{equation}
Fix a system $\Delta\subset\Phi$ of simple roots, i.e., an integer root basis 
(see Subsection 2.5). 
\begin{proposition} \label{prop3} For every $\alpha\in\Delta$ there exists a 
compact subset $K_{\alpha}\subset G$ such that each element of $\gg$ is 
conjugate to an element $v\in\gg$ of the following type by an element of $K_{\alpha}$:
\begin{equation} v=\sum_{\phi\in\Delta\setminus\alpha}c_{\phi}h_{\phi}+
\sum_{\beta\in\Phi}a_{\beta}e_{\beta}.\label{indv}\end{equation}
\end{proposition}

\begin{proof} Recall that the collection $\{ h_{\phi}\}_{\phi\in\Delta}$ 
is a basis of $\gh$, and its union with $\{ e_{\beta}\}_{\beta\in\Phi}$ is a 
basis of $\gg$. Fix a positive scalar product on $\gg$ for which 
the latter basis of $\gg$ is orthogonal. Let $\pi:\gg\to\gh$ denote 
the projection along $E$. We write each $v\in\gg$ as follows:
$$v=\pi v+\sum_{\beta\in\Phi}a_\beta(v)e_{\beta}.$$
 For every  $\alpha\in\Delta$, $v\in\gg\setminus0$ and $\delta>0$ set 
$$q_{\alpha}(v)=\frac{|a_{\alpha}(v)|}{|v-a_{\alpha}(v)e_{\alpha}|}\in[0,+\infty], 
\ V_{\delta}(\alpha)=\{ v\in\gg\setminus0 \ | \ q_{\alpha}(v)>\delta\}:$$
\begin{equation}
\gg\setminus(\cup_{\delta}V_{\delta}(\alpha))=H_{\alpha}=\{a_{\alpha}(v)=0\}=
\gh\oplus
\oplus_{\beta\in\Phi\setminus\alpha}\rr e_{\beta}.\label{vde}\end{equation}

\medskip

{\bf Claim 1.} {\it For every $\alpha\in\Delta$ there exist a $\delta_0>0$ and a 
finite subset $P_{\alpha}\subset G$ such that each element of $\gg\setminus0$ is 
 conjugate to an element of $V_{\delta_0}(\alpha)$ by an element of $P_{\alpha}$.}

\begin{proof} Fix an $\a\in\Delta$. The vector spaces $Ad_gH_{\alpha}$, $g\in G$, 
have zero intersection. Indeed, otherwise their intersection would be a nontrivial 
$Ad_G$- invariant subspace in $\gg$, - a contradiction to the simplicity of 
$\gg$. Therefore, there exists a finite collection of elements $g_1,\dots,
g_k\in G$ such that the images $Ad_{g_j}H_{\alpha}$ have zero intersection. 
Then there exists a $\delta_0>0$ such that 
$\cup_jAd_{g_j}V_{\delta_0}(\alpha)=\gg\setminus0$, 
by (\ref{vde}). This proves the claim for $P_{\alpha}=\{g_1^{-1},\dots,
g_k^{-1}\}$. 
\end{proof}

Fix an $\alpha\in\Delta$ and the corresponding $\delta_0$ and $P_{\alpha}$ from 
Claim 1. Let us fix a 
$u\in\gg\setminus0$ and show that it is conjugate to an expression (\ref{indv}) 
by an element from a compact subset in $G$. Applying conjugation by a 
$g\in P_{\alpha}$ one can achieve that $u\in V_{\delta_0}(\alpha)$ (Claim 1). 
Let $\pi u=\sum_{\phi\in\Delta}c_{\phi}h_{\phi}$, 
$s=\frac{c_{\alpha}(u)}{a_{\alpha}(u)}$: $|s|\leq\delta_0^{-1}$, by definition. 
Then the conjugate $Ad_{\exp(se_{-\alpha})}u$ has the type (\ref{indv}), by 
(\ref{sl2eh}). This proves Proposition \ref{prop3} for 
$K_{\alpha}=\exp([-\delta_0^{-1},\delta_0^{-1}]e_{-\alpha})P_{\alpha}$. 
\end{proof}

Now let us prove Lemma \ref{lemdec} by induction on the rank $r$ of $\gg$.

Induction base: $r=1$. Then Lemma \ref{lemdec} follows from Proposition 
\ref{prop3}. 

Induction step. Let us prove Lemma \ref{lemdec} for the given Lie algebra 
$\gg$, assuming it is already proved for the algebras of smaller ranks. Fix 
an $\alpha\in\Delta$ and the corresponding set  
$K_{\alpha}$ from Proposition \ref{prop3}. Let $\Phi_{\alpha}\subset\Phi$ denote the collection 
of those roots that are integer linear combinations of the roots from $\Delta\setminus
\alpha$. The linear subspace $\gt_{\alpha}\subset\gg$ generated by the vectors 
$h_{\phi}$ with $\phi\in\Delta\setminus\alpha$ and $e_{\beta}$ with $\beta\in
\Phi_{\alpha}$ is a semisimple Lie algebra, and hence, a direct sum of 
simple Lie algebras (\cite{br}, p.765, Lemma 4.2).  Let $T_{\alpha}\subset G$ be 
the connected (virtual) Lie subgroup corresponding to the subalgebra $\gt_{\alpha}$. 
Lemma \ref{lemdec} is  already 
proved for each simple summand  of $\gt_{\alpha}$ (the induction hypothesis), 
and hence, for the whole $\gt_{\alpha}$ and the group 
$Ad(T_{\alpha})\subset Aut(\gt_{\alpha})$. Thus, there exists a compact 
subset $K'\subset T_{\alpha}$ such that each element of $\gt_{\alpha}$ is 
conjugate to an element of $E\cap\gt_{\alpha}$ by an element of $K'$. 

Fix a $v\in\gg\setminus0$. Let us show that $v$ is conjugate to an element of 
$E$ by an element of some compact subset in $G$. 
Without loss of generality we assume that 
$v$ has the type (\ref{indv}), applying a conjugation by an element of 
$K_{\alpha}$. Let $u\in\gt_{\alpha}$ denote the $\gt_{\alpha}$- part of $v$: 
$$u=\sum_{\phi\in\Delta\setminus\alpha}c_{\phi}h_{\phi}+
\sum_{\beta\in\Phi_{\alpha}}a_{\beta}e_{\beta}.$$
The element $u$ is conjugate to an element of $E$ by 
an element $g\in K'\subset T_{\alpha}$. The difference $v-u$ is a linear 
combination of  vectors 
$e_{\beta}$ with $\beta\notin\Phi_{\alpha}$. The conjugation sends each 
of these vectors $e_{\beta}$ to a linear combination of 
the vectors $e_{\gamma}$ with $\gamma\equiv\beta$ modulo 
the integer linear combinations of the roots from $\Phi_{\alpha}$.  
Therefore, $\gamma\neq0$, $Ad_g(v-u)\in E$ and hence, $Ad_gv\in E$. 
This proves Lemma \ref{lemdec} for 
the compact subset $K=K'K_{\alpha}\subset G$. 
 \end{proof}
 
 Now let us prove Theorem \ref{thsk}. Fix 
a positive definite scalar product on $\gg$. Fix a regular 
element $\wt x\in\gh$. Take an arbitrary $z\in\gg$. Let $g\in K$ be an element such that 
$Ad_gz\in E$ (which exists by Lemma \ref{lemdec}). Set
$$x'=Ad_g^{-1}\wt x, \ y'=Ad_g^{-1}(\ad_{\wt x}|_E)^{-1}Ad_g z, \ 
L=\max_{g\in K\cup K^{-1}}
||Ad_g||,$$ 
\begin{equation}x=x'\sqrt{\frac{|y'|}{|x'|}}, \ y=y'\sqrt{\frac{|x'|}{|y'|}},  
c=\sqrt{L^3|\wt x|||(\ad_{\wt x}|_E)^{-1}||}\label{resc}\end{equation}
These $x$, $y$ and $c$ satisfy 
the Solovay-Kitaev inequality (\ref{insk}). Indeed, $[x,y]=z$ by construction and 
Remark \ref{remsk1} applied to $Ad_gz$ instead of $z$; 
 $|x|=|y|=\sqrt{|x'||y'|}\leq c\sqrt{|z|}$, by elementary inequalities. 
 Theorem \ref{thsk} is proved.
\end{proof}

\begin{proof} {\bf of Lemma \ref{lwsk}.} Let $\gg$ be a real semisimple Lie 
algebra. Fix a positive definite scalar product on $\gg$ and extend it to 
a Hermitian norm on $\gg_{\cc}$. The complexification $\gg_{\cc}$ satisfies 
the Solovay-Kitaev inequality (Theorem \ref{thsk}). Let $c>0$ denote 
the corresponding constant. Thus, for every $z\in\gg$ there exist 
$x,y\in\gg_{\cc}$ such that 
$$[x,y]=z, \ |x|=|y|\leq c\sqrt{|z|}; \ x=x_1+ix_2, \ y=y_1+iy_2.$$
Taking real part yields 
$$[x_1,y_1]+[-x_2,y_2]=z, \ |x_j|,|y_j|\leq|x|=|y|\leq c\sqrt{|z|}.$$
Now for every $j$ let us replace $x_j$ and $y_j$ by their rescalings 
with equal moduli, as in (\ref{resc}). The new $x_j$ and $y_j$ thus obtained satisfy the weak 
Solovay-Kitaev inequality with the above constant $c$. Lemma \ref{lwsk} is 
proved.
\end{proof}

\subsection{Approximations with weak Solovay-Kitaev property. Proof of Theorem 
\ref{gapprox1}}

First we give a proof of Theorem \ref{gapprox1}. It is similar to 
the proof of Theorem \ref{tsk} given in \cite{sk1, sk2, sk3}. The  proof of the 
boundedness of 
derivatives in Theorem \ref{tsk} (see Remark \ref{remsk}) will be briefly discussed 
at the end of the subsection. 

Let $G$ be a Lie group whose Lie algebra satisfies the weak Solovay-Kitaev 
inequality (\ref{wsk}), $c$ be the constant from (\ref{wsk}). The group $G$ 
is equipped with a left-invariant Riemannian metric. Let $\G$ be a 
free group with $M\geq2$ generators, $\rho:\G\to G$ be a dense representation. 
For every subset $\Omega$ in either $\G$, or $G$, set 
$$\Omega''=\{[x_1,y_1][x_2,y_2] \ | \ x_i,y_i\in\Omega\}.$$
 For the proof of Theorem \ref{gapprox1} it suffices to show 
that the group $G$ is $\var(x)$- approximable on the unit ball $D_1\subset G_0$ 
by $\rho(\G)$ with bounded derivatives, where the function $\var(x)$ and 
the length majorants $l_m=l_m(\rho,D_1)$ are the same
as in (\ref{enka}) and (\ref{lkgensk}). In (\ref{ome3}) we construct the corresponding 
subsets $\Omega_m=\Omega_{D_1,m}\subset\G$ by induction on 
$m$ so that each element of $\Omega_{m+1}$ is a product of an element of 
$\Omega_m$ and some (iterated) commutators of elements of 
$\cup_{k\leq m}\Omega_k$. To show that $\rho(\Omega_m)$ contains a net on 
$D_1$ with appropriate radius, we use 
the following lemma. It implies that if $\Omega\subset G_0$ is a $\delta$- net on 
$D_1$, then the product $\Omega\Omega''$ contains a $c''\delta^{\frac32}$- net on 
 $D_1$, $c''>0$ is a constant depending only on the metric of $G$, the constant 
 $\delta$ is arbitrary (small enough). 

\begin{lemma} \label{ldelta}
Let $G$ be a connected Lie group whose Lie algebra satisfies the weak Solovay-Kitaev 
inequality (\ref{wsk}). 
Then there exist constants $c',c''>0$ such that for every $\delta>0$ small enough and each  
$\delta$- net $\Omega$ on $D_{c'\sqrt{\delta}}\subset G$ the set
$\Omega''$ contains a $c''\delta^{\frac32}$- net on $D_{\delta}$. (The constant $c'$ may 
be chosen arbitrarily close to the corresponding constant $c$ from inequality (\ref{wsk}).) 
\end{lemma}

\begin{proof} Fix an arbitrary $c'>c$. For every small $\delta>0$ and every 
$z\in D_{\delta}$ set 
\begin{equation}v=\log z\in\gg:  \ |v|=dist(z,1)(1+o(1))\leq\delta(1+o(1)), \text{ as } 
\delta\to0,\label{v}\end{equation}
\begin{equation}
v=[u_1,v_1]+[u_2,v_2], \ u_i,v_i\in\gg, \ |u_i|=|v_i|\leq c\sqrt{|v|}\leq c\sqrt\delta(1+o(1))
<c'\sqrt\delta, 
\label{uivi}\end{equation} 
whenever $\delta$ is small enough (depending on $G$ and the choice of $c'$), 
by (\ref{wsk}) and (\ref{v}). Set
\begin{equation}\wt x_i=\exp{u_i}, \ \wt y_i=\exp{v_i}: \ \ \wt x_i,\wt y_i
\in D_{c'\sqrt{\delta}},\label{newmew}\end{equation} 
by (\ref{uivi}) and the left invariance of the metric. 

\medskip

{\bf Claim 2.} $dist(z,[\wt x_1,\wt y_1][\wt x_2,\wt y_2])=O(\delta^{\frac32})$, 
{\it as $\delta\to0$. The ``O'' is uniform in $z$, $\wt x_i$, $\wt y_i$.}

\begin{proof} 
We use the following  asymptotic relations between products and 
commutators in Lie group and sums and brackets in its Lie algebra:
\begin{equation}dist(\exp(a+b),\exp a\exp b)=O(|a||b|), \ 
\text{as} \ a,b\to0,\label{expab}\end{equation}
\begin{equation}dist(\exp([u,v]),[\exp u,\exp v])=O(|u|^3+|v|^3), \ \text{as} \ u,v\to0.
\label{expuv}\end{equation}
Formula (\ref{expab}) is obvious. 

\begin{proof} {\bf of (\ref{expuv}).} One has
$$[\exp u,\exp v]=\exp u\exp v\exp(-u)\exp(-v)=\exp(Ad_{\exp u}v)\exp(-v),$$
$$Ad_{\exp u}v=v+[u,v]+O(|u|^3+|v|^3).$$
Substituting the latter formula to the above right-hand side  
yields
$$[\exp u,\exp v]=\exp(([u,v]+O(|u|^3+|v|^3))+v)\exp(-v)=\exp([u,v]+O(|u|^3+|v|^3)),$$
by formula (\ref{expab}) applied to 
$a=[u,v]+O(|u|^3+|v|^3)$ and $b=v$. This proves (\ref{expuv}). 
\end{proof}

Set $a=[u_1,v_1]$, $b=[u_2,v_2]$:  $z=\exp(a+b)$; $|a|,|b|=O(\delta)$, 
by (\ref{v}) and (\ref{uivi}). One has 
\begin{equation}dist(z,\exp([u_1,v_1])\exp([u_2,v_2]))=O(\delta^2),\label{duv}
\end{equation}
by (\ref{expab}). Set $u=u_i$, $v=v_i$: $|u|=|v|=O(\sqrt\delta)$, 
by (\ref{uivi}). Then 
$$dist(\exp([u_i,v_i]),[\wt x_i,\wt y_i])=O(\delta^{\frac32}),$$
by (\ref{expuv}). 
Substituting the latter equality to (\ref{duv}) yields the claim. 
\end{proof}

{\bf Claim 3.} {\it Let $\delta>0$, $\Omega$ be a  $\delta$- net on $D_{c'\sqrt\delta}$, 
$z\in D_{\delta}$. Let $v$, $u_i$, $v_i$ and $\wt x_i$, $\wt y_i$ be the same, as in 
(\ref{uivi}) and (\ref{newmew}).   
Let $x_i, y_i\in\Omega$ be some $\delta$- approximants of $\wt x_i$ and 
$\wt y_i$ respectively. Then} 
\begin{equation}dist(z,[x_1,y_1][x_2,y_2])=O(\delta^{\frac32}), \text{ as } \delta\to0; 
\text{ the ``O'' is uniform in } z, \ \Omega, \ x_i, \ y_i.
\label{xy}\end{equation}

\begin{proof} The mapping $G\times G\to G$, $(g,h)\mapsto[g,h]$, 
restricted to $D_{\Delta}\times D_{\Delta}$  has derivative with norm bounded by  
$O(\Delta)$, as $\Delta\to0$.  Set $\Delta(\delta)=c'\sqrt\delta+\delta$. One has 
$x_i,\wt x_i,y_i,\wt y_i\in D_{\Delta(\delta)}$ and $dist(x_i,\wt x_i),dist(y_i,\wt y_i)<\delta$ 
for each $i=1,2$, by definition and (\ref{newmew}).  Therefore, 
$dist([x_i,y_i],[\wt x_i,\wt y_i])=O(\delta\Delta(\delta))=O(\delta^{\frac32})$. 
This together with Claim 2 implies (\ref{xy}). 
\end{proof}
 
There exists a constant $c''>0$ (depending only on $c'$ and the metric of $G$) 
such that 
the ``O'' in (\ref{xy}) is less than $c''\delta^{\frac32}$, whenever $\delta$ is small 
enough. Then $[x_1,y_1][x_2,y_2]\in\Omega''$ is a $c''\delta^{\frac32}$- 
approximant of $z$, by Claim 3. This proves Lemma \ref{ldelta}. 
\end{proof}

Let $c',c''$ be the constants from Lemma \ref{ldelta} (applied to the identity component 
$G_0$ of $G$). Choosing them big enough, without loss of generality 
everywhere below we can assume that 
\begin{equation}c', c''>1.\label{c'>1}\end{equation} 
Fix a $0<\delta<1$ small enough that  satisfies the statements 
 of Lemma \ref{ldelta} and the  inequality
 \begin{equation}c''\delta^{\frac32}<\delta.\label{c32}\end{equation}
Consider a sequence of positive numbers $\delta_m>0$ defined recursively as 
follows: 
 \begin{equation}\delta_1=\delta, \ \delta_{m+1}=c''\delta_m^{\frac32}.  
 \label{deltaj}\end{equation}
 Then $\delta_m\searrow0$ superexponentially fast, 
 by (\ref{c32}). Fix a finite subset $\Omega\subset\G$ so that 
 $$\rho(\Omega) \text{ is a } \delta- \text{ net on }D_1\subset G_0.$$
 We define  sequences of subsets $\Omega_m,\wt\Omega_m\subset\G$ by induction 
 in $m$ as follows:
 \def\omod{\Omega_1}
 \def\omm{\Omega_m}
 \begin{equation}\Omega_1=\Omega, \ \wt\Omega_1=\{ w\in\Omega_1 \ | \ 
 \rho(w)\in D_{2c'\sqrt{\delta_1}}\subset G_0\},\label{ome1}\end{equation}
 \begin{equation}\Omega_2=\omod\wt\Omega_1'', 
 \wt\Omega_2=\{ w\in\Omega_2 \ | \ \rho(w)\in D_{2c'\sqrt{\delta_2}}\subset G_0\},
 \label{ome2}\end{equation}
 \begin{equation}\Omega_{m+1}=\omm\wt\Omega_m'', \ 
 \wt\Omega_{m+1}=\wt\Omega_{m-1}''\wt\Omega_m'' \text{ for every } m\geq2.\label{ome3}
 \end{equation}
 We show that the sequence of collections $\Omega_m$, the set $K=D_1\subset G_0$ 
  and the numbers
 \begin{equation}l_m=9^{m-1}l_1, \ l_1=\max_{w\in\Omega}|w|,\label{delm}
 \end{equation}
 satisfy the statements of Definition \ref{grapp} (the $\var(x)$- 
 approximability on $D_1$ with bounded derivatives), whenever 
 $\delta$ is small enough. To do this, it suffices to show that 
 \begin{equation}|w|\leq l_m \text{ for every } w\in\omm \text{ and each } m\in\nn,
 \label{lmw}\end{equation}
 \begin{equation}\text{the subset } \rho(\omm)\subset G_0 \text{ contains a } 
 \delta_m- \text{ net on } D_1 \text{ for 
 all } m\in\nn,\label{7.14}\end{equation}
 \begin{equation} \text{there exists a } \hat c>0 
 \text{ (depending on } \delta \text{ and } \Omega) \text{ such that } 
 \delta_m<\var(\hat c l_m) \text{ for all } m\in\nn,\label{dm}\end{equation}
 \begin{equation} \text{the subset } \rho(\cup_m\omm)\subset G_0 \text{ is bounded,} 
 \label{bounds}\end{equation}
 there exists a neighborhood $V\subset G^M=R(\G,G)$ of $\rho$ where the 
mappings from (\ref{rw}): 
 \begin{equation}R_w:V\to G, \ R_w(\phi)=\phi(w),  \ w\in\cup_m\Omega_m, 
 \text{ have uniformly bounded derivatives.}\label{der}\end{equation}
 Statements (\ref{7.14}) and (\ref{dm}) imply that the set $\rho(\omm)$ contains an  
 $\var(\hat c l_m)$- net on $D_1$. This together with (\ref{lmw}), (\ref{bounds}) 
 and (\ref{der}) proves the $\var(x)$- approximability of $G$ on $D_1$ by 
 $\rho(\G)$ with bounded derivatives. This proves Theorem \ref{gapprox1}. 
 
 Statements (\ref{lmw})-(\ref{der}) are proved below. Statement (\ref{lmw}) will be 
 proved for arbitrary $\delta$, while statements (\ref{7.14})-(\ref{der}) will be proved 
 for every $\delta$ small enough (as it will be specified at the beginning of the proof 
 of each one of these statements). 
 
 \begin{proof} {\bf of (\ref{lmw}).} We prove (\ref{lmw}) and the next 
 auxiliary inequality  by induction on $m$: 
 \begin{equation} |w|\leq l_m \text{ for every } w\in\wt\Omega_m.\label{wtomm}\end{equation}
 
 Induction base. For $m=1$ inequalities (\ref{lmw}), (\ref{wtomm}) follow 
 from definition and (\ref{ome1}).
 
 \def\omj{\Omega_j}
 \def\omjl{\Omega_{j+1}}
 
 Induction step. Let (\ref{lmw}), (\ref{wtomm}) 
  hold  for all $m\leq j$. Let us prove them for $m=j+1$. One has 
  \begin{equation}|w|\leq8l_k \text{ for every } w\in
  \Omega_k''\cup\wt\Omega_k'', \ k\leq j, 
  \label{KJ}\end{equation}
 since $w=[x_1,y_1][x_2,y_2]$, $x_i,y_i\in\Omega_k\cup\wt\Omega_k$, 
 and $|x_i|,|y_i|\leq l_k$  (the induction hypothesis). 
 
 First let us  prove (\ref{lmw}). For every $w\in\Omega_{j+1}$ one has 
 $$w=w_1w_2, \ w_1\in\omj, \ w_2\in\wt\Omega_j'', \ |w_1|\leq l_j, \ |w_2|\leq 8l_j,$$ 
 by (\ref{ome2}), (\ref{ome3}), (\ref{KJ}) and the induction hypothesis of 
 (\ref{lmw}). Therefore, $|w|\leq 9l_j=l_{j+1}$.  Inequality (\ref{lmw}) for $m=j+1$ is 
 proved.
 
 Now let us prove inequality (\ref{wtomm}). Let $w\in\wt\Omega_{j+1}$. 
  If $j+1=2$, then $\wt\Omega_{j+1}\subset\omjl$, see (\ref{ome2}), 
 and inequality (\ref{lmw}) (already proved) implies (\ref{wtomm}). Let now $j+1\geq3$. 
 Then 
 $$w=w_1w_2, \ w_1\in\wt\Omega_{j-1}'', \ w_2\in\wt\Omega_j'', \ |w_1|\leq 8l_{j-1}, 
 \ |w_2|\leq 8l_j,$$
 by (\ref{KJ}). This together with the inequality
 $8l_{j-1}+8l_j=(1-\frac19)l_j+8l_j< 9l_j=l_{j+1}$
 proves (\ref{wtomm}) for $m=j+1$. The induction step is over. Inequalities (\ref{lmw}) 
 and (\ref{wtomm}) are proved.\end{proof}
 
 \begin{proof} {\bf of (\ref{7.14}).} One has 
 \begin{equation}2c'\sqrt\delta<1, \ \delta_m<c'\sqrt\delta_m, \ 
 \delta_{m-2}>c'\sqrt\delta_m \text{ for every } m\in\nn,\label{7.19}\end{equation}
 whenever $\delta$ is small enough. We prove (\ref{7.14}) for those $\delta$ that 
 satisfy the statements of Lemma \ref{ldelta} and (\ref{7.19}) by induction on $m$. 
 We use and prove simultaneously (by induction) that 
 \begin{equation} \text{the set } \rho(\wt\Omega_m'') \text{ contains a } 
 \delta_{m+1}- 
 \text{ net on } D_{\delta_m} \text{ for every } m.\label{7.20}\end{equation}
 
 \def\sqdl{\sqrt{\delta_1}}
 \def\sqdd{\sqrt{\delta_2}}
 \def\sqdm{\sqrt{\delta_m}}
 
 Induction base: m=1,2. The set $\rho(\Omega_1)$ is a $\delta_1$- net on $D_1$ 
 by definition. This proves (\ref{7.14}) for $m=1$. 
 Let us prove (\ref{7.20}) for $m=1$. The set $\rho(\wt\Omega_1)$ consists 
 of those elements of the above $\delta_1$- net that lie in the ball  
 $D_{2c'\sqrt{\delta_1}}$ (by definition). Therefore, $\rho(\wt\Omega_1)$ 
 contains a $\delta_1$- net on 
 $D_{2c'\sqrt{\delta_1}-\delta_1}\supset D_{c'\sqdl}$ 
 (the middle inequality in (\ref{7.19})). 
 Hence, the set $\rho(\wt\Omega_1'')$ contains a $\delta_2=
 c''\delta_1^{\frac32}$- net on $D_{\delta_1}$ (Lemma \ref{ldelta}). 
This proves (\ref{7.20}) for $m=1$. 
 
The set 
 $\rho(\Omega_2)=\rho(\omod)\rho(\wt\Omega_1'')$ contains a $\delta_2$- net on $D_1$, 
 since $\rho(\omod)$ is a $\delta_1$- net on $D_1$, $\rho(\wt\Omega_1'')$ 
 contains a $\delta_2$- net on $D_{\delta_1}$ (as was shown above) and by 
 Proposition \ref{pnet}. This proves (\ref{7.14}) for $m=2$. The set 
 $\rho(\wt\Omega_2'')$ contains a $\delta_3$- net on $D_{\delta_2}$, analogously to 
 the previous similar statement on $\rho(\wt\Omega_1'')$. This proves 
 statement (\ref{7.20}) for $m=2$. 
 
 Induction step: $m\geq3$. Assume that statements (\ref{7.14}), (\ref{7.20}) 
 hold for indices less than $m$. Let us prove them for the  given 
 $m$. The collection 
 $\rho(\omm)=\rho(\Omega_{m-1})\rho(\wt\Omega_{m-1}'')$
 contains a $\delta_m$- net on $D_1$, since $\rho(\Omega_{m-1})$ contains a 
 $\delta_{m-1}$- net on $D_1$ and $\rho(\wt\Omega_{m-1}'')$ contains a $\delta_m$- 
 net on $ D_{\delta_{m-1}}$ (the induction hypothesis). This 
 proves (\ref{7.14}). Similarly, the set $\rho(\wt\Omega_m)=
 \rho(\wt\Omega_{m-2}'')\rho(\wt\Omega_{m-1}'')$ contains a $\delta_m$- net on 
 $D_{\delta_{m-2}}$, since $\rho(\wt\Omega_{m-2}'')$ contains a $\delta_{m-1}$- net 
 on $D_{\delta_{m-2}}$ and $\rho(\wt\Omega_{m-1}'')$ contains a $\delta_m$- net on 
 $D_{\delta_{m-1}}$. 
 One has $\delta_{m-2}>c'\sqrt{\delta_m}$, by (\ref{7.19}). Hence, the set 
 $\rho(\wt\Omega_m)$ contains a $\delta_m$- net on $D_{c'\sqrt{\delta_m}}$.  
 Thus, the set 
 $\rho(\wt\Omega_m'')$ contains a $\delta_{m+1}=c''\delta_m^{\frac32}$- net on 
 $D_{\delta_m}$ (Lemma \ref{ldelta}). 
 This proves  (\ref{7.20}). The induction step is over. Statements (\ref{7.14}) and 
 (\ref{7.20}) are proved.
 \end{proof}

\begin{proof} {\bf of (\ref{dm}).} We prove (\ref{dm}) for every $\delta$ 
satisfying the inequality 
$(c'')^2\delta<1$. Set  $q=-\ln((c'')^2\delta)>0$. More precisely, we show that 
\begin{equation}\delta_m<\var(\hat c l_m)=e^{-(\hat c l_m)^{\kappa}} 
\text{ for every } 
m\in\nn, \text{ where } \kappa=\frac{\ln1.5}{\ln9}, \ \hat c=\frac{q^{\frac1{\kappa}}}{l_1}.
\label{7.15}\end{equation}
Indeed, the sequence $\ln\delta_m$ is the orbit of $\ln\delta_1=\ln\delta$ 
under the affine mapping $L(x)=\frac32x+\ln c''$, which has a fixed point  
$x_0=-2\ln c''<0$, see (\ref{c'>1}). Hence, 
\begin{equation}\ln\delta_m=x_0+(\frac32)^{m-1}(\ln\delta-x_0)=x_0-
(\frac32)^{m-1}q, 
\label{7.16}\end{equation}
since $\ln\delta-x_0=\ln\delta+2\ln c''=-q$ by definition. Recall that 
$x_0<0$. Therefore,  
\begin{equation}
\ln\delta_m<-(\frac32)^{m-1}q=-9^{\kappa(m-1)}q=-(\frac{l_m}{l_1})^{\kappa}q=
-(\hat cl_m)^{\kappa},
\label{dmq}\end{equation}
by definition and (\ref{7.16}). Exponentiating  
(\ref{dmq}) yields (\ref{7.15}) and proves (\ref{dm}). 
\end{proof}
\def\hgc{\hat G_{\cc}}
\def\hgco{\hat G_{\cc,0}}

\begin{proof} {\bf of (\ref{bounds}) and (\ref{der}).} We consider a  
neighborhood $\hgc\supset G$ of $G$ in its complexification equipped 
with the natural structure of (local) complex Lie group. We denote $\hgco$ 
the identity component of $\hgc$, which is a connected neighborhood of $G_0$. 
The left-invariant Riemannian metric on $G$ 
extends to a left-invariant Riemannian metric on $\hgc$. By definition, 
\begin{equation}\omm=\omod\Omega_{m-1}', \ 
\Omega_k'=\wt\Omega_1''\dots\wt\Omega_k'', \text{ for } m\geq2, \ k\geq1.
\label{prodd}\end{equation}
The collections $\omm$, $\omm'$ depend on $\delta$ and $\Omega$. In 
(\ref{tau}) we define a continuous function 
$\tau=\tau(\delta)$ in small $\delta$, $\tau(0)=0$. 
We show that for every small $\delta$ and any choice of $\Omega$ 
there exists a neighborhood $U\subset\hgc^M$ of 
$\rho\in\rgg=G^M\subset\hgc^M$ such that  
each mapping  $R_w:\rgg\to G$, $R_w(\phi)=\phi(w)$, 
see (\ref{rw}), with $w\in\cup_m\omm'$ extends to a holomorphic mapping
\begin{equation}R_w:U\to D_{\tau(\delta)}\subset\hgco.
\label{ude}\end{equation}
Then the Cauchy bound for derivative of holomorphic 
function implies  (\ref{der}). 

In what follows all the balls $D_q\subset\hgco$ under question are complex. 
Fix constants $\wt c,\sigma>0$ (which obviously exist) such that 
for every $q\leq\sigma$ the commutator mapping $G\times G\to G$: 
\begin{equation} (x,y)\mapsto[x,y] \text{ extends to   
a holomorphic mapping } D_q\times D_q\to D_{\wt c q^2}\Subset \hgco.
\label{p7.4}\end{equation}

\begin{proof} {\bf of (\ref{ude}).} Let $c'$ be the constant from (\ref{ome1}), 
$\wt c$, $\sigma$ be the constants from (\ref{p7.4}).  
Consider the sequence $\tau_m=\tau_m(\delta)>0$ defined recursively as follows 
and their sum $\tau$:
\begin{equation}\tau_1=\tau_2=4c'\sqrt\delta, \ \tau_j=4\wt c\max\{\tau_{j-2}^2,
\tau_{j-1}^2\} \text{ for } j\geq3; \ \ \tau=\tau(\delta)=\sum_{j=1}^{\infty}\tau_j.
\label{tau}
\end{equation}
Recall that $\hgc$ is an open subset of a complex 
manifold $\mathcal M$ with a holomorphic multiplication operation 
$\hgc\times\hgc\to \mathcal M$. It is equipped with a left-invariant Riemannian metric.  

There exists a $\delta_0>0$ such that 
$\tau_j(\delta_0)\searrow0$ superexponentially fast in $j\geq2$. 
The sum $\tau(\delta)$ is an increasing continuous function in $0\leq\delta\leq\delta_0$,  
since so are $\tau_j(\delta)$, $\tau(0)=0$. Fix this $\delta_0$ so that in addition, 
for every $0\leq\delta\leq\delta_0$ 
\begin{equation}\tau(\delta)<\min\{\sigma,\frac16dist(1,\partial\hgco)\}: \  
\ \overline D_{\tau_j}\Subset\overline D_{5\tau}\Subset\hgco.
\label{td}\end{equation}
Fix a $\delta\in(0,\delta_0]$. 
Take a neighborhood $U\subset\hgc^M$ of $\rho$ where  
the mappings $R_w$, $w\in \wt\Omega_1\cup\wt\Omega_2$, are holomorphic and 
such that 
\begin{equation}R_w(U)\subset D_{\tau_1}=D_{\tau_2}\subset\hgco 
\text{ for every } w\in 
\wt\Omega_1\cup\wt\Omega_2.\label{1,2}\end{equation}
This $U$ exists by (\ref{ome1}) and (\ref{ome2}). We prove 
(\ref{ude}) for this $U$. 

{\bf Claim 4.} {\it For every 
$m\in\nn$ and $w\in\wt\Omega_m\cup\wt\Omega_m''$ the mapping $R_w$ is holomorphic on $U$ and}  
\begin{equation}R_w(U)\subset D_{\tau_m}\subset\hgco \text{ for every } w\in
\wt\Omega_m,
 \label{nomd}\end{equation}
 \begin{equation}R_w(U)\subset 
D_{2\wt c\tau_{m}^2}\subset\hgco \text{ for every } w\in\wt\Omega_m''.
\label{nomdn}\end{equation}
\begin{proof} 
We prove (\ref{nomd}) and (\ref{nomdn}) simultaneously by induction using the 
 inequalities 
\begin{equation}4\wt c\tau_1^2<\tau_2, \ 2\wt c\tau_{m-2}^2+2\wt c\tau_{m-1}^2\leq\tau_{m}<\tau 
\text{ for every } m\geq3.
\label{tau2}\end{equation}
Indeed, for each $m\geq3$ one has 
$4\wt c\tau^2_{m-2}, 4\wt c\tau^2_{m-1}\leq\tau_m<\tau$, by (\ref{tau}). 
One has $4\wt c\tau_1^2\leq\tau_3<\tau_2$, by the latter inequality and decreasing. 
The two latter inequalities together imply 
(\ref{tau2}). 

Induction base for (\ref{nomd}). 
For $m=1,2$ statement (\ref{nomd}) follows from (\ref{1,2}). 

Induction base and step for (\ref{nomdn}). Let we have already proved 
(\ref{nomd}) for a given $m$. Let us prove (\ref{nomdn}) for the same $m$. 
Each $w\in\wt\Omega_m''$ is a commutator product $[x_1,y_1][x_2,y_2]$, 
$x_1,y_1,x_2,y_2\in\wt\Omega_{m}$. One has  
$R_{x_i}(U),R_{y_i}(U)\subset D_{\tau_m}$, by (\ref{nomd}). Therefore, 
each $R_{[x_i,y_i]}$ is a holomorphic mapping $U\to D_{\wt c\tau_{m}^2}$, 
by (\ref{p7.4}) and (\ref{td}). This implies (\ref{nomdn}). 

Induction step for (\ref{nomd}): $m\geq3$. Let 
statement (\ref{nomd}) hold  for smaller indices. Then statement (\ref{nomdn}) 
also holds for the same indices, as was shown above. By definition, 
$$\wt\Omega_m=\wt\Omega_{m-2}''\wt\Omega_{m-1}''.$$
This together with (\ref{nomdn}) (applied to the indices 
$m-2$ and $m-1$) and (\ref{tau2}) implies (\ref{nomd}) for the given  $m$. 
The induction step is over. Claim 4 is proved. 
 \end{proof}

\begin{corollary} \label{corclaim} 
For every $m\in\nn$ and $w\in\omm'$ the mapping $R_w$ is  
holomorphic on $U$ and  
\begin{equation}R_w(U)\subset D_{T_m}, \ T_m=\tau_1+\dots+\tau_m.\label{omt}
\end{equation}
\end{corollary}
\begin{proof} 
Induction base: $m=1$. Each  $w\in\omod'=\wt\Omega_1''$ defines a 
holomorphic mapping $R_w:U\to D_{2\wt c\tau_1^2}$, by (\ref{nomdn}). 
One has $2\wt c\tau_1^2<\tau_2=\tau_1=T_1$, by definition and 
(\ref{tau2}). This shows that $R_w(U)\subset D_{\tau_1}=D_{T_1}$ and proves 
(\ref{omt}).

Induction step: $m\geq2$. Let us prove (\ref{omt}) assuming it is proved for smaller 
indices. By definition, 
\begin{equation}\omm'=\Omega_{m-1}'\wt\Omega_m''.\label{qu}\end{equation}
Each $w\in\Omega_{m-1}'\cup\wt\Omega_{m}''$ defines a holomorphic mapping 
$R_w:U\to\hgco$,
\begin{equation}R_w(U)\subset D_{T_{m-1}} \text{ for every } w\in\Omega_{m-1}', 
\ R_w(U)\subset
D_{2\wt c\tau_{m}^2}\subset D_{\tau_m} \text{ for every }  w\in\wt\Omega_m'',\label{aq}\end{equation}
by the induction hypothesis, (\ref{nomdn}) and the inequality 
$2\wt c\tau_m^2\leq2\wt c\tau_{m-1}^2\leq\tau_m$ (the decreasing of the sequence 
$\tau_j$ for $j\geq2$ and (\ref{tau2})). One has 
$T_{m-1}+\tau_m=T_m$. This together with (\ref{qu}) and (\ref{aq}) proves the 
induction step and (\ref{omt}). 
\end{proof}

Corollary \ref{corclaim} immediately implies (\ref{ude}). 
\end{proof}

We take $\delta$ so small that the closed ball 
$\overline D_{\tau(\delta)}\subset\hgco$ is covered 
by a holomorphic chart of $\hgc$. The mappings (\ref{ude})  become 
uniformly bounded holomorphic mappings $U\to\cc^n$ in the latter chart. 
Take an arbitrary smaller (real) 
neighborhood $V\subset G^M$ of $\rho$ such that $\overline V\subset U$. The 
mappings (\ref{ude}) are uniformly bounded with derivatives on $V$ by the 
previous 
statement and the Cauchy estimate for derivatives of bounded holomorphic 
functions. 
 This together with (\ref{prodd}) implies (\ref{bounds}) and (\ref{der}). 
\end{proof}

 Statements (\ref{lmw})-(\ref{der}) imply Theorem \ref{gapprox1}. 
 The proof of Theorem \ref{gapprox1} is complete.

\begin{proof} {\bf of Theorem \ref{tsk} with the boundedness of derivatives (Remark 
\ref{remsk}).} 
To show that a Lie group with a Lie algebra satisfying the (strong) Solovay-Kitaev 
inequality is $\var'(x)$- approximable (see (\ref{var'})) with bounded derivatives, 
we replace the above $\Omega''$ by 
$$\Omega''=\{[x,y] \ | \ x,y\in\Omega\}.$$ 
Then an analogue of Lemma \ref{ldelta} holds for this $\Omega''$. 
Afterwards the proof repeats the previous one (of Theorem \ref{gapprox1}) with 
obvious changes.
\end{proof}

\section{Acknowledgments}

I am grateful to \'E.Ghys who had attracted my attention to the 
problem. I wish to thank him and also \'E.Breuillard, J.-F.Quint, B.Sevennec, 
G.Tomanov, E.B.Vinberg  for helpful conversations. I am grateful to the 
referees for very valuable remarks inspiring major improvements. 
A significant part of the proof of the Main Technical Lemma 
was obtained during my stay at RIMS (Kyoto). The paper was partly
written while I was visiting the Universities of Toronto and Stony Brook. 
I wish to thank these institutions for their hospitality and support.

\end{document}